\documentclass[10pt,reqno]{amsart}
\usepackage[a4paper, top=2cm, bottom=2cm, left=3cm, right=3cm, heightrounded, marginparwidth=3.5cm, marginparsep=.2cm, centering]{geometry}
\usepackage{appendix}
\usepackage{dsfont}
\usepackage{xcolor}
\usepackage[utf8]{inputenc}
\usepackage[colorlinks=true,hyperindex,pagebackref,linktocpage=true]{hyperref}
\usepackage{cleveref}
\usepackage{mathtools}
\usepackage{amsfonts}
\usepackage{amssymb}
\usepackage{tikz-cd}
\usepackage[T1]{fontenc}
\usepackage{pifont}
\usepackage{amsmath}
\usepackage{fontawesome}

\hypersetup{
	colorlinks,
	linkcolor={violet!60!black},
	citecolor={cyan!60!black},
	urlcolor={orange!60!black}
}
\usepackage{stmaryrd}
\usepackage{mathrsfs}  
\usepackage{varwidth}
\theoremstyle{plain}
\newtheorem{theorem}{Theorem}[section]
\newtheorem{theorem*}{Theorem}
\newtheorem{remark*}{Remark}
\newtheorem{conjecture*}{Conjecture}
\newtheorem{proposition}[theorem]{Proposition}

\newtheorem{lemma}[theorem]{Lemma}
\newtheorem{corollary}[theorem]{Corollary}
\newtheorem{conjecture}[theorem]{Conjecture}

\theoremstyle{definition}
\newtheorem{definition}[theorem]{Definition}
\newtheorem{example}[theorem]{Example}
\newtheorem{remark}[theorem]{Remark}

\newcommand{\nc}{\newcommand}
\nc{\on}{\operatorname}

\nc{\Q}{\mathbb{Q}}
\nc{\Z}{\mathbb{Z}}
\nc{\cl}{\mathrm{cl}}

\nc{\fraka}{{\mathfrak a}} \nc{\bba}{{\mathbf a}}
\nc{\frakb}{{\mathfrak b}}
\nc{\frakc}{{\mathfrak c}}
\nc{\frakd}{{\mathfrak d}}
\nc{\frake}{{\mathfrak e}}
\nc{\frakf}{{\mathfrak f}}
\nc{\frakg}{{\mathfrak g}}
\nc{\frakh}{{\mathfrak h}}
\nc{\fraki}{{\mathfrak i}}
\nc{\frakj}{{\mathfrak j}}
\nc{\frakk}{{\mathfrak k}}
\nc{\frakl}{{\mathfrak l}}
\nc{\frakm}{{\mathfrak m}}
\nc{\frakn}{{\mathfrak n}}
\nc{\frako}{{\mathfrak o}}
\nc{\frakp}{{\mathfrak p}}
\nc{\frakq}{{\mathfrak q}}
\nc{\frakr}{{\mathfrak r}}
\nc{\fraks}{{\mathfrak s}}
\nc{\frakt}{{\mathfrak t}}
\nc{\fraku}{{\mathfrak u}}
\nc{\frakv}{{\mathfrak v}}
\nc{\frakw}{{\mathfrak w}}
\nc{\frakx}{{\mathfrak x}}
\nc{\fraky}{{\mathfrak y}}
\nc{\frakz}{{\mathfrak z}}
\nc{\frakA}{{\mathfrak A}}
\nc{\frakB}{{\mathfrak B}}
\nc{\frakC}{{\mathfrak C}}
\nc{\frakD}{{\mathfrak D}}
\nc{\frakE}{{\mathfrak E}}
\nc{\frakF}{{\mathfrak F}}
\nc{\frakG}{{\mathfrak G}}
\nc{\frakH}{{\mathfrak H}}
\nc{\frakI}{{\mathfrak I}}
\nc{\frakJ}{{\mathfrak J}}
\nc{\frakK}{{\mathfrak K}}
\nc{\frakL}{{\mathfrak L}}
\nc{\frakM}{{\mathfrak M}}
\nc{\frakN}{{\mathfrak N}}
\nc{\frakO}{{\mathfrak O}}
\nc{\frakP}{{\mathfrak P}}
\nc{\frakQ}{{\mathfrak Q}}
\nc{\frakR}{{\mathfrak R}}
\nc{\frakS}{{\mathfrak S}}
\nc{\frakT}{{\mathfrak T}}
\nc{\frakU}{{\mathfrak U}}
\nc{\frakV}{{\mathfrak V}}
\nc{\frakW}{{\mathfrak W}}
\nc{\frakX}{{\mathfrak X}}
\nc{\frakY}{{\mathfrak Y}}
\nc{\frakZ}{{\mathfrak Z}}
\nc{\bbA}{{\mathbb A}}
\nc{\bbB}{{\mathbb B}}
\nc{\bbC}{{\mathbb C}}
\nc{\bbD}{{\mathbb D}}
\nc{\bbE}{{\mathbb E}}
\nc{\bbF}{{\mathbb F}} \nc{\bbf}{{\mathbf f}}
\nc{\bbG}{{\mathbb G}}
\nc{\bbH}{{\mathbb H}}
\nc{\bbI}{{\mathbb I}}
\nc{\bbJ}{{\mathbb J}}
\nc{\bbK}{{\mathbb K}}
\nc{\bbL}{{\mathbb L}}
\nc{\bbM}{{\mathbb M}}
\nc{\bbN}{{\mathbb N}}
\nc{\bbO}{{\mathbb O}}
\nc{\bbP}{{\mathbb P}}
\nc{\bbQ}{{\mathbb Q}}
\nc{\bbR}{{\mathbb R}}
\nc{\bbS}{{\mathbb S}}
\nc{\bbT}{{\mathbb T}}
\nc{\bbU}{{\mathbb U}}
\nc{\bbV}{{\mathbb V}}
\nc{\bbW}{{\mathbb W}}
\nc{\bbX}{{\mathbb X}}
\nc{\bbY}{{\mathbb Y}}
\nc{\bbZ}{{\mathbb Z}}
\nc{\calA}{{\mathcal A}}
\nc{\calB}{{\mathcal B}}
\nc{\calC}{{\mathcal C}}
\nc{\calD}{{\mathcal D}}
\nc{\calE}{{\mathcal E}}
\nc{\calF}{{\mathcal F}}
\nc{\calG}{{\mathcal G}}
\nc{\calH}{{\mathcal H}}
\nc{\calI}{{\mathcal I}}
\nc{\calJ}{{\mathcal J}}
\nc{\calK}{{\mathcal K}}
\nc{\calL}{{\mathcal L}}
\nc{\calM}{{\mathcal M}}
\nc{\calN}{{\mathcal N}}
\nc{\calO}{{\mathcal O}}
\nc{\calP}{{\mathcal P}}
\nc{\calQ}{{\mathcal Q}}
\nc{\calR}{{\mathcal R}}
\nc{\calS}{{\mathcal S}}
\nc{\calT}{{\mathcal T}}
\nc{\calU}{{\mathcal U}}
\nc{\calV}{{\mathcal V}}
\nc{\calW}{{\mathcal W}}
\nc{\calX}{{\mathcal X}}
\nc{\calY}{{\mathcal Y}}
\nc{\calZ}{{\mathcal Z}}

\nc{\scrA}{{\mathscr A}}
\nc{\scrB}{{\mathscr B}}
\nc{\scrC}{{\mathscr C}}
\nc{\scrD}{{\mathscr D}}
\nc{\scrE}{{\mathscr E}}
\nc{\scrF}{{\mathscr F}}
\nc{\scrG}{{\mathscr G}}
\nc{\scrH}{{\mathscr H}}
\nc{\scrI}{{\mathscr J}}
\nc{\scrJ}{{\mathscr I}}
\nc{\scrK}{{\mathscr K}}
\nc{\scrL}{{\mathscr L}}
\nc{\scrM}{{\mathscr M}}
\nc{\scrN}{{\mathscr N}}
\nc{\scrO}{{\mathscr O}}
\nc{\scrP}{{\mathscr P}}
\nc{\scrQ}{{\mathscr Q}}
\nc{\scrR}{{\mathscr R}}
\nc{\D}{{\on{D}}}
\nc{\Perv}{{\on{Perv}}}

\nc{\bnu}{{\bar{ \nu}}}
\nc{\olO}{\bar{\calO}}

\nc{\al}{{\alpha}} 
\nc{\be}{{\beta}}
\nc{\ga}{{\gamma}} \nc{\Ga}{{\Gamma}}
\nc{\hGa}{\hat{\Gamma}}
\nc{\ve}{{\varepsilon}} 
\nc{\la}{{\lambda}} \nc{\La}{{\Lambda}}
\nc{\om}{\omega} \nc{\Om}{\Omega} 
\nc{\sig}{{\sigma}} \nc{\Sig}{{\Sigma}}
\nc{\dR}{{\mathrm{dR}}}
\nc{\Perf}{{\mathrm{Perf}}}
\nc{\PSch}{{\mathrm{PSch}}}
\nc{\GrW}{{\calF\ell_{\calG,\bbW}}}
\newcommand{\Tf}[2]{{\widehat{#1}_{/#2}}}
\newcommand{\Lgrm}{{\widehat{\bbW\calM^{\leq \mu}_{\calD,O_{\breve{E}_k}}}}}
\newcommand{\Lgr}{{\widehat{\bbW\Gr_\calD}}}
\newcommand{\Lshtm}{{\widehat{\bbW\Sht^{\leq \mu}_{\calD,O_{\breve{E}_k}}}}}
\newcommand{\Lsht}{{\widehat{\bbW\Sht_{\calD}}}}
\newcommand{\calGd}{{\widehat{L^+_{\bbW}\calG}}_\calD}
\newcommand{\calGphid}{{\widehat{L^+_{\bbW}\calG}}_{\varphi^*\calD}}

\nc{\Shtdos}{{\Sht_\calG(b)}}
\nc{\Auxshte}[1]{{\bbW\Sht_\calG(#1)}}
\newcommand{\Shtm}[1]{{\Sht}_{\calG,{#1}}^{\leq \mu}(b)}
\nc{\adlvdos}{X^{\leq \mu}_\calG(b)}

\nc{\Auxshtm}{{\bbW\Sht^{\leq\mu}_{\calG,O_{\breve{E}}}(b)}}

\nc{\Auxshten}[1]{{\bbW\Sht_{\calG_{#1}}(b_{#1})}}
\nc{\Auxsht}{{\bbW\Sht_\calG(b)}}
\newcommand{\Grm}[1]{{\calM^{\leq \mu}_{\calG,#1}}}
\newcommand{\GrWmfe}{{\calA_{\calG,\mu}}}

\nc{\Witt}{{\mathrm{Witt}}}

\newcommand{\Yint}[2]{{\calY}^{#1}_{#2}}

\newcommand{\Hub}[1]{{(#1,#1^+)}}
\newcommand{\Yinf}[1]{{\calY^{#1}_{[0,\infty)}}}

\nc{\Perfd}{{\mathrm{Perfd}}}
\nc{\Perff}{\mathrm{Perf}_{\bar{\bbF}_p}}
\nc{\Frob}{\mathrm{Frob}}
\nc{\Div}{{\mathrm{Div}}}
\nc{\Sets}{{\mathrm{Sets}}}
\nc{\Sht}{{\mathrm{Sht}}}
\nc{\Gm}{{\mathbb{G}_m}}
\nc{\colim}{{\on{colim}}}
\nc{\et}{\mathrm{\acute{e}t}}

\nc{\dom}{\mathrm{\on{dom}}}
\nc{\Rep}{\mathrm{\on{Rep}}}
\nc{\Isoc}{\mathrm{\on{Isoc}}}

\makeatletter
\DeclareFontEncoding{LS1}{}{}
\DeclareFontSubstitution{LS1}{stix}{m}{n}
\DeclareMathAlphabet{\rhomalpha}{LS1}{stixscr}{m}{n}
\makeatother

\nc{\Spa}{\on{{Spa}}}
\nc{\Spd}{\on{{Spd}}}
\nc{\tnb}{\psi_{\rm tame}}
\nc{\oM}{\overline{{M}}}
\nc{\op}{{\on{op}}}
\nc{\ad}{{\on{ad}}}
\nc{\alg}{{\on{alg}}}
\nc{\Ad}{{\on{Ad}}}
\nc{\Adm}{{\on{Adm}}} \nc{\aff}{{\on{af}}}
\nc{\Aut}{{\on{Aut}}}
\nc{\Bun}{{\on{Bun}}}
\nc{\cha}{{\on{char}}}
\nc{\der}{{\on{der}}}
\nc{\Nm}{{\on{Nm}}}
\nc{\ab}{{\on{ab}}}
\nc{\scn}{{\on{sc}}}
\nc{\Der}{{\on{Der}}}
\nc{\diag}{{\on{diag}}}
\nc{\End}{{\on{End}}}
\nc{\Fl}{{\calF\!\ell}}
\nc{\Tr}{{\on{Transp}}}
\nc{\TR}{{\calT\!\calR}}
\nc{\Gal}{{\on{Gal}}}
\nc{\Gr}{{\on{Gr}}}
\nc{\Hk}{{\on{Hk}}}
\nc{\rH}{{\on{H}}}
\nc{\Hom}{{\on{Hom}}}
\nc{\IC}{{\on{IC}}}
\nc{\id}{{\on{id}}}
\nc{\Id}{{\on{Id}}}
\nc{\ind}{{\on{ind}}}
\nc{\Ind}{{\on{Ind}}}
\nc{\Lie}{{\on{Lie}}}
\nc{\Pic}{{\on{Pic}}}
\nc{\pr}{{\on{pr}}}
\nc{\Res}{{\on{Res}}}
\nc{\res}{{\on{res}}} \nc{\Sat}{{\on{Sat}}}
\nc{\spc}{{\on{sc}}}
\nc{\drv}{{\on{der}}}
\nc{\sgn}{{\on{sgn}}}
\nc{\Spec}{{\on{Spec}\,}}
\nc{\Spf}{\on{Spf}} 
\nc{\Sph}{\on{Sph}}
\nc{\St}{{\on{St}}}
\nc{\tr}{{\on{tr}}}
\nc{\Mod}{{\mathrm{-Mod}}}
\nc{\Hilb}{{\on{Hilb}}} 
\nc{\Ext}{{\on{Ext}}} 
\nc{\vs}{{\on{Vec}}}
\nc{\ev}{{\on{ev}}}
\nc{\nO}{{\breve{\calO}}}
\nc{\tS}{{\tilde{S}}}
\nc{\spe}{{\on{sp}}}
\nc{\loc}{{\on{loc}}}
\nc{\pre}{{\on{pre}}}

\nc{\dimt}{{\on{dim.trg}}}

\nc{\co}{\colon}
\nc{\dia}{{\diamond}}

\nc{\nscrR}{{\mathscr{R}^{\on{nr}}}}

\nc{\GL}{{\on{GL}}}

\nc{\Gl}{\on{Gl}} 
\nc{\GSp}{{\on{GSp}}}
\nc{\gl}{{\frakg\frakl}}
\nc{\SL}{{\on{SL}}} 
\nc{\SU}{{\on{SU}}} 
\nc{\SO}{{\on{SO}}}
\nc{\PGL}{{\on{PGL}}}

\nc{\Conv}{{\on{Conv}}}
\nc{\Dom}{{\on{Dom}}}
\nc{\red}{{\on{red}}}
\nc{\Red}{{\on{Red}}}
\nc{\act}{{\on{act}}}
\nc{\nr}{{\on{nr}}}
\nc{\ctf}{{\on{ctf}}}

\nc{\str}{{\on{-}}} 
\nc{\os}{{\bar{s}}}
\nc{\oeta}{{\bar{\eta}}}

\nc{\hookto}{\hookrightarrow}
\nc{\longto}{\longrightarrow}
\nc{\leftto}{\leftarrow}
\nc{\onto}{\twoheadrightarrow}
\nc{\lonto}{\twoheadleftarrow}

\nc{\pot}[1]{ [\hspace{-0,5mm}[ {#1} ]\hspace{-0,5mm}] }
\nc{\rpot}[1]{ (\hspace{-0,7mm}( {#1} )\hspace{-0,7mm}) }
\nc{\smallpot}{{ <\hspace{-1,0mm}<}}

\numberwithin{equation}{section}

\begin{document}
	
\title[Geometric components of moduli of $p$-adic shtukas.]{On the geometric connected components of moduli spaces of $p$-adic shtukas and local Shimura varieties.}
	
	\author[I. Gleason]{Ian Gleason}
	\address{ Department of Mathematics, Faculty of Science, National University of Singapore.  Block S17, 10 Lower Kent Ridge Road, Singapore 119076}
	\email{ianandreigf@nus.edu.sg}

\begin{abstract}
	We study topological properties of moduli spaces of $p$-adic shtukas and local Shimura varieties. 
On one hand, we construct and study the specialization map for moduli spaces of $p$-adic shtukas at parahoric level whose target is an affine Deligne--Lusztig variety. 
On the other hand, given a $p$-adic shtuka datum $(G,b,\mu)$, with $G$ unramified over $\bbQ_p$ and such that $(b,\mu)$ is HN-irreducible, we determine the set of geometric connected components of infinite level moduli spaces of $p$-adic shtukas.
In other words, we understand $\pi_0(\Sht_{G,b,\mu,\infty}\times \Spd \bbC_p)$ with its right $G(\bbQ_p)\times G_b(\bbQ_p)\times W_E$-action. 
	As a corollary, we prove new cases of a conjecture of Rapoport and Viehmann. 
\end{abstract}

\maketitle
\begin{center}
    {\large \bfseries
{Sur les composantes connexes g\'eom\'etriques
des espace de modules de chtoucas $p$-adiques et vari\'et\'es de
Shimura locales.}}
  \end{center}
 \small Nous étudions les propriétés topologiques des espaces de modules de chtouca $p$-adiques et des variétés de Shimura locales. 
D'une part, nous construisons et étudions l'application de spécialisation pour les espaces de modules de shtukas $p$-adiques au niveau parahorique dont la but du morphisme est une variété affine de Deligne–Lusztig. 
D'autre part, étant donné une donn\'ee de chtouca $p$-adique $(G,b,\mu)$, avec $G$ non ramifié sur $\bbQ_p$ et tel que $(b,\mu)$ soit HN-irréductible, nous déterminons l'ensemble des composantes connexes géométriques des espaces de modules de niveau infini de shtukas $p$-adiques. 
En d'autres termes, nous comprenons $\pi_0(\Sht_{G,b,\mu,\infty}\times \Spd \bbC_p)$ avec son action de $G(\bbQ_p)\times G_b(\bbQ_p)\times W_E$ \`a droite. 
En corollaire, nous prouvons de nouveaux cas d'une conjecture de Rapoport et Viehmann.
\normalsize

\tableofcontents

\subsection*{Introduction.}
In \cite{RV14}, Rapoport and Viehmann propose that there should be a theory of $p$-adic local Shimura varieties. They conjecture the existence of towers of rigid-analytic spaces whose cohomology groups “understand” the local Langlands correspondence for general $p$-adic reductive groups. In this way, these towers of rigid-analytic varieties would “interact” with the local Langlands correspondence in a manner similar to how Shimura varieties “interact” with the global Langlands correspondence. Moreover, they conjecture many properties and compatibilities that these towers should satisfy (see \cite[$\mathsection$ 5]{RV14}). 

In the last decade, the theory of local Shimura varieties has gone through a drastic transformation with Scholze's introduction of perfectoid spaces and the theory of diamonds \cite{Sch17}. 
In \cite{SW20}, Scholze and Weinstein construct the sought-after towers of rigid-analytic spaces and generalize them to what are now known as moduli spaces of $p$-adic shtukas (or mixed characteristic local shtukas) \cite[\S 23.1]{SW20}. 
Moreover, since then, many of the expected properties and compatibilities for local Shimura varieties have been verified and generalized to moduli spaces of $p$-adic shtukas. The study of the geometry and cohomology of local Shimura varieties and moduli spaces of $p$-adic shtukas is still a very active area of research due to its connection to the local Langlands correspondence. 
One of the main aims of this article is to study the set of geometric connected components of moduli spaces of $p$-adic shtukas (see \Cref{infinitelevelmoduli})
\[\pi_0(\Sht^{\on{geo}}_{G,b,\mu,\infty})\]
attached to a $p$-adic shtuka datum $(G,b,\mu)$ (as in \Cref{some-section-defining-shtukas}), and to describe the right action of the group $G({\bbQ_p})\times G_b({\bbQ_p})\times W_E$ on $\pi_0(\Sht^{\on{geo}}_{G,b,\mu,\infty})$. 
This set controls the first cohomology group of moduli spaces of $p$-adic shtukas. 
The upshot is that connected components of moduli spaces of $p$-adic shtukas are completely described by local class field theory (see \Cref{thm2:mainTheorem} for a precise statement). 
As a consequence of our results, we settle \cite[Conjecture 4.30]{RV14} in the case of unramified groups (see \Cref{thm2:mainTheorem}).  \\

Let us recall the formalism of local Shimura varieties and moduli spaces of $p$-adic shtukas. 
Let $\bbC_p\supseteq \bbQ_p$ denote a completed algebraic closure of $\bbQ_p$ and let $\breve{\bbQ}_p\subseteq \bbC_p$ denote the completion of the maximal unramified extension of $\bbQ_p$ in $\bbC_p$.
A local $p$-adic shtuka datum over ${\bbQ_p}$ is a triple $(G,b,\mu)$ where $G$ is a reductive group over ${\bbQ_p}$, $\mu$ is a conjugacy class of geometric cocharacters $\mu:\bbG_m\to G_{\bar{\bbQ}_p}$, and $b$ is an element of Kottwitz' set $B(G,\mu)$ \cite[$\mathsection$ 6]{Kot97}. 
Whenever $\mu$ is minuscule, we say that $(G,b,\mu)$ is a local Shimura datum (see \cite[Definition 5.1]{RV14}). 
We let $E/{\bbQ_p}$ denote the reflex field of $\mu$ and we let $\breve{E}=E\cdot \breve{\bbQ}_p$ be the compositum inside $\bbC_p$. 
Associated to $(G,b,\mu)$, there is a tower of diamonds over $\on{Spd}\breve{E}$, denoted $(\Sht_{G,b,\mu,K})_K$, where $K\subseteq G({\bbQ_p})$ ranges over compact open subgroups of $G({\bbQ_p})$ \cite[$\mathsection$ 23.3]{SW20}. 
Moreover, whenever $\mu$ is minuscule, $\Sht_{G,b,\mu,K}$ is represented by the diamond associated to a unique smooth rigid-analytic space $\bbM_K$ over $\breve{E}$. 
The tower $(\bbM_K)_K$ is the local Shimura variety \cite[Definition 24.1.3]{SW20}. 

Associated to $b\in B(G,\mu)$ there is a reductive group $G_b$ over ${\bbQ_p}$\footnote{In the literature, the group that we denote $G_b$ is often denoted by $J_b$.} \S \ref{subsection-Gb}. After base change to a completed algebraic closure, each individual space $\Sht_{G,b,\mu,K}^{\on{geo}}:=\Sht_{G,b,\mu,K}\times \Spd \bbC_p$ comes equipped with continuous and commuting right actions by $G_b({\bbQ_p})$ and the Weil group $W_E$. Moreover, the tower receives a right action by the group $G({\bbQ_p})$ by using correspondences. When we take the limit as $K\subseteq G(\bbQ_p)$ shrinks we obtain the space at infinite level, which we denote $\Sht^{\on{geo}}_{G,b,\mu,\infty}$. 
Overall, this space comes equipped with a right action by the group $G({\bbQ_p})\times G_b({\bbQ_p})\times W_E$ \S \ref{definingsection}. 

The group $G({\bbQ_p})\times G_b({\bbQ_p})\times W_E$ acts continuously on $\pi_0(\Sht^{\on{geo}}_{G,b,\mu,\infty})$ and one of our main theorems (see \Cref{thm2:mainTheorem}) describes explicitly this action under two technical assumptions on the triple $(G,b,\mu)$. 
The first assumption is that $G$ is an unramified reductive group i.e. $G$ is a quasi-split connected reductive group over $\bbQ_p$ whose base change to some unramified extension $\bbQ_{p^s}$ becomes split (e.g. all split groups are unramified). 
The second assumption is that the pair $(b,\mu)$ is HN-irreducible (Hodge--Newton irreducible).  
For $G=\on{GL}_n$ this condition asks in rough terms that the Hodge polygon determined by $\mu$ and the Newton polygon determined by $b$ do not meet at a ``breaking point'' (see \Cref{defi:BGmu} for the precise definition). 
For unramified groups, our result is optimal in a sense which we discuss later in this introduction (see also \Cref{whereHn-comesin}). 
It is also likely that the condition that $G$ is unramified can be removed (see \Cref{Remark-unramifiedness-rmove}).

Before stating the main theorem of $\mathsection$\ref{thiisghesecondsection} we set some notation. 
Let $(G,b,\mu)$ be local $p$-adic shtuka datum with $G$ an unramified reductive group over ${\bbQ_p}$. 
Let $G^{\on{der}}$ denote the derived subgroup of $G$, $G^{\on{sc}}$ denote the simply connected cover of $G^{\on{der}}$, consider the image of $G^{\on{sc}}({\bbQ_p})$ in $G({\bbQ_p})$ and let ${G(\bbQ_p)_{\circ}}=G({\bbQ_p})/\on{Im}(G^\scn(\bbQ_p))$. 
The group ${G(\bbQ_p)_{\circ}}$ is a locally profinite topological group and it is the maximal abelian quotient of $G({\bbQ_p})$ when this latter is considered as an abstract group (see \Cref{drevidesubgroupsbusiness}). 
Let ${\on{Art}}_{E}:W_{E}\to E^\times$ be the reciprocity character from local class field theory. 
In $\mathsection$\ref{extendedversionsofmpassection} we associate to $\mu$ and to $b$ continuous maps of topological groups ${\on{Nm}}_{\mu}^{E,\circ}:E^\times \to {G(\bbQ_p)_{\circ}}$ and $\on{det}_b^\circ:G_b({\bbQ_p})\to {G(\bbQ_p)_{\circ}}$ respectively. 
\begin{theorem*}[\Cref{thm2:mainTheorem}]
	\label{theoremintheintro}
	Let $(G,b,\mu)$ be a $p$-adic shtuka datum such that $G$ is an unramified reductive group over $\bbQ_p$ and such that the pair $(b,\mu)$ is HN-irreducible. 
	Let $E$ denote the reflex field of $\mu$.
Then the following hold.
	\begin{enumerate}
		\item The right $G(\bbQ_p)$-action on $\pi_0(\Sht^{\on{geo}}_{G,b,\mu,\infty})$ is trivial on $\on{Im}(G^\scn(\bbQ_p))$ and the corresponding $G^{\circ}$-action is simply-transitive. 
		\item If $s\in  \pi_0(\Sht^{\on{geo}}_{G,b,\mu,\infty})$ and $j\in G_b(\bbQ_p)$ then 
			\[s\star_{G_b} j=s\star_{{G(\bbQ_p)_{\circ}}} {\on{det}}_b^\circ(j^{-1}).\]		
		\item If $s\in  \pi_0(\Sht^{\on{geo}}_{G,b,\mu,\infty})$ and $\gamma \in W_{{E}}$ then 
			\[s\star_{{W_{E}}} \gamma =s\star_{{{G(\bbQ_p)_{\circ}}}} [{\on{Nm}}^{E,\circ}_{\mu}\circ {\on{Art}}_{{E}}(\gamma)].\] 
	\end{enumerate}
\end{theorem*}

\begin{remark*}
	In the particular case in which $G^\der=G^\scn$ the above theorem is established by constructing an equivariant isomorphism
	\[\pi_0(\Sht^{\on{geo}}_{G,b,\mu,\infty})\simeq \pi_0(\Sht^{\on{geo}}_{G^\ab,b^\ab,\mu^\ab,\infty}) \]
	where $(G^\ab,b^\ab,\mu^\ab)$ is the $p$-adic shtuka datum attached to the maximal abelian quotient of $G$.
	This case settles \cite[Conjecture 4.30]{RV14} of Rapoport and Viehmann in the case in which the group is unramified (see \Cref{corollaryRVconj}). 
\end{remark*}

\begin{remark*}
	Since moduli spaces of $p$-adic shtukas (as with most moduli spaces) do not have an explicit presentation, saying concrete things about their geometry is usually hard.
	This is to be expected since the reason we study moduli spaces of $p$-adic shtukas is to get a better understanding of the local Langlands correspondence, which is itself a very deep statement. 
	Although we do not discuss explicit applications of our theorem in this article, we believe that our result is not only hard but also powerful. 
	To convince the reader about this, we recall that the study of connected components of affine Deligne--Lusztig varieties \cite{CKV15} was one of the key ingredient in Kisin's work on integral models of Shimura varieties \cite{Kisin-mod-p-points}.
	As we clarify in $\mathsection$\ref{sectionspecialization}, the geometry of moduli spaces of $p$-adic shtukas is intimately related to the geometry of affine Deligne--Lusztig varieties. 
This relation and the methods developed in this article are exploited in our collaboration with Lim and Xu \cite{gleason2023connectedcomponentsaffinedelignelusztig} to give a new method of computing the connected components of affine Deligne--Lusztig varieties. 
\end{remark*}

\begin{remark*}
	\label{whereHn-comesin}
	The above theorem is optimal for unramified groups in the following sense. 
	One can prove that the action of $G(\bbQ_p)$ on $\pi_0(\Sht^{\on{geo}}_{G,b,\mu,\infty})$ only factors through ${G(\bbQ_p)_{\circ}}$ when $(b,\mu)$ is HN-irreducible. 
	Moreover, we expect that combining the methods of \cite{gaisin2022nonsemistablelociheckestacks} and \cite{MR4331441} with the methods in the present article one can express the general formula for the $G(\bbQ_p)\times G_b(\bbQ_p)\times W_E$-action on $\pi_0(\Sht^{\on{geo}}_{G,b,\mu,\infty})$ in terms of a parabolic induction of $\pi_0(\Sht^{\on{geo}}_{M,b_M,\mu_M,\infty})$ for HN-irreducible data $(M,b_M,\mu_M)$ associated to Levi subgroups $M\subseteq G$ appearing in the Hodge-Newton decomposition of $(b,\mu)$.
\end{remark*}

Let us comment on previous results in the literature. 
Before a full theory of local Shimura varieties was available, the main examples of local Shimura varieties one could work with were the ones obtained as the rigid-generic fiber of a Rapoport--Zink space \cite{RZ96} (i.e. a moduli space of $p$-divisible groups with additional structure). 
The most celebrated examples of Rapoport--Zink spaces are the Lubin--Tate tower and the tower of covers of Drinfeld's upper half space. 
In \cite{dJ95} de Jong, as an application of his theory of fundamental groups, computes the connected components of the Lubin--Tate tower for $\on{GL}_n({\bbQ_p})$. 
In \cite{MR2441699}, Strauch computes by a different method the connected components of the Lubin--Tate tower for $\on{GL}_n(F)$ and an arbitrary finite extension $F$ of ${\bbQ_p}$ (including ramification). 
In \cite{Che13}, Chen constructs $0$-dimensional local Shimura varieties and studies their geometry, these are the first examples of local Shimura varieties that do not come from a Rapoport--Zink space. 
In \cite{Chen}, Chen describes the connected components of the local Shimura varieties that arise from Rapoport--Zink spaces of EL and PEL type associated to more general unramified reductive groups. 
Our result goes beyond the previous ones in that the only condition imposed on $G$ is that it is unramified. 
In this way, our result is the first to cover very general families of local Shimura varieties that can not be constructed from a Rapoport--Zink space. 
In particular, our result is new for local Shimura varieties associated to reductive groups of exceptional types.  \\ 

The central strategy of Chen builds on and heavily generalizes the central strategy used by de Jong. 
Two key inputs for our strategy which come from Chen's work are the use of her ``generic'' crystalline representations \cite[\S 4, 5]{Chen} and her collaboration with Kisin and Viehmann on computing the connected components of affine Deligne--Lusztig varieties \cite{CKV15}. 
The present article takes these two inputs as given. 

We build on the central strategy employed by de Jong and Chen, but the versatility of Scholze's theory of diamonds \cite{Sch17} and the functorial construction of local Shimura varieties allow us to make big simplifications and streamline the proof. 
Since our arguments take place in Scholze's category of diamonds rather than the category of rigid-analytic spaces, our argument works even for moduli spaces of $p$-adic shtukas that are not local Shimura varieties. 
In these (non-representable) cases, the result is new even for $G=\on{GL}_n$.  \\

The cost of working within the framework of diamonds is that many ``classical'' constructions cannot be applied directly. 
For example, if one wished to follow the strategy of de Jong and Chen naively, one would have to develop de novo the theory of de Jong fundamental groups in the context of locally spatial diamonds.  
One observation is that the fundamental groups employed by Chen and de Jong are not really necessary.
Instead we consider stabilizer subgroups of the action of $G(\bbQ_p)$ on $\pi_0(\Sht^{\on{geo}}_{G,b,\mu,\infty})$ and related spaces. 
Morally, these stabilizers subgroups would correspond to the image of the monodromy map from the fundamental group to $G(\bbQ_p)$. 
More seriously, Chen heavily relies on the fact that the local Shimura varieties studied in \cite{Chen} are obtained as the rigid-generic fiber of a smooth formal scheme. 
This fact that Chen relies on easily implies, by \cite[Theorem 7.4.1]{dJ95b}, that the connected components of the rigid-generic fiber are in bijection to the connected components of the special fiber of the formal scheme in question.  
Since our moduli spaces no longer admit formal schemes as integral models, we tackle this point with a different strategy. 
Our new main contribution to the strategy is to use specialization maps in the context of diamonds and v-sheaves. 
To use these specialization maps in a rigorous way, we developed a formalism whose details were worked out in the separate paper \cite{Gle22}. \\

Let us sketch the central strategy to prove \Cref{theoremintheintro}. 
Once one knows that $\pi_0(\Sht^{\on{geo}}_{G,b,\mu,\infty})$ is a right ${G(\bbQ_p)_{\circ}}$-torsor, computing the actions by $W_E$ and $G_b({\bbQ_p})$ in terms of the ${G(\bbQ_p)_{\circ}}$-action can be reduced to the case of a torus. 
Indeed, this follows from functoriality of the rule 
\[(G,b,\mu)\mapsto \pi_0(\Sht^{\on{geo}}_{G,b,\mu,\infty}),\]
z-extension techniques and the determinant map from $\on{det}:G\to G^\ab$ from $G$ to its maximal abelian quotient $G^\ab$. 
To do this one can use down-to-earth diagram chases (see $\mathsection$\ref{extendedversionsofmpassection} and $\mathsection$\ref{subsubsectionactions}). 
Moreover, the case in which $G$ is a torus can be further reduced to Lubin--Tate theory as is done in \cite{Che13} (see also \Cref{theactionofJ} and \Cref{theactionofW}). \\

Let us sketch how to prove that $\pi_0(\Sht^{\on{geo}}_{G,b,\mu,\infty})$ is a right ${G(\bbQ_p)_{\circ}}$-torsor in the simplest case. 
For this, let $G$ be semisimple and simply connected. 
Our theorem then says that $\Sht^{\on{geo}}_{G,b,\mu,\infty}$ is connected.
Fix $x\in \pi_0(\Sht^{\on{geo}}_{G,b,\mu,\infty})$, we show that $G(\bbQ_p)$ acts transitively on $\pi_0(\Sht^{\on{geo}}_{G,b,\mu,\infty})$  and that the stabilizer $G_x\subseteq G(\bbQ_p)$ is the whole group.
Transitivity is proven in \Cref{transitiveofgder} and follows closely the ideas of Chen. 
For the latter part let $K\subseteq G({\bbQ_p})$ be a hyperspecial subgroup of $G$. 
We claim that it is enough to prove two things: first that $G_x$ is open and second that $G({\bbQ_p})=K\cdot G_x$. 
Indeed, $K$ surjects onto $G({\bbQ_p})/G_x$ so that this space is discrete and compact, therefore finite.
By a theorem of Margulis \cite[Chapter II, Theorem 5.1]{Marg}, since we assumed $G$ to be simply connected, the only open subgroup of finite index is the whole group, so $G_x=G({\bbQ_p})$. 
To prove $G_x$ is open we use the results of \cite{Chen} on “generic” crystalline representations, this is one of the crucial points where the HN-irreducibility of $(b,\mu)$ enters the picture. 
To prove the second input, i.e. that $G({\bbQ_p})=K\cdot G_x$, one is reduced to proving that $\Sht^{\on{geo}}_{G,b,\mu,K}$, the $K$-level moduli space of $p$-adic shtukas, is connected. 
This is where our theory of specialization maps gets used and where our strategy substantially deviates from the work of Chen. 
Addressing this point is also the heart of the paper.
This leads to discuss the main theorem of $\mathsection$\ref{sectionspecialization}. \\

Let us fix some notation. 
We fix again $G$ a reductive group over ${\bbQ_p}$ (no longer assumed to be unramified).
We fix $(G,b,\mu)$ a $p$-adic shtuka datum (no longer assumed to be HN-irreducible). 
We fix a parahoric model $\calG$ of $G$ defined over $\bbZ_p$ and let $K=\calG(\bbZ_p)$, this is a compact open subgroup of $G(\bbQ_p)$.
Let $O_{\breve{E}}\subseteq \breve{E}$ denote the ring of integers.
In this circumstance, Scholze and Weinstein construct a v-sheaf $\Shtm{O_{\breve{E}}}$ defined over $\on{Spd}(O_{\breve{E}})$ whose generic fiber is $\Sht_{G,b,\mu,K}$ (see \Cref{defi:moduliofshtukas}, \Cref{defi:moduliofshtukasbound} and \cite[\S 25]{SW20}).   
Attached to the tuple $(G,b,\mu,\calG)$ we also have an affine Deligne--Lusztig variety (see \Cref{definition-aDLV} or \cite{MR3966812} for a survey article) that we denote $X^{\leq \mu}_\calG(b)$.
We have the following result.

\begin{theorem*}[\Cref{thm:specializtheorem}]
	\label{intro:thm2}
For every nonarchimedean field extension $F$ of $\breve{E}$ contained in $\bbC_p$ the following hold.  
\begin{enumerate}
\item There is a natural continuous specialization map
\[\on{sp}:|\Sht_{G,b,\mu,K}\times \Spd F|\to |X^{\leq \mu}_\calG(b)|.\]
This map is specializing and a spectral map of locally spectral topological spaces. It is a quotient map and $G_b(\bbQ_p)$-equivariant.
\item The specialization map induces a bijection on connected components 
\[\on{sp}:\pi_0(\Sht_{G,b,\mu,K}\times \Spd F)\to \pi_0(X^{\leq \mu}_\calG(b)).\]
In particular, we have a bijection of connected components
\[\on{sp}:\pi_0(\Sht^{\on{geo}}_{G,b,\mu,K})\to \pi_0(X^{\leq \mu}_\calG(b)).\]
\end{enumerate}
\end{theorem*}

	Using \Cref{intro:thm2} above and known results in the study of connected components of affine Deligne--Lusztig varieties \cite{CKV15}, \cite{Nie18} \cite{HZ20} one can finish the proof of \Cref{theoremintheintro}. 
	Indeed, if $G$ unramified and $(b,\mu)$ is HN-irreducible $\pi_0(X^{\leq \mu}_\calG(b))$ is identified with certain subset of Borovoi's \cite{borovoi1989algebraic} fundamental group $\pi_1(G)$ \cite[Theorem 1.1]{CKV15}, \cite[Theorem 8.1]{HZ20}, \cite[Theorem 1.1]{Nie18}.
	If we go back to the assumptions of \Cref{theoremintheintro} and assume further that $G$ is semisimple and simply connected, we get that $\pi_1(G)=\{e\}$, which finishes the (sketch of) the proof of \Cref{theoremintheintro} in this particular case. 
	The proof of the general case is not very different, but it requires more patience.

	\begin{remark*}
		The first proof of \Cref{intro:thm2} appeared in one of the early versions of \cite{Gle22} in the case in which $\calG$ is a hyperspecial group.
		For editorial reasons, this theorem was removed from \cite{Gle22} and became part of the present article. 
		Although, the most immediate interest in proving \Cref{intro:thm2} was its consequence to \Cref{theoremintheintro} for which the hyperspecial case sufficed, we pursued the greater generality in this article in anticipation to our collaboration with Lim and Xu \cite{gleason2023connectedcomponentsaffinedelignelusztig} where our \Cref{intro:thm2} plays a crucial role. 
	\end{remark*}

	The proof of \Cref{intro:thm2} uses the following ingredients.
	\begin{enumerate}
		\item The machinery from integral $p$-adic Hodge theory as discussed in \cite{SW20} (see \S \ref{p-adichodgetheory}).
		\item The formalism of kimberlites as developed in \cite{Gle22} (see \S \ref{Sepcializationsections}).
		\item The main results in our collaborations \cite{AGLR22} and \cite{gleason2024tubula} (see \S \ref{Grassmanianssection}). 
	\end{enumerate}
	The construction and abstract properties of the specialization map (continuous, specializing and spectral) is an application of the theory of kimberlites developed in \cite{Gle22}. 
	In very rough terms, the theory of kimberlites addresses the question: what does it mean to be a formal scheme within Scholze's formalism of diamonds and v-sheaves?
	The theory of kimberlites selects axioms on a v-sheaf which approximate the behavior of the v-sheaves that are obtained from applying the $\diamondsuit$-functor to a formal scheme (see $\mathsection$ \ref{Sepcializationsections}).  
	In our case, we prove that the integral model $\Shtm{O_{\breve{E}}}$ proposed by Scholze and Weinstein is a smelted kimberlite (see \Cref{definitionkimberlite}, \Cref{thm:shtukasisOrapian}). 
	At heart, the main reason that $\Shtm{O_{\breve{E}}}$ ``behaves'' like a formal scheme lies on the work of Kedlaya \cite{Ked20} and Ansch\"utz' work \cite[Theorem 1.1]{Ans22} on extending vector bundles and $\calG$-torsors over the punctured spectrum of $A_{\on{inf}}$ (see \Cref{thm:anschutzgoodloci} for the version of Ansch\"utz' result that we use). 
	This already produces a specialization map 
	\[\on{sp}:|\Sht_{G,b,\mu,K}|\to |(\Shtm{O_{\breve{E}}})^\red|.\]
	Here, $(\Shtm{O_{\breve{E}}})^\red$ is the ``reduced special fiber'' or the image of the reduction functor (see \ref{equationreductionfunctor}).  
	In general, the reduced special fiber of a v-sheaf can be quite abstract and does not necessarily admit the structure of a scheme.
	Nevertheless, we construct an identification $X^{\leq \mu}_\calG(b)\simeq (\Shtm{O_{\breve{E}}})^\red$ (see \Cref{pro:adlvasreducedfunctor}). 
	The main insight is as follows.
	Recall that in rough terms $\Shtm{O_{\breve{E}}}$ parametrizes triples $(\calT,\Phi,\rho)$ where $(\calT,\Phi)$ is a shtuka with $\calG$-structure and $\rho:\calT\to \calG_b$ is $\varphi$-equivariant trivialization over $\calY_{[r,\infty]}$ for large enough $r$ (see \Cref{defi:moduliofshtukas} for details). 
	The observation is that $(\Shtm{O_{\breve{E}}})^\red$ is the locus in which $\rho$ is defined over $\calY_{(0,\infty]}$ and meromorphic along the locus $\{|p|=0\}$ \cite[Definition 5.3.5]{SW20}. 
Once this is established, constructing the isomorphism $X^{\leq \mu}_\calG(b)\simeq (\Shtm{O_{\breve{E}}})^\red$ becomes formal.
This gives a specialization map 
	\[\on{sp}:|\Sht_{G,b,\mu,K}|\to |X^{\leq \mu}_\calG(b)|.\]
	A boon of having specialization maps is that one can construct ``formal neighborhoods'' at closed points.
	These are the v-sheaf theoretic analogues of formal completions (see \Cref{formalnbbhoosoughtotbedefined}).
	Roughly, one can think of these formal neighborhoods as the subsheaves of $\Shtm{O_{\breve{E}}}$ whose points map to a given fixed point $x\in |X^{\leq \mu}_\calG(b)|$ under the specialization map.
To prove surjectivity of the specialization map and relate the connected components of the generic fiber with the connected components of the reduced special fiber, we analyze these formal neighborhoods at closed points. 
Indeed, to prove surjectivity one shows that the generic fibers of the formal neighborhoods are all non-empty. 
To show that $\on{sp}$ induces bijections of connected components one shows that the generic fibers of the formal neighborhoods are geometrically connected.
The main input to achieve this is the construction of a correspondence that relate formal neighborhoods of $\Shtm{O_{\breve{E}}}$ to formal neighborhoods of a simpler space. We clarify this below. \\ 

Before stating our last main theorem we setup some terminology and formulate a conjectural statement that is philosophically aligned with Grothendieck--Messing theory. 
Let $\Gr_\calG$ denote the Beilinson--Drinfeld Grassmannian over $\Spd \bbZ_p$ (see \Cref{defi:Grassmannian} \cite[\S 21.2]{SW20}).
The generic geometric fibers of the map $\Gr_\calG\to \Spd \bbZ_p$ are isomorphic to Scholze's affine $B_{dR}$-Grassmannian \cite[\S 19.1]{SW20}, and its reduced special fiber $\Gr_{\calG}^\red$ is the Witt vector affine flag variety \cite{Zhu17}, \cite{BS17}.
Let 
\[\Grm{O_{E}}\subseteq \Gr_{\calG,O_E}\]
denote the local model studied in \cite{AGLR22} (see \Cref{thelocalmodeldefi}) and let $\GrWmfe=(\Grm{O_{E}})^{\on{red}}$ denote its reduced special fiber.
By \cite[Theorem 6.16]{AGLR22}, the space $\GrWmfe\subseteq \Gr_{\calG,O_E}^\red$ is the $\mu$-admissible locus in the Witt vector affine flag variety \cite{KR00}, \cite[Definition 3.11]{AGLR22}. 
We let $F\supseteq \breve{E}$ be a nonarchimedean field extension with ring of integers $O_F$ and algebraically closed residue field $k_F$. 

\begin{conjecture*}
\label{intro:conjecturetubular}
For every closed point $x\in |(X^{\leq \mu}_\calG(b))_{{k_F}}|$ there exist a pair $(y,\Theta)$ where $y$ is a closed point $y\in |(\GrWmfe)_{k_{F}}|$ and 
\[\Theta:\Tf{\Shtm{O_F}}{x}\to \Tf{\Grm{O_F}}{y}\]
is an isomorphism of v-sheaves.
Here $\Tf{\Shtm{O_F}}{x}$ and $\Tf{\Grm{O_F}}{y}$ denote the formal neighborhoods as in \Cref{formalnbbhoosoughtotbedefined}.
\end{conjecture*}

The weaker version that we can prove is as follows (see \Cref{thm:comparetub}).

\begin{theorem*}[\Cref{thm:comparetub}]
	\label{intro:comparisontubular}
	Let the notation be as in \Cref{intro:conjecturetubular}, there is v-sheaf in groups $\widehat{L^+_{\bbW}\calG}$ (see \Cref{defi:moduliofmatrices}) over $\Spd k_F$ that is connected and satisfies the following. 
	For every closed point $x\in |(X^{\leq \mu}_\calG(b))_{{k_F}}|$ there exists a closed point $y\in |(\GrWmfe)_{k_{F}}|$ and a correspondence 
\begin{center}
\begin{tikzcd}
&  X \arrow{rd}{f} \arrow{dl}{g}  &  \\
\Tf{\Shtm{O_F}}{x} &  & \Tf{\Grm{O_F}}{y}
\end{tikzcd}
\end{center}
where $f$ and $g$ are both $\widehat{L^+_{\bbW}\calG}$-bundles. In particular, since $\Tf{\Grm{O_F}}{y}$ is non-empty and connected for every closed point $y\in |(\GrWmfe)_{k_{F}}|$ it follows that $\Tf{\Shtm{O_F}}{x}$ is also non-empty and connected for all $x\in |(X^{\leq \mu}_\calG(b))_{{k_F}}|$. 
\end{theorem*}

\begin{remark*}
	As with \Cref{intro:thm2}, the first version of \Cref{intro:comparisontubular} appeared in \cite{Gle22} in the case in which $\calG$ is hyperspecial.
For editorial reasons this theorem was removed from \cite{Gle22} and became part of the present article. 
The proof of \Cref{intro:comparisontubular} in the generality pursued here relies on the main theorem of our collaboration with Louren\c{c}o \cite{gleason2024tubula}. We thank him heartily for sharing his ideas with us on that project. 
\end{remark*}

\begin{remark*}
A previous version of the material surrounding \Cref{intro:comparisontubular} contained a flawed proof of \Cref{intro:conjecturetubular} which used to be a stepping stone to prove \Cref{intro:thm2}. 
The flaw was found and communicated to us by Pappas and Rapoport while they were working on \cite{PR24}. 
Soon after, we found a different argument to show \Cref{intro:thm2} by exploiting the correspondence described in \Cref{intro:comparisontubular} and avoiding the use of the difficult \Cref{intro:conjecturetubular}.
The correspondence described in \Cref{intro:comparisontubular} has been used by Pappas and Rapoport in their works \cite{PR24} and \cite{PR22}.
\end{remark*}

\begin{remark*}
During the revision process of this article there has been a lot of progress in proving \Cref{mainconjecture} whenever $\mu$ is minuscule. 
Notably, \cite{PR22} for the local abelian type case, \cite{bartling2022mathcalgmudisplayslocalshtuka}, \cite{ito2023deformation} for the hyperspecial case and \cite{takaya2024moduli} for the unramified group case. 
Although we do not have evidence for this conjecture outside the minuscule case, we still expect this to be true.
\end{remark*}

\begin{remark*}
There are plenty of cases in which $\mu$ is minuscule and $\Shtm{O_F}$ is known to be representable by a formal scheme \cite{PR22}. 
On those cases our \Cref{intro:comparisontubular} can be used to show that the formal scheme representing  $\Shtm{O_F}$ is normal.
\end{remark*}

The goal of $\mathsection$\ref{sectionspecialization} is to prove \Cref{intro:thm2} and \Cref{intro:comparisontubular}, and it is the heart of the paper. 
In $\mathsection$\ref{p-adichodgetheory} we collect the facts from integral $p$-adic Hodge theory required to define and study the specialization map for moduli spaces of $p$-adic shtukas. 
In $\S$\ref{Sepcializationsections} we explain the theory of kimberlites by highlighting the main constructions and concepts introduced in \cite{Gle22}. 
In $\S$\ref{Grassmanianssection} we summarize the results of our collaborations \cite{AGLR22} and \cite{gleason2024tubula} that we use in the present article.
In $\S$\ref{section-moduli-ofshtukas} we start the study of specialization maps for moduli spaces of $p$-adic shtukas.
Here we construct the specialization map and the identification of reduced special fibers with the affine Deligne--Lusztig varieties. 
In $\S$\ref{tubularnighborhoodssection} we finish proving \Cref{intro:thm2} and \Cref{intro:comparisontubular}. 

The goal of $\mathsection$\ref{thiisghesecondsection} is to prove \Cref{theoremintheintro}.
We start $\S$\ref{Moduliofshtukassection} by discussing the conjecture of Rapoport and Viehmann on the structure of $\pi_0(\Sht^{\on{geo}}_{G,b,\mu,\infty})$. We also recall the definition of the space $\Sht^{\on{geo}}_{G,b,\mu,\infty}$ together with its three actions and the formalism of degree $1$ divisors of Fargues--Scholze \cite{FS21} which replaces the formalism of Weil descent datum.
In $\S$\ref{thesimplestcase} we prove \Cref{theoremintheintro} in the case in which $G^\der$ is simply connected. 
In $\S$\ref{generalcase} we prove \Cref{theoremintheintro} in full generality by reducing it to the case in which $G^\der$ is simply connected. 

\subsection*{Acknowledgements.}
We thank our PhD advisor, Sug Woo Shin, for his interest, his insightful questions and suggestions at every stage of the project, and for his generous constant encouragement and support during the PhD program. 
We thank David Hansen for insightful comments and his interest in our work.
We thank Johannes Ansch\"utz, Jo\~ao Louren\c{c}o and Timo Richarz for the collaborations that this paper relies on. 
We thank Georgios Pappas and Michael Rapoport for bringing to our attention a serious flaw on an attempt we had to prove \Cref{intro:conjecturetubular} and for their interest in our work. 
We thank the anonymous referees for a very careful read of the article, and several comments that improved the quality of the paper.

The author would also like to thank Alexander Bertoloni, Rahul Dalal, Gabriel Dorfsman-Hopkins, Laurent Fargues, Zixin Jiang, Dong Gyu Lim, Sander Mack-Crane, Gal Porat, Peter Scholze, Koji Shimizu, Jared Weinstein for conversations related to the material. \\

This work was financially supported by the Doctoral Fellowship from the “University of California Institute for Mexico and the United States” (UC MEXUS) by the “Consejo Nacional de Ciencia y Tecnolog\'ia” (CONACyT), and by Deutsche Forschungsgemeinschaft (DFG, German Research Foundation) via Peter Scholze's Leibniz price.

	\section{Notation.}
	\label{sec:Notation}
	Let $p\in \bbZ$ be a prime number. 
	Given $R$ a perfect ring in characteristic $p$, we let $\on{Fr}:R\to R$ be given by $\on{Fr}(r)=r^p$. 
	We let $\Frob:\Spec R\to \Spec R$ denote $\Spec(\on{Fr})$. 
	We let $\bbW(R)$ denote the ring of $p$-typical Witt vectors. 
	We let $\phi:\bbW(R)[\frac{1}{p}]\to \bbW(R)[\frac{1}{p}]$ be the canonical lift of $\on{Fr}$.
	We let $\varphi=\Spec(\phi)$.

	We let $\Perfd$, $\Perf$ and $\Perff$ denote the category of perfectoid spaces over $\bbZ_p$, $\bbF_p$ and $\bar{\bbF}_p$ respectively \cite[Definition 3.19]{Sch17}.
	We endow $\Perf$ with the v-topology \cite[Definition 8.1]{Sch17}. 
	We denote by $\PSch$ the category of perfect schemes in characteristic $p$ endowed with the scheme-theoretic v-topology \cite[Definition 2.1]{BS17}. 
	Unless we say otherwise, the geometric objects we consider are all either perfect schemes in characteristic $p$ or small v-stacks \cite[Definition 12.1, 12.4]{Sch17}.  
	We denote by $\widetilde{\Perf}$ the category of small v-stacks \cite[Definition 12.4]{Sch17} and by $\widetilde{\PSch}$ the category of small scheme-theoretic v-stacks. 

	If $S\in \Perf$, by an \textit{untilt} of $S$ we mean a pair $(S^\sharp,i)$ where $S^\sharp\in \Perfd$ and $i:(S^\sharp)^\flat\simeq S$ is an isomorphism \cite[Definition 3.9, Corollary 3.20]{Sch17}. 
	If the context is clear, we simply write $S^\sharp$ for an untilt of $S$.
	Our convention will be to denote an untilt by a pair $(S^\sharp,\iota)$ when we are giving a definition or when we need explicitly information about $\iota$. Otherwise we will drop $\iota$ from the notation, typically while formulating statements and proving them.
	If $S=\Spa(R,R^+)$ and $S^\sharp$ is an untilt of $S$ we let $R^\sharp=\Gamma(S^\sharp,\calO_{S^\sharp})$, and we still refer to $R^\sharp$ as an untilt of $R$.

	Recall the following definition introduced in \cite[$\mathsection$ 18.1]{SW20}.
	\begin{definition}
		\label{definitionspdA}
		If $(A,A^+)$ is a complete Huber pair over $\bbZ_p$ (respectively an adic space over $\bbZ_p$) we denote by 
		\[\Spd(A,A^+):\Perf^\op\to \Sets \text{ (respectively } X^\Diamond: \Perf^\op\to \Sets \text{)}\]
	the functor with formula
	\[\Spd(A,A^+)(S)=\{((S^\sharp,\iota),f) \}/\simeq \text{ (respectively } X^\Diamond(S)=\{((S^\sharp,\iota),f) \}/\simeq \text{)}\]
	where $(S^\sharp,\iota)$ is an untilt of $S$ and $f:S^\sharp\to \Spa(A,A^+)$ (respectively $f:S^\sharp\to X$) is a map of adic spaces.
	\end{definition}
	Whenever $\Hub{A}$ is Huber pair such that $A^+=A^\circ$ we let $\Spa A=\Spa(A,A^+)$ and $\Spd A=\Spd(A,A^+)$. 
	For any perfectoid Huber pair $(R,R^+)$ we let $R^+_\red:=R^+/R^{\circ \circ}$.
	For any perfectoid Huber pair in characteristic $p$ we let 
	\[\Frob:\Spa(R,R^+)\to \Spa(R,R^+)\]
	be given as $\Frob=\Spa(\on{Fr})$.
	This formally extends to all small v-sheaves.
	Explicitly for $X=\Spd(A,A^+)$, 
	\[\on{Frob}_{X}:X\to X\]
	has formula
\[\on{Frob}_{X}[((S^\sharp,i),f)]=((S^\sharp,\on{Frob}^{-1}_S\circ i),f).\]
	If $X$ is a small v-sheaf \cite[Definition 12.1]{Sch17} (respectively a scheme) we let $|X|$ denote its associated topological space as in \cite[Proposition 12.7]{Sch17} (respectively its underlying topological space).
	If $T$ is a topological space, we denote by $\underline{T}:\Perf^\op\to \Sets$ the functor with formula
	\[\underline{T}(S)=\{f:|S|\to T\mid f\, \textrm{is continuous} \}.\]
	Given an object $X$ and a group object $G$ together with a left (respectively right) action of $G$ on $X$, 
	we use 
	\[g\star_G x \text{       (resp.        } x\star_G g \text{)}\]
	to denote the action. 
	If the context is clear, we drop $G$ from the notation and simply write $\star$.

	Throughout the text we fix an algebraic closure $\bar{\bbQ}_p$ of $\bbQ_p$, and we let $\bbC_p$ denote its $p$-adic completion.
	We let $\bar{\bbF}_p$ denote the residue field of $\bbC_p$ which is an algebraic closure of $\bbF_p$.
	Let $\breve{\bbQ}_p:=\bbW(\bar{\bbF}_p)[\frac{1}{p}]$, we regard $\breve{\bbQ}_p$ as a subfield of $\bbC_p$ after fixing an embedding.
	We let $\Gamma_{\bbQ_p}$ be the Galois group of $\bbQ_p$, we regard it as a topological group.
	Moreover, we identify it with the group of continuous automorphisms of $\bbC_p$. 
	We let $W_{\bbQ_p}\subseteq \Gamma_{\bbQ_p}$ be the Weil group, which we regard as the subgroup of continuous automorphisms of $\bbC_p$ that act on $\breve{\bbQ}_p$ by $\phi^s$ for some $s$. 

	The group $W_{\bbQ_p}$ has its standard left action on $\bbC_p$ which induces a right action on $\Spd \bbC_p$ with formula
	\[\Spd \bbC_p \times W_{\bbQ_p} \to \Spd \bbC_p\]
	\[((S^\sharp,i),f)\star_{\on{std}} w \mapsto ((S^\sharp,i),\Spa(w)\circ f).\]

	For the rest of the article we let $G$ denote a connected reductive group over $\Spec \bbQ_p$ and we let $\calG$ denote a parahoric model of $G$ defined over $\bbZ_p$, we refer the reader to \cite{BT72}, \cite{BT84}, \cite{MR4520154} for background on the theory of parahoric group schemes. 

	Throughout the text we will use the Tannakian formalism when dealing with $G$-torsors (and $\calG$-torsors) \cite{saavedra_rivano_categories_tannakiennes,deligne_categories_tannakiennes,Broshi}. 
	Namely, if $H$ is an affine algebraic group over a ring $R$ that is either a field or a Dedekind domain we let $\Rep_H$ denote the $\otimes$-exact category of algebraic representations of $H$ on finite projective $R$-modules.

	If $\calT$ is a $\otimes$-exact category, we can form the groupoid of $\otimes$-exact functors from $\Rep_H$ to $\calT$.
	We refer to the objects in this groupoid by the fixed phrase ``objects in $\calT$ with $H$-structure''. 
	We will use this mostly when $\calT=\on{Vec}_X$ i.e. the category of vector bundles on a scheme or an analytic adic space $X$ as in \cite[Appendix to Lecture 19]{SW20}. 
	We will also use this when $\calT$ is a Tannakian category \cite[\S 2.8]{deligne_categories_tannakiennes}.

	\section{Specialization maps for moduli spaces of local sthukas.}
	\label{sectionspecialization}
	\subsection{Recollections on integral $p$-adic Hodge theory.}
	\label{p-adichodgetheory}
	\subsubsection{The geometry of $\calY$.}
	\label{geometryo-f-Y}
	Recall \cite[$\mathsection$ 11.2]{SW20} that one can associate to any $S\in \Perf$ an analytic sous-perfectoid adic space \cite[\S 6.3]{SW20} over $\Spa \bbZ_p$ denoted by $``S\dot{\times} \Spa \bbZ_p"$ that represents the diamond $S\times \Spd \bbZ_p$. 
	The formula for this space when $S=\Spa(R,R^+)$ and $\varpi\in R^+$ is a pseudo-uniformizer is 
	\[``S\dot{\times} \Spa \bbZ_p" = \Spa(\bbW(R^+))\setminus \{[\varpi]=0\}.\]
	Since we will use some variants of these spaces parametrized by intervals in $[0,\infty]$ we will use the following notation instead.
	Our notation agrees with the one used in \cite[$\mathsection$ 12.2]{SW20}.
\begin{definition}
	\label{defi:Y}
	Given a perfectoid Huber pair $(R,R^+)$ in characteristic $p$ and a pseudo-uniformizer $\varpi\in R^+$, we let $\calY_{[0,\infty)}^{R^+}$ denote $\Spa{\bbW(R^+)}\setminus V([\varpi])$. Here $[\varpi]$ denotes a Teichm\"uller lift of $\varpi$, and $\bbW(R^+)$ is given the $(p,[\varpi])$-adic topology. We let $\Yint{R^+}{[0,\infty]}$ denote $\Spa{\bbW(R^+)}\setminus V(p,[\varpi])$.
\end{definition}
We review the geometry of $\Yint{R^+}{[0,\infty]}$. 
Fix a pseudo-uniformizer $\varpi\in R^+$.
One defines a continuous map \[\kappa_{\varpi}: |\Yint{R^+}{[0,\infty]}|\to [0,\infty]\]
such that $\kappa_\varpi(y)=r$ if and only if for all positive rational numbers with $\frac{m}{n}\leq r\leq \frac{M}{N}$ the inequalities 
\[|p|_y^M\leq |[\varpi]|_y^N\ \text{   and    } |[\varpi]|_y^n\leq |p|_y^m\] 
hold.
We often leave the pseudo-uniformizer implicit and omit it from the notation.
Given an interval $I\subseteq [0,\infty]$ we denote by $\Yint{R^+}{I}$ the open subset corresponding to the interior of $\kappa_{\varpi}^{-1}(I)$. 
For intervals of the form $[0,\frac{h}{d}]$ where $h$ and $d$ are positive integers the space $\Yint{R^+}{[0,\frac{h}{d}]}$ is represented by $\Spa\Hub{R_{h,d}}$ corresponding to the rational localization, $\{x\in \Spa{\bbW(R^+)}\mid |p^h|_x \leq |[\varpi]^d|_x\neq 0\}$.
In this case, we can compute a ring of definition $R_{h,d}^0\subseteq R_{h,d}$ explicitly as the $[\varpi]$-adic completion of $\bbW(R^+)[\frac{p^h}{[\varpi]^d}]$, then $R_{h,d}$ is simply $R_{h,d}^0[\frac{1}{[\varpi]}]$. 
A direct computation shows that $R_{h,d}$ does not depend on the choice of integral elements $R^+\subseteq R$. 
In particular, the exact category of vector bundles over $\Yinf{R^+}$ does not depend on the choice of $R^+$ either \cite[Theorem 5.2.8]{SW20} \cite[Theorem 2.7.7, Remark 2.5.23]{KL15}.

Recall the ``algebraic version'' of $\Yint{R^+}{[0,\infty]}$, which we denote $\on{Y}_{[0,\infty]}^{R^+}$ and define as $\on{Spec}(\bbW(R^+))\setminus V(p,[\varpi])$. 
Since $\bbW(R^+)\subseteq \calO_{\calY^{R^+}_{[0,\infty]}}$ and since $p$ and $[\varpi]$, do not vanish simultaneously on $\calY_{[0,\infty]}^{R^+}$ we get a map of locally ringed spaces $f:\calY_{[0,\infty]}^{R^+}\to \on{Y}^{R^+}_{[0,\infty]}\subseteq \on{Spec}(\bbW(R^+))$.

Recall the following GAGA-type result of Kedlaya for vector bundles on $\on{Y}^{R^+}_{[0,\infty]}$.

\begin{theorem}\textup{(\cite[Theorem 3.8]{Ked20})}
	\label{thm:gagatypekedlay}
	Suppose $(R,R^+)$ is a perfectoid Huber pair in characteristic $p$. 
	The natural morphism of locally ringed spaces $f:\calY^{R^+}_{[0,\infty]}\to \on{Y}_{[0,\infty]}^{R^+}$ induces, via the pullback functor 
	\[f^*:\on{Vec}_{\on{Y}_{[0,\infty]}^{R^+}}\to \on{Vec}_{\calY^{R^+}_{[0,\infty]}}\]
	a $\otimes$-exact equivalence of $\otimes$-exact categories. 
\end{theorem}

\begin{remark}	
	\label{rem:exactness}	
 	Although the reference does not explicitly claim that this equivalence is exact, one can simply follow the proof loc. cit. exchanging the word ``equivalence'' by ``exact equivalence'' since every arrow involved in the proof is an exact functor.
\end{remark}	

Recall that $\calG$ denotes a parahoric group scheme, in particular it is smooth over $\Spec \bbZ_p$. 
Recall that $\on{Rep}_\calG$ denotes the $\otimes$-exact category of algebraic representations on finite free $\bbZ_p$-modules. 
One can define the category of $\calG$-torsors on a scheme $X$ (respectively on a sous-perfectoid adic space $Y$) as the category of $\otimes$-exact functors with source category $\on{Rep}_\calG$ and target category $\on{Vec}_X$ (respectively $\on{Vec}_Y$). 
By \cite[Theorem 19.5.1, Theorem 19.5.2]{SW20} and \cite{Broshi}, this definition agrees with other more natural definitions of $\calG$-torsors that do not play a role in this article.
With this definition of $\calG$-torsors one can immediately generalize \Cref{thm:gagatypekedlay}. 

\begin{corollary}
	\label{g-torsorsgaga}
	Let the notation be as in \Cref{thm:gagatypekedlay}, then pullback $f^*$ induces an equivalence from the category of $\calG$-torsors over $\on{Y}_{[0,\infty]}^{R^+}$ to the category of $\calG$-torsors over $\calY^{R^+}_{[0,\infty]}$.
\end{corollary}

Recall that if $X$ is an analytic adic space and $Z\subseteq X$ is a closed Cartier divisor \cite[Definition 5.3.2]{SW20} one can define what it means for a global section of $\Gamma(X\setminus Z, \calO_{X\setminus Z})$ to be meromorphic along $Z$ \cite[Definition 5.3.5]{SW20}. 
Given an untilt $R^\sharp$ of $R$ there is a canonical surjection $\bbW(R^+)\to R^{\sharp +}$ whose kernel is generated by an element $\xi\in \bbW(R^+)$ primitive of degree $1$ \cite[Lemma 6.2.8]{SW20}. 
The element $\xi$ defines a closed Cartier divisor over $\calY_{[0,\infty]}^{R^+}$ \cite[Proposition 11.3.1]{SW20} and also defines a Cartier divisor on the scheme $\on{Y}_{[0,\infty]}^{R^+}$. 

Let us fix some notation. Fix $S=\Spa\Hub{R}$ with $S\in \Perf$ and $S^\sharp=\Spa(R^\sharp,R^{\sharp,+})$ an untilt of $S$. Let $\xi$ be a generator for the kernel of the surjection $\bbW(R^+)\to R^{\sharp,+}$.
We let \[(\on{Vec}^{\xi\neq 0}_{\calY^{R^+}_{[0,\infty]}})^{\on{mer}} \text{   (respectively    } (\on{Vec}^{\xi\neq 0}_{Y^{R^+}_{[0,\infty]}})\text{)}\] denote the category whose objects are vector bundles over $\calY^{R^+}_{[0,\infty]}$ (respectively vector bundles over $Y_{[0,\infty]}^{R^+}$) and morphisms are vector bundle maps over $\calY^{R^+}_{[0,\infty]}\setminus S^\sharp$ that are meromorphic along the ideal $S^\sharp$ (respectively vector bundle maps over $Y_{[0,\infty]}^{R^+}\setminus V(\xi)$). 
One gets the following direct generalization of \Cref{thm:gagatypekedlay}. 

\begin{corollary}
	\label{cor:gagatypekedlymeoro}
	Let the notation be as above, and $f$ as in \Cref{thm:gagatypekedlay}, then pullback $f^*$ induces an equivalence of categories
	\[f^*:(\on{Vec}^{\xi\neq 0}_{Y^{R^+}_{[0,\infty]}})\to (\on{Vec}^{\xi\neq 0}_{\calY^{R^+}_{[0,\infty]}})^{\on{mer}}.\]
\end{corollary}

\begin{remark}
	One can also formulate and show a version of \Cref{g-torsorsgaga} for $\calG$-torsors defined over $\on{Y}_{[0,\infty]}^{R^+}$ with morphisms defined over $\on{Y}_{[0,\infty]}^{R^+}\setminus V(\xi)$.
	Indeed, it suffices to pass to functor categories $\on{Fun}_{\on{ex}}^{\otimes}(\on{Rep}_\calG,-)$ with values on the two categories that appear in \Cref{cor:gagatypekedlymeoro}.
\end{remark}

Kedlaya proves another important statement.
\begin{theorem}\textup{(\cite[Lemma 2.3, Theorem 2.7, Remark 3.11]{Ked20})}
	\label{thm:puncturedtoAinfked}
	\label{detail:checkkedlayasproof}
	Let $j:\on{Y}_{[0,\infty]}^{R^+}\to \on{Spec}(\bbW(R^+))$ denote the natural open immersion, the following statements hold.
\begin{enumerate}
	\item The pullback functor $j^*:\on{Vec}_{\on{Spec}(\bbW(R^+))}\to \on{Vec}_{\on{Y}_{[0,\infty]}^{R^+}}$ is fully-faithful. 
	\item If $R^+$ is a valuation ring then $j^*$ is an equivalence.
	\item For all vector bundles $\calV\in \on{Vec}_{\on{Y}_{[0,\infty]}^{R^+}}$ the quasi-coherent sheaf $j_*\calV$ over $\on{Spec}(\bbW(R^+))$ satisfies that the adjunction map $j^*j_*\calV\to \calV$ is an isomorphism.
\end{enumerate}
\end{theorem}
We will need a small modification of \Cref{thm:puncturedtoAinfked}.
\begin{definition}
	\label{defi:prodpoints}
	Given a set $I$ and a collection of tuples $\{(C_i,C^+_i),\varpi_i\}_{i\in I}$ we construct a perfectoid adic space $\Spa{(R,R^+)}$. 
	Here each $C_i$ is an algebraically closed nonarchimedean field in characteristic $p$, the $C_i^+$ are open and bounded valuation subrings of $C_i$, and $\varpi_i\in C_i^+$ is a choice of pseudo-uniformizer. 
	We let $R^+:=\prod_{i\in I} C^+_i$, we let $\varpi=(\varpi_i)_{i\in I}$, we endow $R^+$ with the $\varpi$-adic topology and we let $R:=R^+[\frac{1}{\varpi}]$. 
	Any space constructed in this way will be called a \textit{product of points}.
\end{definition}
The following statement is implicitly used and proved in (\cite[Theorem 25.1.2]{SW20}).
\begin{proposition}
	\label{pro:puncturedAinfprodofpoints}
	Let $\Spa(R,R^+)$ be the product of points associated to $\{(C_i,C_i^+), \varpi_i\}_{i\in I}$ as in \Cref{defi:prodpoints}. The pullback functor $j^*:\on{Vec}_{\on{Spec}(\bbW(R^+))}\to \on{Vec}_{\on{Y}_{[0,\infty]}^{R^+}}$ gives an equivalence of categories of vector bundles with fixed rank. 
	In other words, for $\calE\in \on{Vec}_{\on{Y}_{[0,\infty]}^{R^+}}$ the quasi-coherent sheaf $j_*\calE$ is a vector bundle.
\end{proposition}
Fix $\xi\in \bbW(R^+)$ primitive of degree $1$ as before and recall that both $\on{Spec}(\bbW(R^+))$ and $Y_{[0,\infty]}^{R^+}$ are qcqs schemes.
Consequently, the equivalence of vector bundles of \Cref{pro:puncturedAinfprodofpoints} generalizes to an equivalence similar in form to the one described in \Cref{cor:gagatypekedlymeoro}. 
Namely, it is an equivalence of those categories whose objects are as in \Cref{pro:puncturedAinfprodofpoints}, but whose morphisms are allowed to have poles along $\xi$ on both categories.  

Interestingly, extending $\calG$-torsors from $\on{Y}_{[0,\infty]}^{R^+}$ to $\on{Spec}(\bbW(R^+))$ adds yet another layer of complexity. 
Indeed, the equivalences of \Cref{thm:puncturedtoAinfked} and \Cref{pro:puncturedAinfprodofpoints} are not exact equivalences, so the Tannakian formalism can't be used naively. 
As a matter of fact, only the pullback functor $j^*$ is exact. 
Ansch\"utz gives a detailed study of the problem of extending $\calG$-torsors along $j$ for parahoric group schemes $\calG$ \cite{Ans22}.
\begin{theorem}\textup{(\cite[Theorem 1.1]{Ans22})}
	\label{pro:prodpointextendg-tors}
	\label{thm:anschutzgoodloci}
Let $\Spa(R,R^+)$ be a product of points. 
Every $\calG$-torsor $\calT$ over $\on{Y}_{[0,\infty]}^{R^+}$ extends along $j:\on{Y}^{R^+}_{[0,\infty]}\to \on{Spec}(\bbW(R^+))$ to a unique $\calG$-torsor over $\on{Spec}(\bbW(R^+))$. 
\end{theorem}
\begin{proof}
	The reference only states explicitly the case $R^+=O_C$ with $O_C\subseteq C$ the ring of integers in an algebraically closed non-Archimedean field $C$. 
	Nevertheless, the reference \cite{Ans22} already provides the technical tools to conclude more generally.
	Indeed, if $R^+=C^+$ with $C^+\subseteq C$ a more general open and bounded valuation subring we can argue by pointing out that the proof of \cite[Proposition 9.2]{Ans22} goes through in this generality, that the first part of \cite[Corollary 9.3]{Ans22} also holds in this case and by appealing to \cite[Proposition 11.5]{Ans22}.

	The general case can be done as follows.
	We need to prove that the functor $j_*\calT:\on{Rep}_\calG\to \on{Vec}_{\on{Spec}(\bbW(R^+))}$ is right-exact, since it is always left-exact. 
	By the case $R^+=C^+$ discussed above, it suffices to show that if a morphism of finite free modules $g:\calV_1\to \calV_2$ over $\on{Spec}(\bbW(R^+))$ satisfies that each base change $g_i:\calV_{1,i}\to \calV_{2,i}$ to $\on{Spec}(\bbW(C^+_i))$ is surjective for every $i\in I$, then $g$ is also surjective. 
	Taking determinant bundles we can reduce to the case that $\calV_2$ is free of rank $1$. 
	After taking trivializations we have $n$ sections $f_1,\cdots, f_n\in \bbW(R^+)$ and we need to prove that they generate the unit ideal. 
	Consider the family of subsets $\{I_m\}_{1\leq m\leq n}$ defined by 
	\[I_m=\{i\in I\mid f_m\in \bbW(C^+_i)^\times\}.\]

	Observe that $\pi_0(\Spec \bbW(R^+))\simeq \beta I$ is the Stone–\v{C}ech compactification of $I$.
	Let $e_{I_m}$ denote the idempotent associated to the closed open subset $\beta I_m\subseteq \beta I$.
	Observe that $e_{I_m}$ is in the ideal generated by $f_m$ for every $1\leq m\leq n$. 
	Since each $\bbW(C^+_i)$ is a local ring and the $\{f_m\}_{1\leq m\leq n}$ generate the unit ideal in $\bbW(C^+_i)$ for each fixed $i\in I$, the union $\bigcup_{i=1}^n I_m$ has to be $I$ and in particular the set $\{e_{I_m}\}_{1\leq m\leq n}$ generates the unit ideal. 
	Consequently, the set $\{f_m\}_{1\leq m\leq n}$ also generates the unit ideal. 
\end{proof}

In what follows we will define several geometric objects all of which are v-sheaves or v-stacks.
The way to show that these objects are v-sheaves or v-stacks is to use systematically the following descent result. 
\begin{proposition}\textup{(\cite[Proposition 19.5.3]{SW20})}
	\label{pro:descendalongcurve}
	Let $S$ be a perfectoid space in characteristic $p$ and let $U\subseteq ``S\dot{\times} \Spa \bbZ_p"$ be an open subset. 
	For a map of perfectoid spaces $f:S'\to S$, let $\mathcal{C}_{S'}$ denote the category of $\calG$-torsors over 
	\[``S'\dot{\times} \Spa \bbZ_p"\times_{``S\dot{\times} \Spa \bbZ_p"} U.\] 
	The rule $S'\mapsto \mathcal{C}_{S'}$, as a fibered category over $\on{Perf}_S$, is a v-stack.
\end{proposition}

\subsubsection{Lattices and $p$-adic shtukas.}
\label{subsection-lattices-and-padicshtukas}
For this subsection we let $\Spa(R,R^+)$ denote an affinoid perfectoid space in characteristic $p$, let $\varpi\in R^+$ a choice of pseudo-uniformizer, let $(R^\sharp,\iota)$ be an untilt of $R$ and $\xi_{R^\sharp}$ a generator for the kernel of the map $\bbW(R^+)\to R^{\sharp,+}$.
\begin{definition}
	\label{defi:latticesGstructure}
	We define the groupoid of $\on{B}^+_{\on{dR}}(R^\sharp)$-lattices with $\calG$-structure to have as objects pairs $(\calT,\psi)$ where $\calT$ is a $\calG$-torsor over $\calY_{[0,\infty)}^{R^+}$ and $\psi:\calT\to\calG$ is an isomorphism over $\calY_{[0,\infty)}^{R^+}\setminus V(\xi_{R^\sharp})$ that is meromorphic along $(\xi_{R^\sharp})$. 
	Morphisms are defined in the natural way.
Note that morphisms in this category, if they exist, are unique.
Compare with \cite[Definition 20.3.1]{SW20}.
\end{definition}

Analogously, we consider the groupoid of $p$-adic shtukas with $\calG$-structure over $\bbZ_p$. 
Recall that the spaces $\on{Spec}\bbW(R^+)$, $\calY_{[0,\infty)}^{R^+}$, $\on{Y}_{[0,\infty]}^{R^+}$ and $\calY_{[0,\infty]}^{R^+}$ come equipped with a Frobenius action which we denote by $\varphi=\Spa(\phi)$ (or $\varphi=\Spec(\phi)$), induced from the ring homomorphism $\phi:\bbW(R^+)\to \bbW(R^+)$ discussed in \S\ref{sec:Notation}.
 \begin{definition}
	\label{defi:singleshtukasGstructure}
	We define the groupoid of $p$-adic shtukas with $\calG$-structure with one paw (or leg) over $\on{Spa}(R^\sharp,R^{\sharp,+})$.
Objects are pairs $(\calT,\Phi)$ where $\calT$ is a $\calG$-torsor over $\Yinf{R^+}$ and $\Phi:\varphi^{*}\calT\to\calT$ is an isomorphism over $\Yinf{R^+}\setminus V(\xi_{R^\sharp})$ meromorphic along $(\xi_{R^\sharp})$. 
Morphisms are given by $\varphi$-equivariant isomorphism of $\calG$-torsors.
Compare with \cite[Definition 11.4.1]{SW20}.
\end{definition}

\begin{definition}
	\label{defi:phi-modules}
	If $B$ is a ring on which $\phi$ acts, by a $\varphi$-module over $\Spec B$ (resp. $\Spa (B,B^+)$) we mean a pair $(\calE,\Phi)$ where $\calE$ is a vector bundle over $\Spec B$ (resp. $\Spa(B,B^+)$) together with an isomorphism $\Phi:\varphi^*\calE \to \calE$.
	Similarly, if we have spaces $X\subseteq Y$ and an automorphism $\varphi:Y\to Y$ with the property that $\varphi(X)\subseteq X$ we define a $\varphi$-module over $X$ to be a pair $(\calE,\Phi)$ where $\calE$ is a vector bundle over $X$ and $\Phi:\varphi^*\calE\to \calE$ is an isomorphism. 
	Finally, with the setup as above, by a $\varphi$-module with $\calG$-structure we mean a $\otimes$-exact functor from $\on{Rep}_\calG$ to the category of $\varphi$-modules.
	Compare with \cite[Definition 12.3.3]{SW20}.
\end{definition}
\begin{example}
	\label{bundlesonffcurvegobrr}
	Let $S=\Spa(R,R^+)$. 
	Since the action of $\varphi$ on $\calY^{R^+}_{(0,\infty)}$ is free and totally discontinuous \cite[Page 136]{SW20} the category of $\varphi$-modules over $\calY^{R^+}_{(0,\infty)}$ is equivalent to the category of vector bundles on the relative Fargues--Fontaine curve $X_{\on{FF},S}=\calY^{R^+}_{(0,\infty)}/\varphi^\bbZ$ \cite[Definition II.1.15]{FS21}.	
\end{example}

	Recall the categories of isocrystals $\Isoc_{\bar{\bbF}_p}$ and of isocrystals with $\calG$-structure \cite{Kot85}, \cite[$\mathsection$ 3]{Kot97}.
	The objects of $\Isoc_{\bar{\bbF}_p}$ are pairs $(V,\Phi)$ where $V$ is a finite dimensional $\breve{\bbQ}_p$-vector space and 
\[\Phi:V\to V\]
is a $\phi$-linear isomorphism. 

As usual, isocrystals with $G$-structure as in \cite[$\mathsection$ 3]{Kot97} \cite[Definition III.2.1]{FS21} are $\otimes$-exact functors
\[\calF:\on{Rep}_G\to \Isoc_{\bar{\bbF}_p}.\]
Recall Kottwitz' set of $\phi$-conjugacy classes in $G(\breve{\bbQ}_p)$ \cite[$\mathsection$ 1.4]{Kot97},
\[B(G):=\frac{G(\breve{\bbQ}_p)}{\on{Ad}_\phi G(\breve{\bbQ}_p)}.\]
To any element $b\in G(\breve{\bbQ}_p)$ one can attach an isocrystal with $G$-structure $V_b$, and $V_{b_1}$ is isomorphic to $V_{b_2}$ if and only if $b_1$ and $b_2$ represent the same class in $B(G)$.
Moreover, since $G$ is connected and reductive, every isocrystal with $G$-structure is isomorphic to $V_b$ for some $b$ \cite[$\mathsection$ 3.1]{Kot97}.

Isocrystals give rise to $\varphi$-modules by considering $V$ as a vector bundle over $\Spec \breve{\bbQ}_p$ (respectively over $\Spa \breve{\bbQ}_p$) and by interpreting $\Phi$ as a linear isomorphism of the form 
\[\Phi:\varphi^*V\to V.\]
Furthermore, if $\Spa\Hub{R}\in \Perff$ then we have a $\varphi$-equivariant map $\calY^{R^+}_{(0,\infty]}\to \Spa \breve{\bbQ}_p$.
Pullback along this map defines a $\otimes$-exact functor from the category of isocrystals to the category of $\varphi$-modules.
In particular to any isocrystal with $G$-structure $\calF$ we can associate $G_\calF$ which is a $\varphi$-module with $G$-structure over $\calY^{R^+}_{(0,\infty]}$.

\begin{remark}
	\label{modulephistructureg}
	Given $b\in G(\breve{\bbQ}_p)$ and $\Spa\Hub{R}\in \Perff$ we use $(\calG_b,\Phi_b)$ or $(G_b,\Phi_b)$ to denote the $\varphi$-module with $G$-structure on $\calY^{R^+}_{(0,\infty]}$ (or $\calY^{R^+}_{(0,\infty)}$) associated to the isocrystal with $G$-structure $V_b$. 
	Under the equivalence explained in \Cref{bundlesonffcurvegobrr} the pair $(\calG_b,\Phi_b)$ over $\calY^{R^+}_{(0,\infty)}$ corresponds to a vector bundle on $\calY^{R^+}_{(0,\infty)}/\varphi^\bbZ$ which we denote by $\calE_b$. 
	This is the same convention as in \cite[$\mathsection$ III.2.1]{FS21}.
\end{remark}

\begin{definition}
	\label{defi:isogenies}
	Given a $\varphi$-module with $\calG$-structure $(\calE,\Phi_\calE)$ over $\calY^{R^+}_{(0,\infty)}$ as in \Cref{defi:phi-modules} and a shtuka $(\calT,\Phi_\calT)$ as in \Cref{defi:singleshtukasGstructure} we define an \emph{isogeny} from $(\calT,\Phi_\calT)$ to $(\calE,\Phi_\calE)$ to be an equivalence class of pairs $(r,f)$ with $r\in \bbR$ and $f:(\calT,\Phi_\calT)\to (\calE,\Phi_\calE)$ a $\varphi$-equivariant isomorphism defined over $\calY^{R^+}_{[r,\infty)}$. 
	Two pairs $(r_1,f_1)$ and $(r_2,f_2)$ are equivalent if there is a third pair $(r_3,f_3)$ with $r_3>r_1,r_2$ and $f_1=f_3=f_2$ when restricted to $\calY^{R^+}_{[r_3,\infty)}$. 
\end{definition}

\begin{remark}
	Recall that for all $r>0$ the subspace $\calY^{R^+}_{[r,\infty)}\subseteq \calY^{R^+}_{(0,\infty)}$ contains a fundamental domain for the $\varphi$-action. One can use this to attach to a shtuka with $\calG$-structure $(\calT,\Phi_\calT)$ a $\calG$-torsor on $X_{\on{FF},S}$ that we may denote $\calE_\calT$. 
	As mentioned in \Cref{bundlesonffcurvegobrr}, the category of $\varphi$-module with $\calG$-structure over $\calY^{R^+}_{(0,\infty)}$ is equivalent to the category of $\calG$-torsors over $X_{\on{FF},S}$. 
	If $\calE$ is the $\calG$-torsor over $X_{\on{FF},S}$ corresponding to the $\varphi$-module with $\calG$-structure $(\calE,\Phi)$, then isogenies from $(\calT,\Phi_\calT)$ to $(\calE,\Phi_\calE)$ as in \Cref{defi:isogenies} are in natural bijection with isomorphisms of $\calG$-torsors between $\calE_\calT$ and $\calE$ over $X_{\on{FF},S}$. 
\end{remark}

In what follows, we prove some technical lemmas that intuitively speaking allow us to “deform” lattices and shtukas with $\calG$-structure. 
For any $r\in [0,\infty)$ let $B^{R^+}_{[r,\infty]}=\on{H}^0(\Yint{R^+}{[r,\infty]},\calO_{\Yint{R^+}{[r,\infty]}})$, and consider the ring ${R_\red^+}:=(R^+/\varpi)^{\on{perf}}=R^+/R^{\circ \circ}$ endowed with the discrete topology. 

Fix $r=\frac{q_1}{q_2}$. 
Since $[\varpi]=0$ in $\Spa \bbW(R^+_\red)[\frac{1}{p}]$ and $\Yint{R^+}{[r,\infty]}$ is the rational subset of those $x\in\Spa{\bbW(R^+)}$ for which $|[\varpi^{q_2}]|_x\leq |p^{q_1}|_x\neq 0$ holds, we have a family of ring maps 
\[
(-)_\red:B^{R^+}_{[r,\infty]}\to \bbW({R_\red^+})[{1}/{p}]
\] 
that is compatible with the natural ring maps $B^{R^+}_{[r,\infty]}\to B^{R^+}_{[r',\infty]}$ for $r\leq r'$. 
By abuse of notation we also denote $(-)_\red:R^+\to R^+_\red$ and $(-)_\red:\bbW(R^+)\to \bbW(R^+_\red)$ the natural ring maps.
Note that we have the following commutative diagram
\begin{center}
\begin{tikzcd}
	\bbW(R^+) \arrow{r} \arrow{rd}{(-)_\red}  & B^{R^+}_{[r,\infty]} \arrow{d}{(-)_\red} \\
				       &  \bbW(R^+_\red).
\end{tikzcd}
\end{center}
\begin{lemma}
	\label{lem:smallintegraldecomp}
	\label{detail:makingrbigger}
	Let $s\in B^{R^+}_{[r,\infty]}$ and suppose that ${s}_\red\in \bbW(R^+_\red)[\frac{1}{p}]$ lies in $ \bbW(R^+_\red)$. 
	Then there is a tuple $(m,a,b,\varpi_s)$ with $m\in \bbN$ a number $r\leq m$, $a\in \bbW(R^+)$, $b\in B^{R^+}_{[m,\infty]}$ and a pseudo-uniformizer $\varpi_s\in R^+$ such that $s=a+b$ in $B^{R^+}_{[m,\infty]}$ and $b\in [\varpi_s]\cdot B^{R^+}_{[m,\infty]}$. 
	Moreover, if $\varpi\in R^+$ is any pseudo-uniformizer we may choose $\varpi_s$ so that $\varpi\in \varpi_s \cdot R^+$.
\end{lemma}
\begin{proof}
	Choose $m\in \bbN$ with $r\leq m$, 
	we compute $B^{R^+}_{[m,\infty]}$ explicitly. 
	If $L_0$ denotes the $p$-adic completion of $\bbW(R^+)[\frac{[\varpi]}{p^m}]$, then $B^{R^+}_{[m,\infty]}=L_0[\frac{1}{p}]$. 
	Any element $s\in B^{R^+}_{[m,\infty]}$ is of the form $s=\frac{1}{p^n}\cdot \ell$ where $\ell\in L_0$.
	In turn, any element $\ell$ is in the image of the map
	\[\theta:\bbW(R^+)\langle T \rangle \to L_0 \text{ with } T\mapsto \frac{[\varpi]}{p^m}.\]
	Let $f(T)\in \bbW(R^+)\langle T \rangle$ with image $\ell$, and write $f(T)=f_0+T\cdot \Sigma_{i=1}^\infty f_i T^{i-1}$ with $f_0,f_i\in \bbW(R^+)$.
	Let $d(T)=\Sigma_{i=1}^\infty f_i T^{i-1}$ so that $f(T)=f_0+T\cdot d(T)$, then $\theta(f)=f_0+[\varpi](\frac{1}{p^m}\cdot \theta(d(T)))$.
Since the second term is divisible by $[\varpi]$ in $B^{R^+}_{[m,\infty]}$ as long as we pick a $\varpi_s$ that divides $\varpi$, we may and do reduce to the case $\ell=f_0$ or in other words 
\[s=\frac{1}{p^n}f_0=\Sigma_{i=0}^\infty [a_i]p^{i-n}.\] 
In this case, ${s}_\red=\Sigma_{i=0}^\infty [{{(a_i)}}_\red]p^{i-n}$ and by hypothesis we have that for $i< n$ ${{(a_i)}_\red}=0$ in ${R_\red^+}$. 
We can choose a pseudo-uniformizer $\varpi_s$ such that all of the $a_i$ for $i\in \{0,\dots n-1\}$ are zero in $R^+/\varpi_s$. 
We can take $a=\sum_{i=n}^\infty [a_i]p^{i-n}$ and $b=\sum_{i=0}^{n-1}[a_i]p^{i-n}$. 
These clearly satisfy the properties we were looking for.
\end{proof}
We will need a small improvement of \Cref{lem:smallintegraldecomp}. 
\begin{lemma}
	\label{lemma:random}
	Let $X$ be an affine scheme smooth over $\Spec \bbZ_p$.
	Suppose we have a commutative diagram
	\begin{center}
	\begin{tikzcd}
		\Spec \bbW(R^+_\red)[\frac{1}{p}] \arrow{r} \arrow{d}  & \Spec B^{R^+}_{[r,\infty]} \arrow{d}{f} \\
	  \Spec \bbW(R^+_\red) \arrow{r} & X 
	\end{tikzcd}
	\end{center}
	over $\Spec \bbZ_p$.
	Then there exists $r\leq m$, $\varpi\in R^+$ a pseudo-uniformizer, and a map $g:\Spec \bbW(R^+)\to X$ fitting in the following commutative diagram
	\begin{center}
	\begin{tikzcd}
		\Spec B^{R^+}_{[m,\infty]}/[\varpi]  \arrow{r} \arrow{d}  & \Spec B^{R^+}_{[r,\infty]} \arrow{d}{f} \\
		\Spec \bbW(R^+) \arrow{r}{g} & X. 
	\end{tikzcd}
	\end{center}
\end{lemma}
\begin{proof}
	The case $X=\bbA^1$ is \Cref{lem:smallintegraldecomp} and one can easily adapt the argument to the case $X=\bbA^n$.
	We now consider the general case, for this we fix a closed immersion $\iota:X\hookrightarrow \bbA^n$.
	Using the case $X=\bbA^n$ we get a commutative diagram
	\begin{center}
	\begin{tikzcd}
		\Spec B^{R^+}_{[m,\infty]}/[\varpi]  \arrow{r} \arrow{d}  & \Spec B^{R^+}_{[r,\infty]} \arrow{d}{f} \\
		\Spec \bbW(R^+)/[\varpi]\ar{d} \ar{dr}{\tilde{g}_\varpi} 							   & X \arrow{d}{\iota}  \\
		\Spec \bbW(R^+) \arrow{r}{\tilde{g}} & \bbA^n 
	\end{tikzcd}
	\end{center}
	for some $\tilde{g}:\Spec \bbW(R^+)\to \bbA^n$. 
	Moreover, since the map of rings $\bbW(R^+)/[\varpi]\to B^{R^+}_{[m,\infty]}/[\varpi]$ is injective and $\Spec B^{R^+}_{[m,\infty]}/[\varpi]$ factors through $X$ then $\tilde{g}_\varpi$ factors through $X$ and defines a map $g_0:\Spec \bbW(R^+)/[\varpi] \to X$.
	Since $\bbW(R^+)$ is $(p,[\varpi])$-adically complete in particular it is also $[\varpi]$-adically complete \cite[Lemme 1.4.14]{FF18}. 
	Finally, since $X$ is smooth we may find a lift $g:\Spec \bbW(R^+)\to X$ making the following diagram commutative

	\begin{center}
	\begin{tikzcd}
		\Spec B^{R^+}_{[m,\infty]}/[\varpi]  \arrow{r} \arrow{d}  & \Spec B^{R^+}_{[r,\infty]} \arrow{d}{f} \\
		\Spec \bbW(R^+)/[\varpi]\ar{d}  \ar{r}{g_0}							   & X   \\
		\Spec \bbW(R^+)  \ar{ru}{g} &  
	\end{tikzcd}
	\end{center}
	This is the diagram that we wished to construct.
\end{proof}
We apply \Cref{lemma:random} in the following particular case.
\begin{lemma}
	\label{lem:smallintegraldecomptorsors}
	Let $\calT_1$ and $\calT_2$ be trivial $\calG$-torsors over $\on{Spec}(\bbW(R^+))$ and let $\lambda:\calT_1\to \calT_2$ be an isomorphism over $\Yint{R^+}{[r,\infty]}$ whose reduction to $\on{Spec}(\bbW(R^+_\red)[\frac{1}{p}])$ extends to $\on{Spec}(\bbW(R^+_\red))$. 
	Then, there is an isomorphism $\widetilde{\lambda}:\calT_1\to \calT_2$ over $\on{Spec}(\bbW(R^+))$, a pseudo-uniformizer $\varpi_\lambda\in R^+$ and a number $r\leq m$ such that the restriction of $\lambda$ and of $\widetilde{\lambda}$ to $\on{Spec}(B^{R^+}_{[m,\infty]})$ agree as elements of $\on{Hom}_{\on{Spec}(B^{R^+}_{[m,\infty]}/[\varpi_\lambda])}(\calT_1,\calT_2)$.
\end{lemma}

\begin{proof}
	Fix trivializations $\iota_i: \calT_i \to \calG$, and consider $g=\iota_2\circ\lambda\circ \iota_1^{-1}$ as an element
	\[g\in {\on{H}}^0(\Yint{R^+}{[r,\infty]},\calG)=\{f\mid f:\Spec B^{R^+}_{[r,\infty]}\to \calG\}.\] 
	Since $\calG$ is affine and smooth over $\Spec \bbZ_p$ we may apply \Cref{lemma:random} to find a lift $g'\in \calG(\bbW(R^+))$ whose image in $\calG(B^{R^+}_{[m,\infty]}/[\varpi_\lambda])$ agrees with the image of $g$ for some $m\geq r$ and some pseudo-uniformizer $\varpi_\lambda$. 
	We conclude by letting $\tilde{\lambda}=\iota_2^{-1}\circ g'\circ \iota_1$ to get the desired isomorphism.
\end{proof}

Let us set some notation.
In what follows given a finite type affine scheme $X$ over $\Spec \bbZ_p$ and a topological ring $R$ over $\bbZ_p$, we topologize the set $X(R)$ by choosing an embedding $\iota:X\to \bbA^n_{\bbZ_p}$ and endowing $X(R)\subseteq R^n$ with the subspace topology.  
One can verify that the topology on $X(R)$ does not depend on $\iota$.

For the remainder of the section given $m\in \bbN$ we let $L^m_0$ denote the $p$-adic completion of $\bbW(R^+)[\frac{[\varpi]}{p^m}]$. 

\begin{remark}
	\label{randomcompletenesstatement}
Recall that $L^m_0$ is a ring of definition for $B^{R^+}_{[m,\infty]}$.
Note that since $[\varpi]=p^m\cdot (\frac{[\varpi]}{p^m})$ in $L^m_0$ and the ring is $p$-adically complete, then it is also $[\varpi]$-adically complete \cite[Lemme 1.4.14]{FF18}.
\end{remark}

\begin{lemma}
	\label{isogeny-andle}
	Let $m\in \bbN$ and pick two elements $\eta,\Phi\in \calG(B^{R^+}_{[m,\infty]})$.
	Suppose that $\eta$ restricts to $\on{Id}$ in $\calG(B^{R^+}_{[m,\infty]}/[\varpi])$.
	Consider the following sequence 
	\[\eta_{n+1}=\Phi\circ \varphi^*\eta_n \circ \Phi^{-1},\]
	starting with $\eta_0=\eta$.
	Then there is $m'\in \bbN$ with $m'\geq m$ sufficiently big such that the following hold
	\begin{enumerate}
		\item $\{\eta_n\}_{n\in \bbN}\subseteq \calG(L^{m'}_0) \subseteq \calG(B^{R^+}_{[m',\infty]})$.
		\item The sequence $\{\eta_n\}$ converges to the identity in $\calG(B^{R^+}_{[m',\infty]})$. 
	\end{enumerate}
\end{lemma}
\begin{proof}
	We fix a closed immersion $\iota:\calG\to \on{GL}_r$ for some $r$, which we may always find by \cite[Lemma 3.2]{Broshi}.
	Moreover, we regard $\on{GL}_r$ embedded as an open subset of the space of square matrices $M_{r\times r}$. 
	We will think of $\Phi$, $\Phi^{-1}$ and $\eta_n$ as elements of $M_{r\times r}(B^{R^+}_{[m,\infty]})$.
	Since $\eta=\eta_0$ restricts to $\on{Id}$ in $\calG(B^{R^+}_{[m,\infty]}/[\varpi])$ we can write it as 
	\[\eta_0=\on{Id}+[\varpi]\cdot\frac{1}{p^k}\cdot M_0 \] 
	where $M_0\in M_{r\times r}(L_0^m)$. 
	We may also write $\Phi=\frac{1}{p^\alpha} N_\alpha$ and $\Phi^{-1}=\frac{1}{p^\beta} N_\beta$ with $N_\alpha, N_\beta\in M_{r\times r}(L^m_0)$. 
	Then 
	\[\eta_1=\on{Id}+[\varpi^p]\cdot \frac{1}{p^{k+\alpha+\beta}}N_\alpha \varphi^*M_0 N_\beta.\]
	Observe that the term $N_\alpha\varphi^*M_0N_\beta\in M_{r\times r}(L^m_0)$.
Let $c=\alpha+\beta$.
	Choosing $m'$ large enough we can ensure that both $\frac{[\varpi]}{p^{k}}\in L^{m'}_0$ and $\frac{[\varpi]}{p^{c}}\in L^{m'}_0$ hold.
Consequently, $\eta_0, \eta_1\in \calG(L^{m'}_0)$.
By induction, we see that $\eta_n$ has the form
\[\eta_n=\on{Id}+[\varpi^{p^n}]\cdot \frac{1}{p^{k+n\cdot c}}N_\alpha \varphi^*M_{n-1} N_\beta=\on{Id}+[\varpi^{p^n-n-1}]\cdot \frac{[\varpi]}{p^k}\cdot  \frac{[\varpi^{n}]}{p^{n\cdot c}}\cdot M_n.\]
for some element $M_n\in M_{r\times r}(L^{m'}_0)$.	
This already shows that $\{\eta_n\}\subseteq \calG(L^{m'}_0)$.
Moreover, the same computation shows that for fixed $d\in \bbN$ there is a sufficiently large $n$ such that $\eta_n=\on{Id}$ on $\calG(L^{m'}_0/[\varpi^d])$.   
Since $L^{m'}_0$ is $[\varpi]$-adically complete the sequence $\eta_n$ converges to $\on{Id}$ in $\calG(L^{m'}_0)$.
Since $L^{m'}_0$ is a ring of definition of $B^{R^+}_{[m',\infty]}$, convergence in $\calG(L^{m'}_0)$ implies convergence in $\calG(B^{R^+}_{[m',\infty]})$.
\end{proof}

The proof of the following lemma is inspired by the computations that appear in \cite[Theorem 5.6]{HV11}, and it is a key input in the proof of \Cref{thm:comparetub}.
\begin{lemma}[Unique liftability of isogenies]
	\label{lem:liftabilityofiso}
	Let $\calT$ be a trivial $\calG$-torsor over $\on{Spec}(\bbW(R^+))$ and let $\calG_b$ denote the trivial $\calG$-torsor endowed with the $\varphi$-module structure over $\Yint{R^+}{(0,\infty]}$ given by an element $b\in \calG(B^{R^+}_{(0,\infty]})$. 
	Let $\Phi:\varphi^{*}\calT\to \calT$ be an isomorphism defined over $\on{Spec}(\bbW(R^+)[\frac{1}{\xi_{R^\sharp}}])$ and $\lambda:\calT\to \calG_b$ a $\varphi$-equivariant isomorphism defined over $B^{R^+}_{[m,\infty]}/[\varpi]$ for some $m\in \bbN$ sufficiently large so that $\xi_{R^\sharp}$ becomes a unit. 
	Then, there is $m'\in \bbN$ sufficiently large and a unique $\varphi$-equivariant isomorphism $\tilde{\lambda}:\calT\to \calG_b$ defined over $\Yint{R^+}{[m',\infty]}$ such that $\tilde{\lambda}$ and $\lambda$ restrict to the same map after basechange to $B^{R^+}_{[m',\infty]}/[\varpi]$.
\end{lemma}
\begin{proof}
	By transport of structure, we assume that $\calG=\calT$, that $\Phi\in \calG(\bbW(R^+)[\frac{1}{\xi_{R^\sharp}}])$, and that $\lambda\in \calG(B^{R^+}_{[m,\infty]}/[\varpi])$. 
	It suffices to find $m'\in \bbN$ and $\tilde{\lambda}\in \calG(B^{R^+}_{[m',\infty]})$ reducing to $\lambda$ in $\calG(B^{R^+}_{[m',\infty]}/[\varpi])$ and satisfying $\Phi=\tilde{\lambda}^{-1}\circ b\circ \varphi^*(\tilde{\lambda})$. 
	Choose an arbitrary lift $\lambda_1\in \calG(B^{R^+}_{[m,\infty]})$ of $\lambda$, and let $\eta_1={\lambda_1}^{-1}\circ b\circ \varphi^*({\lambda_1})\circ\Phi^{-1}$. 
	Observe that $\eta_1=\on{Id}$ in $\calG(B^{R^+}_{[m,\infty]}/[\varpi])$. 
	By \Cref{isogeny-andle} we may find $m'\in \bbN$ large enough so that $\eta_1$ and consequently $\eta_1^{-1}$ lie in $\calG(L^{m'}_0)\subseteq \calG(B^{R^+}_{[m',\infty]})$.

	We construct sequences of maps, $\{\lambda_i:\calG\to \calG_b\}_{i\in \bbN}$ defined over $\Spec B^{R^+}_{[m',\infty]}$ and $\{\eta_i:\calG\to \calG\}_{i\in \bbN}$ defined over $\Spec L^{m'}_0$ and given recursively by the relations 
	\[\lambda_{n+1}=\lambda_{n}\circ \eta_{n} \text{     and     } \eta_{n}=\lambda_{n}^{-1}\circ b\circ \varphi^*(\lambda_n)\circ \Phi^{-1}.\]
	Consider the following computation. 
\begin{align}
	\eta_{n+1} & = \lambda_{n+1}^{-1}\circ b \circ \varphi^*(\lambda_{n+1}) \circ \Phi^{-1} 
	= [\eta_{n}^{-1}\circ \lambda_n^{-1}]\circ b \circ \varphi^*(\lambda_{n+1}) \circ \Phi^{-1} \\
		   & = [\Phi \circ \varphi^*(\lambda_n)^{-1}\circ b^{-1}\circ \lambda_n] \circ \lambda_n^{-1}\circ b \circ \varphi^*(\lambda_{n+1}) \circ \Phi^{-1} \\
		   & = \Phi \circ \varphi^*(\lambda_n)^{-1}\circ\varphi^*(\lambda_{n+1}) \circ \Phi^{-1} 
	 = \Phi\circ\varphi^*(\eta_{n})\circ \Phi^{-1}
\end{align}
By \Cref{isogeny-andle} we may choose $m'$ so that $\{\eta_n\}\subseteq \calG(L^{m'}_0)$ and such that this sequence converges to $\on{Id}$.
This allows us to define $\tilde{\lambda}\in \calG(B^{R^+}_{[m',\infty]})$ as the limit of the $\lambda_i$. 
Taking limits we get the equation 
\[\on{Id}=\eta_\infty=\tilde{\lambda}\circ b\circ \varphi^*(\tilde{\lambda})\circ \Phi^{-1}\]
and we get the relation $\tilde{\lambda}=\lambda_i=\lambda$ in $\calG(B^{R^+}_{[m',\infty]}/[\varpi])$.
	
Let us prove uniqueness. 
Given two lifts $\tilde{\lambda}_i$ of $\lambda$ we let $g= \tilde{\lambda}_1\circ \tilde{\lambda}_2^{-1}$ with $g\in \calG(B^{R^+}_{[m',\infty]})$. 
Now, $\varphi$-equivariance gives $b=g^{-1}\circ b\circ \varphi^*(g)$, and since $g=\on{Id}$ in $\calG(B^{R^+}_{[m',\infty]}/[\varpi])$ then $\varphi^*(g)=\on{Id}$ in $\calG(B^{R^+}_{[m',\infty]}/[\varpi^p])$. 
From the identity $b=g^{-1}\circ b\circ \on{Id}$ in $\calG(B^{R^+}_{[m,\infty]}/[\varpi^p])$ and induction we can prove that $g=\on{Id}$ in $\calG(B^{R^+}_{[m',\infty]}/[\varpi^{p^n}])$ for every $n$. 
By \Cref{separatednessnewton} the identity $\on{Id}=g$ also holds in $\calG(B^{R^+}_{[m,\infty]})$. 
\end{proof}

\begin{lemma}
	\label{separatednessnewton}	
	If $\varpi\in R^+$ is a pseudo-uniformizer and $m\in \bbN$, then (as an abstract ring) $B^{R^+}_{[m,\infty]}$ is $[\varpi]$-adically separated.
\end{lemma}
\begin{proof}
	It suffices to show that if $f\in \bigcap_{i=1}^\infty [\varpi^i]\cdot B^{R^+}_{[m,\infty]}$ then $f=0$.	
	Alternatively, since $\calY^{R^+}_{[m,\infty]}$ is sous-perfectoid \cite[Proposition 13.1.1]{SW20} it suffices to show that for all points $x\in \calY^{R^+}_{[m,\infty]}$ the value $|f|_x=0$ \cite[Theorem 5.2.1]{SW20}.
	This condition can be checked on rank $1$ geometric points $\Spa(C,O_C)\to \Spa(R,R^+)$, so without loss of generality $(R,R^+)=(C,O_C)$.

	In this case $B^{O_C}_{[m,\infty]}=B_I^{+}$ in the notation of \cite[D\'efinition 1.10.2]{FF18} (this definition depends on \cite[D\'efinition 1.4.1, 1.3.2, 1.3.1]{FF18}). 
	Here $I=[\rho,1]$ for some number $\rho\in (0,1)$ associated to $m$, whose precise formula is irrelevant. 
	Indeed, this follows from \cite[Exemple 1.10.3, Proposition 1.10.5]{FF18}.
	We wish to show that $f=0$ assuming that $f\in \bigcap_{i=1}^\infty [\varpi^i]\cdot B^{+}_{I}$.
	Write $f=f_n\cdot [\varpi^n]$.

	Recall the ring $B_I$ \cite[D\'efinition 1.6.2]{FF18} (see \cite[D\'efinition 1.4.1, 1.3.2, 1.3.1]{FF18}).
	Recall that the inclusion of intervals $\{1\}\subset I$ induces a continuous map $B_I\to B_{\{1\}}$ (see \cite[$\mathsection$ 1.6.1]{FF18}). 
	Moreover, recall from \cite[\S 1.10.1]{FF18} that $B^+_I$ (respectively $B^+_{\{1\}})$) is the closure of $B^{b,+}$ \cite[D\'efinition 1.3.2]{FF18} in $B_I$ (respectively in $B_{\{1\}}$).
Overall, we get a commutative diagram of continuous ring homomorphisms 
\begin{center}
\begin{tikzcd}
	B^+_I\arrow{r} \arrow{d}  & B^+_{\{1\}} \arrow{d} \\
	B_I \arrow{r} & B_{\{1\}}.
\end{tikzcd}
\end{center}
Also, the map $B_I\to B_{\{1\}}$ is injective by \cite[Proposition 1.6.15]{FF18}, so in particular the map $B^+_I\to B^+_{\{1\}}$ is also injective.
It suffices to show that $B^+_{\{1\}}$ is $[\varpi]$-adically separated. 
Now $B_{\{1\}}$ is obtained as the completion of $B^b$ \cite[D\'efinition 1.3.2]{FF18} under the norm $|\cdot|_1$ as in \cite[D\'efinition 1.4.1]{FF18}. 
One can verify directly that for all $g\in B^{b,+}$ $|g|_1\leq 1$ holds. 
By continuity, we can conclude that any element $g\in B^+_{\{1\}}$ satisfies that $|g|_1\leq 1$.
In particular, $|f_n|_1\leq 1$ for all $f_n$ as above.
Since $\varpi$ is a pseudo-uniformizer in $C$, then 
\[|[\varpi]|_1=|\varpi|_C=\epsilon \text{ for some } 0<\epsilon<1,\] we can conclude that $|f|_1=|f_n|_1|[\varpi^n]|_1\leq \epsilon^n$ for all $n$. 
This allows us to deduce that $f=0$.
\end{proof}

	\subsection{Specialization maps.}
	\label{Sepcializationsections}
	In this section we recall the theory of specialization maps for v-sheaves as developed in \cite{Gle22}. 
	The theory is an attempt to answer the question: What is the analogue of a formal scheme in the context of v-sheaves?
	This theory approaches the question in steps. 
	\begin{align*}
		\{\on{Small\, v-sheaves}\}& \supseteq \{\on{Specializing\,v-sheaves}\} \\
						\{\on{Specializing\,v-sheaves}\} & \supseteq \{\on{Prekimberlites}\} \\
						\{\on{Prekimberlites}\} & \supseteq \{\on{Valuative\,prekimberlites}\} \\
						\{\on{Valuative\,prekimberlites}\} & \supseteq \{\on{Kimberlites}\} \\
						\{\on{Kimberlites}\} & \supseteq \{\on{Locally\,spatial\,kimberlites}\} 
	\end{align*}
	Each of these categories is a full subcategory of the category of small v-sheaves \cite[Definition 12.1]{Sch17} obtained by adding axioms at each stage. 
	For our purposes it will suffice to discuss valuative prekimberlites, but we will also mention a category that vaguely speaking sits between the category of valuative prekimberlites and the category of kimberlites.  
	We will make this precise in what follows.
	\subsubsection{Specializing v-sheaves}
	Let $\on{CAlg}^{\on{perf},\op}_{\bbF_p}$ denote the category of perfect affine schemes in characteristic $p$.
	Recall from \cite[$\mathsection$ 18.3]{SW20} the $\diamond$-functor 
	\[\diamond:\on{CAlg}^{\on{perf},\op}_{\bbF_p}\to \widetilde{\Perf}\]
	with
	\[(\Spec A)^\diamond := \Spd(A,A).\] 
	In other words, it attaches to a perfect affine scheme over $\bbF_p$ the small v-sheaf associated to the non-analytic adic space $\Spa(A,A)$ where $A$ is regarded as a topological ring endowed with the discrete topology. 
	Recall that $\diamond$ extends to a fully faithful functor $\diamond:\PSch\to \widetilde{\Perf}$ \cite[Proposition 18.3.1]{SW20} \cite[$\mathsection 3$]{Gle22} from the category $\PSch$ of perfect schemes in characteristic $p$ to the category of small v-sheaves. 
	Furthermore, this functor formally extends to a functor $\diamond:\widetilde{\PSch}\to \widetilde{\Perf}$ from the category of small scheme-theoretic v-sheaves to the category of small v-sheaves. 
Furthermore, $\diamond$ admits a right adjoint functor \cite[Definition 3.12]{Gle22} which we call the reduction functor
\begin{equation}
	\label{equationreductionfunctor}
	\red:\widetilde{\Perf}\to \widetilde{\PSch}.
\end{equation}
	We let $X^\Red=(X^\red)^\diamond$, it comes with a canonical map coming from adjunction $X^\Red\to X$.
	\begin{definition}{\cite[Definition 3.20]{Gle22}}
		\label{formallyadicdefi}
		A map $Y\to X$ is \textit{formally adic} if the following diagram is Cartesian. 
		\begin{center}
		\begin{tikzcd}
		 Y^\Red \arrow{r} \arrow{d}  & \arrow{d} Y \\
		 X^\Red \arrow{r} & X. 
		\end{tikzcd}
		\end{center}
	\end{definition}

	Given maps $X\xrightarrow{f} Y \xrightarrow{g} Z$ with $h=g\circ f$, one can use the cancellation and composition properties of Cartesian diagrams to verify that if $g$ is formally adic then $f$ is formally adic if and only if $h$ is formally adic. 

\begin{definition}{\cite[Definition 3.27]{Gle22}}
	\label{formallyseparateddefi}
		A v-sheaf is \textit{formally separated} if the diagonal is a closed immersion and formally adic.	
	\end{definition}

	Recall that to any Huber pair $(A,A^+)$ over $\bbZ_p$ we can attach a v-sheaf $\Spd\Hub{A}$ \Cref{definitionspdA}.

	\begin{definition}{\cite[Definition 4.6]{Gle22}}
		\label{defi:formalizing}
		Given a v-sheaf $X$ and a map $f:\Spa \Hub{R}\to X$ we say that $X$ \textit{formalizes} $f$ if there exists a dashed arrow completing the commutative diagram below.
		\begin{center}
		\begin{tikzcd}
			\Spa (R,R^+)\arrow{r}{f} \arrow{d}  & X  \\
			\Spd (R^+,R^+) \arrow[dashed]{ru}  
		\end{tikzcd}
		\end{center}
		Any such arrow is called a \textit{formalization} of $f$.
		We say $X$ is \textit{v-formalizing} if for any $f$ as above there is a v-cover $g:\Spa (R',R'^+)\to \Spa (R,R^+)$ such that $X$ formalizes $g\circ f$. 
	\end{definition}

	Whenever $X$ is formally separated and $f:\Spa (R,R^+)\to X$ is a map, a formalization of $f$, if it exists, is unique \cite[Proposition 4.9]{Gle22}.
This leads to the first class of v-sheaves that admit a specialization map.

	\begin{definition}{\cite[Definition 4.11]{Gle22}}
		\label{defi:specializingvsheaf}
		We say that a small v-sheaf $X$ is \textit{specializing} if it is formally separated and v-formalizing.	
	\end{definition}

As shown in \cite[Proposition 4.14]{Gle22} if $X$ is a specializing v-sheaf then it has a continuous specialization map 
\[\on{sp}_X:|X|\to |X^\red|.\]

\begin{proposition}
	\label{functoriality-of-special}
The rule 
\[X\mapsto (\on{sp}_X:|X|\to |X^\red|)\]
with source the category of specializing v-sheaves and target the category of maps of topological spaces is a functor.
In other words, given a map of specializing v-sheaves $f:X\to Y$ we get a commutative diagram of topological spaces.
\begin{center}
\begin{tikzcd}
	\mid X\mid \arrow{r}{\mid f\mid} \arrow{d}{{\on{sp}}_X}  & \mid Y\mid \arrow{d}{{\on{sp}}_Y} \\
 \mid X^\red\mid \arrow{r}{\mid f^\red \mid} & \mid Y^\red\mid
\end{tikzcd}
\end{center}
\end{proposition}
\begin{proof}
This is the content of \cite[Proposition 4.14]{Gle22}.	
\end{proof}

	\begin{definition}{\cite[Definition 4.15]{Gle22}}
		\label{thisisaprekimberlitehurray}
		If $X$ is a specializing v-sheaf we say that it is a \textit{prekimberlite} if $X^\red\in \widetilde{\PSch}$ is represented by a perfect scheme and $X^\Red\to X$ is a closed immersion.
		We let $X^{\on{an}}:=X\setminus X^\Red$ and we call this the \textit{analytic locus} of $X$.
	\end{definition}

The following lemma shows that, under certain assumption, being a prekimberlite can be verified Zariski locally. 

\begin{lemma}
	\label{lemmaspecializing}
	Let $X$ be a specializing v-sheaf and let $Y=X^\red$. 
	Suppose that $Y$ is representable by a perfect scheme. 
	Let $U\subseteq Y$ be an open subset and let $V$ denote the unique open subsheaf of $X$ with $|V|=\on{sp}^{-1}(|U|)$. 
	Then the following hold.
	\begin{enumerate}
		\item $V$ is a specializing v-sheaf with $V^\red=U$.
		\item The map $V\to X$ is formally adic.
	\end{enumerate}
	Moreover, if there is an open cover $\{U_i\to Y\}_{i\in I}$ such that $V_i:=\on{sp}_X^{-1}(U_i)$ is a prekimberlite, then $X$ is a prekimberlite. 
\end{lemma}
\begin{proof}
	Let us show $V$ is v-formalizing. 
	Fix $f:\Spa (R,R^+)\to X$ and assume that $f$ formalizes to a map $\Spd(R^+,R^+)\to X$.
	It follows easily from the definition of the specialization map that $f$ factors through $V$ if and only if the map induced by reduction  \[\Spec(R^+/R^{\circ \circ})=\Spec R^+_\red=\Spd(R^+,R^+)^\red\to Y\] 
	factors through $U$ \cite[Proposition 3.18, Definition 4.12]{Gle22}. 
	In this case, $\Spd(R^+,R^+)^\Red\to X$ factors as a composition 
	\[\Spd(R^+_\red,R^+_\red)=\Spd(R^+,R^+)^\Red\to U^\diamond\subseteq V\subseteq X.\]
	Since $|\Spd(R^+,R^+)|=|\Spa\Hub{R}|\cup |\Spd(R^+,R^+)^\Red|$ and $V$ is an open subsheaf, we deduce that if $f$ factors through $V$ then $\Spd(R^+,R^+)\to X$ also factors through $V$. 
	Consequently, $V$ is v-formalizing. 

	Any subsheaf of a formally separated v-sheaf is again formally separated.
	Indeed, this follows by \cite[Lemma 3.30]{Gle22} using that both $(-)^\red$ and $(-)^\diamond$ commute with finite limits. 
	This shows $V$ is formally separated, and hence specializing.

	We claim that $Y^\diamond\cap_X V=U^\diamond$.
	Any map $\Spa(R,R^+)\to Y^\diamond$ is v-locally induced by a map of schemes $\Spec R^+\to Y$.
	If it factors through $V$ it is because the map $\Spec (R^+/R^{\circ \circ})\to Y$ factors through $U$. 
	But since $U$ is open and $\Spec R^+/R^{\circ \circ}$ contains all closed points of $\Spec R^+$, the map $\Spec R^+\to Y$ also factors through $U$ as we needed to show.
	By \cite[Lemma 3.32]{Gle22}, $V^\red=U$ and the map $V\to X$ is formally adic. 

	Let us show the final statement. 
	Since we already assumed that $Y$ is representable, it suffices to show that $Y^\diamond\to X$ is a closed immersion.
	This can be checked v-locally on $X$ by \cite[Proposition 10.11.(i)]{Sch17}. By our argument above, we have a Cartesian diagram
	\begin{center}
	\begin{tikzcd}
		 \coprod_{i\in I} U_i^\diamond \arrow{r} \arrow{d}  & \coprod_{i\in I} V_i \arrow{d} \\
	 Y^\diamond \arrow{r} & X
	\end{tikzcd}
	\end{center}
	and since we assumed each $V_i$ is a prekimberlite the map $\coprod_{i\in I} U_i^\diamond \to \coprod_{i\in I} V_i$ is a closed immersion as we wanted to show. 
\end{proof}

\subsubsection{Prekimberlites and valuative prekimberlites.}
	For prekimberlites one can construct a v-sheaf theoretic specialization map (or Heuer specialization map). 
	Recall our notation $R^+_\red=R^+/R^{\circ \circ}$.
Given $S\in \PSch$ one can construct a v-sheaf $S^{\diamond/\circ\circ}$ \cite[Definition 4.23]{Gle22}\footnote{The symbol $\diamond/\circ\circ$ suggest the similarity with the $\diamond$ functor up to a quotient by the space of topological nilpotent elements. 
In \cite{Gle22} we initially took the minimalistic notation $\diamond/\circ$ instead of $\diamond/\circ\circ$.}, \cite[Definition 5.1]{Heu21} by v-sheafifying the formula
\[S^{\diamond/\circ\circ_{\on{pre}}}:\Spa \Hub{R}\mapsto S(\Spec R^+_\red).\]
As it turns out it is only necessary to sheafify for the analytic topology \cite[Lemma 5.2]{Heu21}. 
When $X$ is a prekimberlite we let $X^{\on{H}}:=(X^\red)^{\diamond/\circ\circ}$. 
We get a v-sheaf theoretic specialization map \cite[$\mathsection 4.4$]{Gle22} 
\[\on{SP}:X\to X^{\on{H}}.\]
\begin{definition}
The v-sheaf theoretic specialization map is constructed as follows. If $\alpha\in X(R,R^+)$ and $\alpha$ is formalizable we let $\tilde{\alpha}\in X(\Spd R^+)$ be its unique formalization.
Applying the reduction functor gives $\tilde{\alpha}^\red\in X^\red(\Spec R^+_\red)$, which is an element in the presheaf $(X^{\red})^{\diamond/\circ\circ_{\on{pre}}}(R,R^+)$.
We get 
\[\on{SP}_{\on{pre}}:X^{\on{frml}}\to (X^\red)^{\diamond/\circ\circ_{\on{pre}}},\] 
where the source is the sub-presheaf of formalizable maps in $X$. Now, $\on{SP}$ is the sheafification of $\on{SP}_{\on{pre}}$.
\end{definition}

The v-sheaf theoretic specialization map allow us to make interesting constructions. 
If $X$ is a prekimberlite and $Z\subseteq X^\red$ is a locally closed subset we can define a prekimberlite $\widehat{X}_{/Z}$ \cite[Proposition 4.21]{Gle22} such that $Z=(\widehat{X}_{/Z})^\red$ and such that it fits in the following Cartesian diagram 
\begin{equation}
	\label{formalnbhoods}
\begin{tikzcd}
	\widehat{X}_{/Z} \arrow{r} \arrow{d}{\on{SP}}  & X \arrow{d}{\on{SP}} \\
 Z^{\diamond/\circ\circ} \arrow{r} & X^{\on{H}}. 
\end{tikzcd}
\end{equation}
When $Z$ is constructible the map $\widehat{X}_{/Z}\to X$ is an open immersion \cite[Proposition 4.22]{Gle22}.
\begin{definition}{\cite[Definition 4.18]{Gle22}}
	\label{formalnbbhoosoughtotbedefined}
	The prekimberlite $\widehat{X}_{/Z}$ obtained from diagram \eqref{formalnbhoods} is called the \textit{formal neighborhood} of $X$ along $Z$. 
\end{definition}
\begin{remark}
	We will mostly apply this construction to the case where $Z\to X$ is the inclusion of a closed point. 
	In this case $\widehat{X}_{/Z}$ is the v-sheaf theoretic analogue of taking the formal completion at the closed point. 
\end{remark}

The specialization map for prekimberlites, $\on{SP}_X:X\to X^{\on{H}}$, is always separated. 
Indeed, the inclusion $X\times_{X^{\on{H}}}X\subseteq X\times X$ is a separated map and since $X$ is formally separated $\Delta:X\to X\times X$ is a closed immersion. 
Consequently, $X\to X\times_{X^{\on{H}}}X$ also is.
\begin{definition}{\cite[Definition 4.30]{Gle22}}
	We say that a prekimberlite $X$ is \textit{valuative} if $\on{SP}:X\to X^{\on{H}}$ is partially proper (\cite[Definition 18.4]{Sch17}).
\end{definition}

For valuative prekimberlites the topological specialization map, $\on{sp}_X:|X|\to |X^\red|$, is specializing \cite[Proposition 4.33]{Gle22}.
Moreover, the valuative property is stable under natural constructions like taking formal neighborhoods or \'etale formal neighborhoods \cite[Proposition 4.34]{Gle22}. 

Although the specialization map for a specializing v-sheaf $\calF$ is defined as a map 
\[\on{sp}_\calF:|\calF|\to |\calF^\red|,\]
one often wishes to study this map after restricting it to certain open subsets. 
For example, if $\calF$ is a prekimberlite it is natural to study 
\[\on{sp}_\calF:|\calF^{\on{an}}|\to |\calF^\red|\]
instead. 

Since we are interested in studying the $p$-adic generic fiber of moduli spaces of $p$-adic shtukas the natural setup is to consider prekimberlites $\calF$ together with a map $f:\calF\to \Spd \bbZ_p$ and study the specialization map 
\[\on{sp}_\calF:|\calF\times_{\Spd \bbZ_p}\Spd \bbQ_p|\to |\calF^\red|.\]
Note that in general the map $f$ might not be formally adic (i.e. $\calF$ is not a $p$-adic prekimberlite), but we always have that 
\[\calF\times_{\Spd \bbZ_p}\Spd \bbQ_p\subseteq \calF^{\on{an}}.\]
This motivates the following definition.
\begin{definition}
	\label{definitionkimberlite}
	Let $\calF$ be a valuative prekimberlite.
	\begin{enumerate}
		\item A \textit{smelted kimberlite} is a pair $\calK=(\calF,\calD)$ such that $\calD$ is a quasiseparated locally spatial diamond and $\calD \subseteq \calF^{\on{an}}$ is open.\footnote{The main cases of interest are when $\calD=\calF^{\on{an}}$ or when $\calD=\calF\times_{\Spd\bbZ_p} \Spd \bbQ_p$.} 
\item The specialization map $\on{sp}_{\mathcal{K}}:|\calD|\to |\calF^\red|$ is defined as the composition $|\calD|\subseteq |\calF|\xrightarrow{\on{sp}_{\calF}} |\calF^\red|$.
	When the context is clear we may write $\on{sp}_\calD$ or merely $\on{sp}$ instead of $\on{sp}_\calK$.
\item We say that $\calF$ is a \textit{kimberlite} if $(\calF,\calF^{\on{an}})$ is a smelted kimberlite and $\on{sp}_{\calF^{\on{an}}}$ is quasicompact.
\item Given a smelted kimberlite $\calK=(\calF,\calD)$ and a locally closed subset $Z\subseteq \calF^\red$ we define the \textit{tubular neighborhood} of $\calK$ along $Z$ as $\calK^{\circledcirc}_{/Z}:=\Tf{\calF}{Z}\cap \calD$. 
	\end{enumerate}
\end{definition}

Recall the definition of locally spectral spaces and spectral maps between them \cite{MR251026}, \cite[Definition 2.1]{Sch17}. 
We have the following key result which we will use later on \cite[Theorem 9, Theorem 4.40]{Gle22}. 

\begin{theorem}
	\label{pro2:spectralmapSch-Spatial}
	Let $\calK=(\calF,\calD)$ be a smelted kimberlite, then 
	\[\on{sp}_\calD:|\calD|\to |\calF^\red|\] is a specializing, spectral map of locally spectral spaces. 
\end{theorem}

\subsection{The specialization map for the $p$-adic Beilinson--Drinfeld Grassmannian.}
\label{Grassmanianssection}
We recall the definition of the $p$-adic Beilinson--Drinfeld Grassmannian that is the most suitable to study its specialization map.
\begin{definition}\textup{(\cite[Definition 20.3.1]{SW20})}
	\label{defi:Grassmannian}
	We let $\Gr_\calG$ denote the v-sheaf 
	\[\Gr_\calG:\Perf^\op\to \Sets\]
	that assigns to an affinoid perfectoid pair $(R,R^+)$ the set
	\[\Gr_\calG (R,R^+)=\{((R^\sharp,\iota),\calT,\psi)\}_{/\simeq}\] 
	where $(R^\sharp,\iota)$ an untilt of $R$ and $(\calT,\psi)$ is a lattice with $\calG$-structure as in \Cref{defi:latticesGstructure}.
\end{definition}

\begin{proposition}
	\label{pro:grassmanianvformal}
	Let $\Spa(R,R^+)$ be a product of points as in \Cref{defi:prodpoints} and let $f:\Spa(R,R^+)\to \Gr_\calG$ be a map.  
	Then $\Gr_\calG$ formalizes $f$ (\Cref{defi:formalizing}). 
	In particular, $f$ is v-formalizing.
\end{proposition}
\begin{proof}
Let $\Spa(R,R^+)$ be a product of points and $f:\Spa(R,R^+)\to \Gr_\calG$ a map. 
By definition, associated to this map we have an untilt $R^\sharp$ and a $\calG$-torsor $\calT$ over $\Yinf{R^+}$ together with a trivialization $\psi:\calT\to \calG$ over $\Yinf{R^+}\setminus V(\xi_{R^\sharp})$ meromorphic along $\xi_{R^\sharp}$. 
We use $\psi$ to glue $\calT$ and $\calG$ along $\Yint{R^+}{[r,\infty)}$ to get a $\calG$-torsor defined over $\calY_{[0,\infty]}^{R^+}$. 
Using \Cref{cor:gagatypekedlymeoro} and \Cref{pro:prodpointextendg-tors} 
we can extend the data $(\calT,\psi)$ to a get a $\calG$-torsor defined over $\on{Spec}(\bbW(R^+))$ together with a trivialization defined over $\on{Spec}(\bbW(R^+)[\frac{1}{\xi_{R^\sharp}}])$. 
	This is enough to define a map $\Spd(R^+,R^+)\to \Gr_\calG$ that restricts to the original one. 
	Indeed, take a second affinoid perfectoid $\Spa(T,T^+)$ and a map $g:\Spa(T,T^+)\to \Spd(R^+,R^+)$, we want to produce a map $\Spa{(T,T^+)}\to \Gr_\calG$ in a functorial way. 
	We may construct an untilt $T^\sharp$ by letting $\xi_{T^\sharp}$ denote the image of $\xi_{R^\sharp}$ under the ring map $g':\bbW(R^+)\to \bbW(T^+)$ induced by $g$.
	Base change along $g'$ gives a $\calG$-torsor over $\on{Spec}(\bbW(T^+))$ together with a trivialization over $\on{Spec}(\bbW(T^+)[\frac{1}{g'(\xi_{R^\sharp})}])$. 
	This restricts to a $\calG$-torsor over $\Yinf{T^+}$ and a trivialization over $\Yinf{T^+}\setminus V(g'(\xi_{R^\sharp}))$ that is meromorphic along $g'(\xi_{R^\sharp})$. 
	This gives our desired natural transformation $\Spd{(R^+,R^+)}\to \Gr_\calG$. 
	Clearly the composition $\Spa(R,R^+)\to\Spd{(R^+,R^+)}\to \Gr_\calG$ agrees with $f$, so this map is a formalization. 
\end{proof}

Recall that associated to our parahoric group scheme $\calG$ one can construct a Witt vector flag variety
\[\calF\ell_{\calG,\bbW}:\PSch^\op \to \Sets\]
with formula 
\[\Spec R\mapsto \{(\calT,\psi)\}_{/\simeq}\] 
where $\calT$ is a $\calG$-torsors over $\Spec \bbW(R)$ and $\psi:\calT\to \calG$ is a trivialization over $\Spec \bbW(R)[\frac{1}{p}]$ \cite{Zhu17}, \cite[Definition 9.4]{BS17}. 

\begin{proposition}\textup{(\cite[\S 20.3]{SW20})}
	\label{pro:grassmanianreduced}
	The v-sheaf $\Gr_\calG$ is specializing, the map $\Gr_\calG\to \Spd {\bbZ}_p$ is formally adic and $\Gr^\red_\calG$ is represented by $\calF\ell_{\calG,\bbW}$.
\end{proposition}
\begin{proof}
	To show that $\Gr_\calG$ is specializing we have to show that it is v-formalizing and formally separated. 
	It is v-formalizing by \Cref{pro:grassmanianvformal}.
	In turn, it is separated by \cite[Theorem 20.3.4, Theorem 21.2.1]{SW20}. 
	By \cite[Proposition 3.29]{Gle22}, to show that it is formally separated it will suffice to show it is formally adic over $\Spd \bbZ_p$.
	Let us now show that $\Gr_\calG\to \Spd \bbZ_p$ is formally adic.
	By \cite[Lemma 3.32]{Gle22} it suffices to show that 
	\[\Gr_\calG\times_{\Spd \bbZ_p} \on{Spec}(\bbF_p)^\diamond=(\GrW)^\diamond.\] 
	Indeed, by \cite{BS17} $\GrW$ is ind-representable and by \cite[Proposition 3.16.(1)]{Gle22} combined with a quasi-compactness argument one can show that ind-representable scheme-theoretic v-sheaves are reduced in the sense of \cite[Definition 3.15]{Gle22}.

To find an identification $\Gr_\calG\times_{\Spd \bbZ_p} \on{Spec}(\bbF_p)^\diamond=(\GrW)^\diamond$ we begin by constructing a map $\GrW\to (\Gr_\calG)^\red$. 
We need to produce a map $\on{Spec}(R)^\dia\to \Gr_\calG$ functorially on the comma category $\PSch_{/\GrW}$. 
An object in $\PSch_{/\GrW}$ is given by an affine scheme $\Spec R$, a $\calG$-torsor $\calT$ over $\on{Spec}(\bbW(R))$ together with a trivialization $\psi:\calT\to \calG$ over $\on{Spec}(\bbW(R)[\frac{1}{p}])$. 
Given an affinoid perfectoid $\Spa(T,T^+)$ and a map $f:\Spa(T,T^+)\to \on{Spec}(R)^\dia$ we need to produce a map $\Spa(T,T^+)\to \Gr_\calG$. 
The morphism $f$ induces the ring map $f':\bbW(R)\to \bbW(T^+)$. 
We can assign to $f$ the characteristic $p$ untilt and assign the $\calG$-bundle $f'^*\calT$ over $\Yinf{T^+}$ with trivialization $f'^*\psi$, and using \Cref{cor:gagatypekedlymeoro} we see that it is meromorphic along $p$. 
This construction is clearly functorial and gives the desired map. 
We prove that for any $(R,R^+)$ we have bijection of sets: 
\[(\GrW)^\diamond(R,R^+)\to \Gr_\calG\times_{\Spd \bbZ_p} \Spd \bbF_p (R,R^+).\] 
To prove injectivity, suppose we are given two maps $g_i:\Spa(R,R^+)\to (\GrW)^\diamond$ in characteristic $p$ whose composition agree. 
It is enough to prove that $g_1=g_2$ after taking a v-cover of $\Spa(R,R^+)$. 
Locally for the v-topology we can assume that both maps factor through morphisms $g_i':\on{Spec}(R^+)\to \GrW$ given by pairs $(\calT_i,\psi_i)$. 
Since the compositions agree, these pairs become isomorphic over $\Yinf{R^+}$, and arguing as in the proof of \Cref{pro:grassmanianreduced} we can conclude that this pairs are already isomorphic over $Y^{R^+}_{[0,\infty]}$. 
Since both $\calT_i$ are defined over $\on{Spec}(\bbW(R^+))$ and the pullback functor $j^*$ of \Cref{thm:puncturedtoAinfked} is fully faithful we can conclude that $g_1'=g_2'$.
	
To prove surjectivity take a map $f:\Spa(R,R^+)\to \Gr_\calG\times_{\Spd \bbZ_p}\on{Spec}(\bbF_p)^\diamond$. 
Since surjectivity can be checked v-locally we can assume that $\Spa(R,R^+)$ is a product of points. 
By the proof of \Cref{pro:grassmanianvformal} we get a $\calG$-torsor $\calT$ over $\on{Spec}(\bbW(R^+))$ and a trivialization over $\on{Spec}(\bbW(R^+)[\frac{1}{p}])$ which gives a map $\on{Spec}(R^+)\to \GrW$ and consequently the required lift to our original map $\Spa(R,R^+)\to (\GrW)^\dia$.
\end{proof}
Let $T\subseteq G_{{\bbQ_p}}$ be a maximal torus and let $T_{\bar{\bbQ}_p}\subseteq B\subseteq G_{\bar{\bbQ}_p}$ be a choice of Borel of the geometric generic fiber of $G$.
Fix a dominant geometric cocharacter of $T$, $\mu\in X^+_*(T)$, with reflex field $E$ and ring of integers $O_E$. 
Recall that to $\mu$ we may attach a ``local model v-sheaf'' $\Grm{O_{E}}$ over $\Spd {O_{E}}$ \cite[Definition 4.11]{AGLR22}, \cite[\S 21.4]{SW20}. 
\begin{definition}
	\label{thelocalmodeldefi}
	The \textit{local model} $\Grm{O_{E}}$ is defined as the v-sheaf closure \cite[Definition 2.3]{AGLR22} of $\on{Gr}^{G,\leq \mu}_{E}$ in $\Gr_\calG\times_{\Spd \bbZ_p} \Spd O_{E}$. 
\end{definition}

\begin{definition}
		Let $X$ be a kimberlite over $\Spd \bbZ_p$ we say it is \textit{flat} if there is a set $I$, a family of perfectoid Huber pairs $\{(R_i^\sharp,R_i^{\sharp+})\}_{i \in I}$ over $\bbQ_p$ and a v-cover over $\Spd \bbZ_p$ 
		\begin{equation}\coprod_{i\in I} \Spd(R_i^{\sharp,+})\to X\end{equation}	
	\end{definition}

Recall $\GrWmfe\subseteq \GrW$ the $\mu$-admissible locus \cite{KR00}, \cite[Definition 3.11]{AGLR22}. 
In our collaboration with Ansch\"utz, Louren\c{c}o, and Richarz \cite{AGLR22} we prove the following statement.

\begin{theorem}\textup{(\cite[Proposition 4.14, Theorem 6.16]{AGLR22})}
	\label{thm:AGLR}
	If $\calG$ is parahoric and $\mu\in X^+_*(T)$, then $\Grm{O_{E}}$ is a flat $p$-adic kimberlite with $(\Grm{O_{E}})^{\on{red}}=\GrWmfe$. 
\end{theorem}

\subsection{Specialization maps for moduli spaces of $p$-adic shtukas.}
\label{section-moduli-ofshtukas}
Fix an element $b\in G(\breve{\bbQ}_p)$ and let \[V_b:\on{Rep}_{G}\to \Isoc_{\bar{\bbF}_p}\] denote the associated isocrystal with $G$-structure (as in \S \ref{subsection-lattices-and-padicshtukas}). 
This induces a $\varphi$-module over $\calY_{(0,\infty)}$ which we denote by $\calG_b$.
\begin{definition}
	\label{defi:moduliofshtukas}
	The \textit{integral moduli space of $p$-adic $\calG$-shtukas} associated to $V_b$, which we denote by $\Sht_\calG(b)$, is the functor 
	\[\Sht_\calG(b):\Perff^\op\to \Sets\]
	\[(R,R^+)\mapsto \{((R^\sharp,\iota),\calT,\Phi,\lambda)\}/_{\simeq}\]
	where $(R^\sharp,\iota)$ is an untilt of $R$, $(\calT,\Phi)$ is a shtuka as in \Cref{defi:singleshtukasGstructure} and 
	$\lambda:\calT\to \calG_b|_{\Yint{R^+}{[r,\infty)}}$ an isogeny as in \Cref{defi:isogenies}. 
\end{definition}
Fix $\Hub{R}$, fix $R^\sharp$ an untilt of $R$ and fix $M\in \calG(\bbW(R^+)[\frac{1}{\xi_{R^\sharp}}])$ to such data we associate a shtuka 
\[\calG_M:=(\calG,\Phi_M)\]
where \[\Phi_M:\varphi^*\calG\to \calG\]
is the only isomorphism conjugate to $M:\calG\to \calG$ under the canonical isomorphism $\varphi^*\calG\simeq \calG$. 
We consider the following auxiliary space.
\begin{definition}
\label{defi:auxshtuka}
Let $\bbW\Sht(b)$ denote the functor 
	\[\bbW\Sht_\calG(b):\Perff^\op\to \Sets\]
	\[{(R,R^+)}\mapsto \{((R^\sharp,\iota),M,\lambda)\}\] 
	where $(R^\sharp,\iota)$ is an untilt of $R$, $M\in \calG(\bbW(R^+)[\frac{1}{\xi_{R^\sharp}}])$ and $\lambda:\calG_M\to \calG_b$ an isogeny. 
Here $\calG_M:=(\calG,\Phi_M)$ as above. 
\end{definition}
We denote by $\mathbb{W}^+\calG$ the sheaf in groups 
\[\mathbb{W}^+\calG(R,R^+)=\calG(\bbW(R^+)).\] 
\begin{proposition}
	\label{pro:shtukavformalizing}
	Consider the map $\bbW\Sht_\calG(b)\to \Shtdos$ given by 
	\[(R^\sharp, M,\lambda)\mapsto (R^\sharp, \calG_M, \Phi_M, \lambda)\]
	The following statements hold.
\begin{enumerate}
	\item The map $\bbW\Sht_\calG(b)\to \Shtdos$ is a $\mathbb{W}^+\calG$-torsor for the v-topology.
	\item $\Auxsht$ is formalizing and $\Shtdos$ is v-formalizing as in \Cref{defi:formalizing}. 
\end{enumerate}
\end{proposition}
\begin{proof}
	Given $N\in \mathbb{W}^+\calG(R,R^+)$ and $(R^\sharp,M,\lambda)\in \Auxsht\Hub{R}$ let 
	\[(R^\sharp,M,\lambda)\star N=(R^\sharp,N^{-1}M\varphi(N),\lambda\circ N).\] 
	This action on $\Auxsht$ makes the map $\Auxsht\to \Shtdos$ equivariant for the trivial $\bbW^+\calG$-action on the target. 
	It suffices to show that the map is v-locally isomorphic to the trivial $\mathbb{W}^+\calG$-torsor.
	Since products of points are a basis for the v-topology it suffices to understand the base changes $\Auxsht\times_\Shtdos \Spa\Hub{R}$ as $\Spa\Hub{R}$ ranges over products of points.
	Let $\Spa\Hub{R}$ be a product of points, and let $(R^\sharp,\calT,\Phi,\lambda)\in \Shtdos\Hub{R}$. 
	Similarly to the proof of \Cref{pro:grassmanianvformal}, we can glue $\calT$ along $\lambda$ over $\Yint{R^+}{[r,\infty)}$ and use \Cref{pro:prodpointextendg-tors} to get (uniquely) a $\calG$-bundle $\calT_\bbW$ over $\on{Spec}(\bbW(R^+))$ with a meromorphic $\Phi_\bbW$ that restricts to $(\calT,\Phi)$.
Now, any $\calG$-bundle on $\on{Spec}(\bbW(R^+))$ is trivial. 
Indeed, it is easy to show that $\on{Spec}(\bbW(R^+))$ splits every \'etale cover (see \cite[Proposition 3.2.2]{PR22}).
The choice of a trivialization $\tau:\calT_\bbW\simeq \calG$ specifies a section $(R^\sharp,M,\lambda)\in \Auxsht\Hub{R}$ where $M\in \calG(\bbW(R^+)[\frac{1}{\xi_{R^\sharp}}]$ is the unique element making the following diagram of maps of $\calG$-torsors over $\Spec \bbW(R^+)[\frac{1}{\xi_{R^\sharp}}]$ commute
\begin{center}
\begin{tikzcd}
	\varphi^*\calT_\bbW \arrow{r}{\varphi^*\tau} \arrow{d}{\Phi_\bbW}  & \calG \arrow{d}{M} \\
	\calT_\bbW \arrow{r}{\tau} & \calG. 
\end{tikzcd}
\end{center}
	After chasing definitions one can see that the natural action of $\mathbb{W}^+\calG$ on the set of trivializations acts compatibly with the action specified above.

	Let us prove that $\Auxsht$ is formalizing (see \Cref{defi:formalizing}). 
	Once we prove this, it follows immediately from surjectivity of the map $\Auxsht \to \Shtdos$ that $\Shtdos$ is v-formalizing. 
	Let $\Spa\Hub{T}\in \Perff$, and $\varpi_T\in T^+$ a pseudo-uniformizer. 
	Let $(T^\sharp,M,\lambda)\in \Auxsht\Hub{T}$, we construct a natural transformation $\Spd{(T^+,T^+)}\to\Auxsht$ (see \Cref{definitionspdA}). 
	Let $\Spa\Hub{L}\in \Perff$, a map $f:\Spa\Hub{L}\to \Spd{(T^+,T^+)}$ induces $f:\bbW(T^+)[\frac{1}{\xi_{T^\sharp}}]\to \bbW(L^+)[\frac{1}{\xi_{L^\sharp}}]$, then we let $M_L=f(M)$. 
	Fix a pseudo-uniformizer $\varpi_L\in L^+$. 
	Note that for all $r\in (0,\infty)$ there is a large enough $r'\in (0,\infty)$ for which the following diagram is commutative 
\begin{center}
\begin{tikzcd}
 \Yint{L^+}{[r',\infty]}\arrow{r} \arrow{d} & \Yint{T^+}{[r,\infty]}\arrow{d} \\
	 \Spa{\bbW(L^+)} \arrow{r} & \Spa{\bbW(T^+)}.
\end{tikzcd}
\end{center}
This map allows us to pullback the isogeny $\lambda$ to $\Spa\Hub{L}$. 
The isogeny constructed in this way does not depend of the choices of $\varpi_T$, $\varpi_L$, $r$ or $r'$. 
\end{proof}

Note that $\Sht_\calG(b)$ satisfies the valuative criterion for partial properness over $\Spd\breve{\bbZ}_p=\Spd \bbZ_p\times \Spd \bar{\bbF}_p$ \cite[Definition 18.4]{Sch17}. 
Indeed, all of the data used to define $\Sht_\calG(b)$ (\Cref{defi:moduliofshtukas}) takes place in the exact category of vector bundles over $\calY_{[0,\infty)}^{R^+}$ which is equivalent (by an exact equivalence) to the category of vector bundles over $\Yinf{R^\circ}$ (see \Cref{geometryo-f-Y} for a related discussion).
\begin{lemma}
	\label{lem:closedembeddingauxsht}
	Let $\calG_1\to \calG_2$ be a closed embedding of parahoric group schemes over $\bbZ_p$. Let $b_1\in \calG_1(\breve{\bbQ}_p)$ and let $b_2$ be the image of $b_1$ in $\calG_2(\breve{\bbQ}_p)$.
	Let $V_{b_1}$ be the isocrystal with $\calG_1$ structure associated to $b_1$ and let $V_{b_2}=V_{b_1}\overset{\calG_1}{\times}\calG_2$ be the isocrystal with $\calG_2$ structure associated to $b_2$. 
	The induced map $\Auxshten{1}\to \Auxshten{2}$ is a closed immersion.
\end{lemma}
\begin{proof}
	Let $\Spa\Hub{T} \in \Perff$ be totally disconnected, and let $(M,\lambda)\in \Auxshten{2}\Hub{T}$ (we suppress the untilt from the notation). 
	By \cite[Definition 10.7]{Sch17}, it suffices to prove that the base change along $\Spa\Hub{T}$ is a closed immersion.  
	Abusing notation, we let $(r,\lambda)$ represent the isogeny. 
	After unraveling the definitions, we can think of $M$ and $\lambda$ as ring maps $\calO_{\calG_2}\to \bbW(T^+)[\frac{1}{\xi_{T^\sharp}}]$ and $\calO_{\calG_2}\to B_{[r,\infty]}^{T^+}$.
	Moreover, $\calO_{\calG_1}=\calO_{\calG_2}/I$ where $I$ is the ideal cutting $\calG_1$ inside of $\calG_2$. 
	The base change $\Spa\Hub{T}\times_\Auxshten{2}\Auxshten{1}$ is the subfunctor of $\Spa \Hub{T}$ of those maps
	$\Spa\Hub{R}\to\Spa\Hub{T}$ for which the induced ring morphisms 
	$M:\calO_{\calG_2} \to \bbW(R^+)[\frac{1}{\xi_{R^\sharp}}]$ and
	$\lambda:		\calO_{\calG_2} \to  B^{R^+}_{[r,\infty]} $
	map elements in $I$ to $0$.
	We can fix $\{i_1,\dots, i_n\}\subseteq I$ a set of generators, and let $m_j\in \bbW(T^+)[\frac{1}{\xi_{R^\sharp}}]$ (respectively $t_j\in B_{[r,\infty]}^{T^+}$) denote the image of $i_j$ under $M$ (respectively $\lambda$). 
	The subfunctor in question corresponds to the loci where all the $m_j$ and $t_j$ are $0$.
	It suffices to show the more general statement that for any $m\in \bbW(T^+)[\frac{1}{\xi_{T^\sharp}}]$ (or $t \in B_{[r,\infty]}^{T^+}$) the subfunctor of points in $\Spa\Hub{T}$ for which the element $m_R\in \bbW(R^+)[\frac{1}{\xi_{R^\sharp}}]$ is $0$ (respectively $t_R\in B_{[r,\infty]}^{R^+}$ is $0$) is a closed subfunctor (i.e. the inclusion map is a closed immersion).

	Fix $m\in \bbW(T^+)[\frac{1}{\xi_{T^\sharp}}]$, replacing $m$ by $(\xi_{T^\sharp}^n\cdot m)$ we may assume $m\in \bbW(T^+)$. 
	Using the Teichm\"uller expansion, we can think of $m$ as an element in $(T^+)^{\bbN}$ and $m$ restricts to $0$ if and only if each entry restricts to $0$. 
	This defines a Zariski closed subset of $\Spa\Hub{T}$ \cite[Definition 5.7]{Sch17}. 

	Now fix $t\in B_{[r,\infty]}^{T^+}\subseteq B_{[r,\infty)}^{T^+}$ and let $Z\subseteq |\Yint{R^+}{[r,\infty)}|$ be the set of valuations with $|t|_z=0$. 
	Let us clarify what we mean by this. 
	If $z\in |\Yint{R^+}{[r,\infty)}|$ and $\Spa(A,A^+)\subseteq \Yint{R^+}{[r,\infty)}$ is an open affinoid subset containing $z$ we can think of $z$ as given by an equivalence class of continuous valuations 
	\[|\cdot|_z:A\to \Gamma_z\cup \{0\}.\]
	If $t_A$ denotes the restriction of $t$ to $\Spa(A,A^+)$, then $z\in Z$ if and only if $|t_A|_z=0$.  
	The structure map $\pi:(\Yint{T^+}{[r,\infty)})^\lozenge\to  \Spd\Hub{T}$ is surjective and $\ell$-cohomologically smooth \cite[Proposition 24.5]{Sch17}, therefore it is universally open \cite[Proposition 23.11]{Sch17}. 
	Recall that $|\Spa\Hub{T}|=|\Spd\Hub{T}|$ and similarly $|\Yint{R^+}{[r,\infty)}|=|(\Yint{R^+}{[r,\infty)})^\lozenge|$ \cite[Lemma 15.6]{Sch17}.
The subfunctor of points we consider consists of those maps to $\Spd\Hub{T}$ that factor through 
\[Z'=|\Spa\Hub{T}|\setminus [\pi(|\Yint{R^+}{[r,\infty)}|\setminus Z)]\] which is a closed subset. 
Moreover, we claim that this set is both closed and generalizing, so it defines a closed immersion into $\Spa\Hub{T}$ (\cite[Lemma 7.6]{Sch17}). 
Indeed, if $f:\Spa\Hub{C}\to \Spa\Hub{T}$ is a geometric point that factors through $Z'$ it is because $f(t)\in B^{C^+}_{[r,\infty]}$ is identically $0$ in this ring. 
Nevertheless, if $C^+\subseteq  C'^+\subseteq O_C$ is another open and bounded valuation ring, the map $B^{T^+}_{[r,\infty]}\to B^{C'^+}_{[r,\infty]}$ factors through $B^{C^+}_{[r,\infty]}$.
This shows that the map $\Spa(C,C'^+)\to \Spa\Hub{T}$ also factors through $Z'$.
\end{proof}

\begin{lemma}
	\label{lemma-random-qc}
	The map $\bbW\Sht_\calG(b)\to \Shtdos$ is quasicompact.
\end{lemma}
\begin{proof}
	From \Cref{pro:shtukavformalizing} we know that $\bbW\Sht_\calG(b)\to \Shtdos$ is a $\bbW^+\calG$-torsor. 
	By \cite[Proposition 10.11.(o)]{Sch17}, to deduce that the map is quasicompact it suffices to show that $\bbW^+\calG\to \Spd \bar{\bbF}_p$ is quasicompact. 
	Now, $\calG$ is a finitely presented affine scheme over $\bbZ_p$. 
	Suppose that $\calG=\Spec \calO_\calG$ with $\calO_\calG=\bbZ_p[x_1,\dots,x_n]/(f_1,\dots,f_n)$.
	Then the set of data $t\in \bbW^+\calG(R,R^+)$ corresponds to the choice of $n$-elements $t_1,\dots,t_n\in \bbW(R^+)=(R^+)^\bbN$ subject to the condition that $f_i(t_1,\dots,t_n)=0$. 
	As in the proof of \Cref{lem:closedembeddingauxsht} we see that $\bbW^+\calG$ is a closed subfunctor of the functor 
	\[\Hub{R}\mapsto ((R^+)^\bbN)^n.\]
	This is an infinite dimensional closed unit ball, which is qcqs over $\Spd \bar{\bbF}_p$.
	Indeed, the base change \[\Spa\Hub{R}\times_{\Spd{\bar{\bbF}_p}} \bbB^{\bbN}\] is representable by $\Spa(R\langle T^{\frac{1}{p^\infty}}_i\rangle_{n\in \bbN}, R^+\langle T^{\frac{1}{p^\infty}}_i\rangle_{n\in \bbN})$ which is an affinoid perfectoid and in particular qcqs over $\Spa\Hub{R}$.
\end{proof}

\begin{proposition}
	\label{pro:shtukaclosedimmersion}
	With notation as in \Cref{lem:closedembeddingauxsht} the map $\Sht_{\calG_1}(b_1)\to \Sht_{\calG_2}(b_2)$ is a closed immersion. 
	Moreover, $\Shtdos\to \Spd{\breve{\bbZ}_p}$ is separated.
\end{proposition}
\begin{proof}
	The second claim follows easily from the first one by letting $\calG_2=\calG_1\times_{\bbZ_p} \calG_1$ and letting the map $\calG_1\to \calG_2$ be the diagonal embedding.
	Let us prove the first claim.
	To do this we use that a closed immersion is the same as an injective and proper map of v-sheaves \cite[Lemma 2.1]{AGLR22}.
	For injectivity, let \[t_i=(R^\sharp,\calT_i,\Phi_i,\lambda_i)\in\Sht_{\calG_1}(b_1)\Hub{R} \text{ with } i\in\{1,2\}\] for some $\Spa\Hub{R}\in \Perff$. Assume that $t_1\overset{\calG_1}{\times}\calG_2$ and $t_2\overset{\calG_1}{\times}\calG_2$ are isomorphic where 
	\[t_i\overset{\calG_1}{\times}\calG_2:=(R^\sharp,\calT_i\overset{\calG_1}{\times} \calG_2, \Phi_i, \lambda_i) \text{ with } i\in \{1,2\}.\] 
	We can assume that $\Spa\Hub{R}$ is a product of points. 
	In this case the $t_i$'s lift to  $\Auxshten{1}$, say given by $T_i\in \Auxshten{1}\Hub{R}$ with $T_i:=(R^\sharp,M_i,\lambda_i)$. 
	Since $t_1\overset{\calG_1}{\times}\calG_2\simeq t_2\overset{\calG_1}{\times}\calG_2$, then $T_1\overset{\calG_1}{\times}\calG_2$ and $T_2\overset{\calG_1}{\times}\calG_2$ are in the same $\calG_2(\bbW(R^+))$-orbit. 
	Now, $\lambda_i\in \calG_1(B_{[r,\infty)}^{R^+})$ so $\lambda_1\circ \lambda_2^{-1}\in \calG_1(B_{[r,\infty)}^{R^+})\cap \calG_2(\bbW(R^+))$, this intersection is $\calG_1(\bbW(R^+))$ since $\bbW(R^+)\subseteq B^{R^+}_{[r,\infty)}$. 
	This and \Cref{lem:closedembeddingauxsht} proves that $T_1$ and $T_2$ are in the same $\mathbb{W}^+\calG_1$-orbit, which proves $t_1=t_2$.
	Let us prove $\Sht_{\calG}(b_1)\to \Sht_{\calG}(b_2)$ is proper \cite[Definition 18.1, Proposition 18.3]{Sch17}. 
	Since the map is injective then it is also separated.
	Since both $\Sht_{\calG}(b_1)$ and $\Sht_{\calG}(b_2)$ satisfy the valuative criterion of partial properness \cite[Definition 18.4]{Sch17} over $\Spd \breve{\bbZ}_p$, the map $\Sht_{\calG}(b_1)\to \Sht_{\calG}(b_2)$ is also partially proper. 
	The only thing left to prove is quasi-compactness. 
	Now, by \Cref{lemma-random-qc} and \Cref{lem:closedembeddingauxsht} the composition $\Auxshten{1}\to \Auxshten{2} \to \Sht_{\calG_2}(b_2)$ is a quasi-compact map.
	Since $\Auxshten{1}\to \Sht_{\calG_1}(b_1)$ is surjective, it follows that $\Sht_{\calG}(b_1)\to\Sht_{\calG_2}(b_2)$ is also quasicompact.
\end{proof}

As with $p$-adic Beilinson--Drinfeld Grassmannians, integral moduli spaces of $p$-adic shtukas admit bounded versions. 
Fix a conjugacy class of geometric cocharacters $\mu\in\{\bbG_m\to G_{\bar{\bbQ}_p}\}/\sim$ with field of definition $E$.
In what follows we let $\breve{E}$ denote the compositum of $E$ and $\breve{\bbQ}_p$ in $\bbC_p$.
Fix a geometric point 
\[x:\Spa\Hub{C}\to \Shtdos\times_{\Spd \breve{\bbZ}_p} \Spd O_{\breve{E}} \text{ given by } x=(C^\sharp,f,\calT,\Phi,\lambda)\]
with $(C^\sharp,f)$ an untilt of $C$ over $O_E$.
Observe that by \Cref{thm:anschutzgoodloci} the $\calG$-torsor $\calT$ is trivial. 
Let us fix a trivialization of $\tau:\calT\to \calG$. 
The morphism $\tau\circ \Phi:\varphi^{*}\calT\to \calG$ defines a map 
\[y_{\tau,x}:\Spa\Hub{C}\to \Gr_\calG\times_{\Spd \bbZ_p}\Spd O_{{E}}.\]
 We say that $x$ has relative position bounded by $\mu$ if $y_{\tau,x}$ factors through $\Grm{O_{E}}$ as in \Cref{thelocalmodeldefi}. 
 By \cite[Proposition 4.13]{AGLR22}, this condition doesn't depend on the choice of $\tau$ 

\begin{definition}
	\label{defi:moduliofshtukasbound}
	We let $\Shtm{O_{\breve{E}}}\subseteq \Shtdos\times_{\Spd \breve{{\bbZ}}_p} \Spd O_{\breve{E}}$ denote the subfunctor of tuples 
	\[((R^\sharp,\iota),f,\calT,\Phi,\lambda)\] 
	for which the shtuka $(\calT,\Phi)$ is point-wise bounded by $\mu$.
\end{definition}

\begin{remark}
	\label{remark-boundednesscond}
One can give a more conceptual reformulation of \Cref{defi:moduliofshtukasbound} by considering Hecke stacks as follows. 
Let $L^+\calG_{\Spd \bbZ_p}$ denote the positive loop group considered in \cite[$\mathsection$ 19.1]{SW20} and \cite[$\mathsection$ 4.2]{AGLR22}. Let \[\on{Hk}_\calG:=[L^+\calG_{\Spd \bbZ_p}\backslash\Gr_\calG]\] 
and \[\on{Hk}_{\calG,O_{{E}}}^{\leq \mu}:= [L^+\calG_{\Spd O_E}\backslash\Grm{O_{E}}]\subseteq \on{Hk}_{\calG}\times_{\Spd \bbZ_p}\Spd O_E.\] 
Then we have a Cartesian diagrams
\begin{center}
\begin{tikzcd}
	\Grm{O_{E}} \arrow{r} \arrow{d}  & \Gr_\calG \times_{\Spd {\bbZ}_p} \Spd O_{{E}} \arrow{d} &
	\Shtm{O_{\breve{E}}}\arrow{r} \arrow{d}  & \Shtdos\times_{\Spd \breve{{\bbZ}}_p} \Spd O_{\breve{E}} \arrow{d} \\
	\on{Hk}_{\calG,O_{{{E}}}}^{\leq \mu} \arrow{r} & \on{Hk}_{\calG}\times_{\Spd \bbZ_p}\Spd O_{{E}} &
 \on{Hk}_{\calG,O_{{\breve{E}}}}^{\leq \mu} \arrow{r} & \on{Hk}_{\calG}\times_{\Spd \bbZ_p}\Spd O_{\breve{E}}. \\
\end{tikzcd}
\end{center}
\end{remark}

\begin{remark}
	Whenever $\calG$ is reductive over $\bbZ_p$, the field of definition $\mu$ is always an unramified extension of $\bbQ_p$ since the group itself $G=\calG_{\bbQ_p}$ splits over an unramified extension.
	In this case, $O_{\breve{E}}={\breve{\bbZ}_p}$. 
\end{remark}

Let $\bbW\Sht_{\calG,O_{\breve{E}}}(b)={{\bbW\Sht}}_\calG(b)\times_{\Spd \breve{\bbZ}_p}\Spd O_{\breve{E}}$, let $\on{Sht}_{\calG,O_{\breve{E}}}(b)=\Sht_\calG(b) \times_{\Spd \breve{\bbZ}_p}\Spd O_{\breve{E}}$ and let $\Auxshtm$ denote the base change of 
\[{{\bbW\Sht}}_{\calG,O_{\breve{E}}}(b) \to \on{Sht}_{\calG,O_{\breve{E}}}(b)\]
along $\Shtm{O_{\breve{E}}}$. 
We can consider the map  
\[\Auxshtm\to \Gr_{\calG,O_{\breve{E}}}\]
constructed as follows. Fix 
	\[(R^\sharp,f,M,\lambda)\in \Auxshtm\Hub{R},\] 
	and recall that $M\in \calG(\bbW(R^+)[\frac{1}{\xi_{R^\sharp}}])$ which we may think of as an automorphism of the trivial $\calG$-torsor defined over $\Spec \bbW(R^+)[\frac{1}{\xi_{R^\sharp}}]$.
	We can consider the tuple 
	\[(R^\sharp,f,\calG_M,\psi_M)\in \Gr_{\calG,O_{\breve{E}}}\Hub{R},\] 
	where $\calG_M$ is the trivial $\calG$-torsor over $\calY^{R^+}_{[0,\infty)}$ and $\psi_M:\calG_M\to \calG$ is the lattice with $\calG$-structure \Cref{defi:latticesGstructure} obtained from restricting to $\calY^{R^+}_{[0,\infty)}\setminus V(\xi_{R^\sharp})$ the map $M:\calG\to \calG$ which is initially defined over $\Spec \bbW(R^+)[\frac{1}{\xi_{R^\sharp}}]$.
\begin{proposition}
	\label{pro:shtukaboundedclosed}
	Let $\mu\in \{\bbG_m\to G_{\bar{\bbQ}_p}\}/\sim$ with field of definition $E$. Then 
	\[\Shtm{O_{\breve{E}}}\to \Shtdos\times_{\Spd \breve{\bbZ}_p} \Spd O_{\breve{E}}\] 
	is a closed immersion. Moreover, $\Shtm{O_{\breve{E}}}$ is v-formalizing as in \Cref{defi:formalizing}.
\end{proposition}
\begin{proof}
	We have a pair of Cartesian diagrams:
		\begin{center}
	\begin{tikzcd}
		\Auxshtm \ar{r}\ar{d} & {{\bbW\Sht}}_{\calG,O_{\breve{E}}}(b) \ar{d} &
		\Auxshtm \ar{r}\ar{d} & {{\bbW\Sht}}_{\calG,O_{{\breve{E}}}}(b) 	\ar{d} \\
			\Grm{O_{\breve{E}}}\ar{r}& \Gr_{\calG,O_{\breve{E}}} &
			\Shtm{O_{\breve{E}}} \ar{r} & {\on{Sht}}_{\calG, O_{\breve{E}}}(b) 
	\end{tikzcd}
		\end{center}
		Since being a closed immersion can be checked v-locally on the target (\cite[Proposition 10.11]{Sch17}), and since ${{\bbW\Sht}}_{\calG,O_{{\breve{E}}}}(b) \to {{\on{Sht}}_{\calG, O_{\breve{E}}}}(b)$ is surjective $\Shtm{O_{\breve{E}}} \to {{\on{Sht}_{\calG,O_{\breve{E}}}}} $ is a closed immersion. 
		Moreover, by \Cref{pro:grassmanianreduced} and \Cref{thm:AGLR} the map $\Grm{O_{\breve{E}}}\to \Gr_{\calG,{O_{\breve{E}}}}$ is formally adic since they are both formally adic over $\Spd O_{\breve{E}}$ (see discussion below \Cref{formallyadicdefi}). 
		We claim that $\Auxshtm$ is formalizing and consequently that $\Shtm{O_{{\breve{E}}}}$ is v-formalizing.
		Indeed, given a map $\Spa\Hub{R}\to \Auxshtm$ we get a map $\Spd (R^+,R^+)\to {{\bbW\Sht}}_{\calG,O_{\breve{E}}}(b)$ by \Cref{pro:shtukavformalizing}.   
		Moreover, the induced map $\Spd (R^+,R^+)\to \Gr_{\calG,{O_{\breve{E}}}}$ factors through $\Grm{O_{\breve{E}}}$ since $\Spa(R,R^+)\to \Gr_{\calG,{O_{\breve{E}}}}$ does.
		This follows from the fact that $\Spd R^+\times_{\Gr_{\calG,{O_{\breve{E}}}}} \Grm{O_{\breve{E}}}$ is a formally adic closed subsheaf of $\Spd (R^+,R^+)$ which by \cite[Lemma 3.31]{Gle22} has to agree with $\Spd (R^+,R^+)$.  
\end{proof}
At the moment we have only proven that $\Shtm{O_{\breve{E}}}$ is a separated (\Cref{pro:shtukaclosedimmersion}) and v-formalizing v-sheaf (\Cref{pro:shtukaboundedclosed}). 
This is not enough to construct a specialization map. 
Indeed, we still have to show that it is formally separated \Cref{formallyseparateddefi} (see \Cref{defi:specializingvsheaf}). 
If $\Shtm{O_{\breve{E}}}$ was formally adic over $\Spd O_{\breve{E}}$, then formal separatedness would follow easily, but this is not the case.
Instead we rely on the following lemma, see \cite[Lemma 3.30]{Gle22} for a proof.
\begin{lemma}
	\label{rem:diamondofreduced}
Let $\calF$ be a small v-sheaf.	
The diagonal $\calF\to \calF\times \calF$ is formally adic if and only if the adjunction morphism $(\calF^\red)^\diamond\to \calF$ is injective. 
\end{lemma}

To show that \Cref{rem:diamondofreduced} holds for $\Shtm{O_{\breve{E}}}$ we identify 
\[(\Shtm{O_{\breve{E}}})^\red:\PSch^\op_{\bar{\bbF}_p}\to \Sets\]
with an affine Deligne--Lusztig variety, and we identify the adjunction map
\[(\Shtm{O_{\breve{E}}})^\Red\to \Shtm{O_{\breve{E}}}\]
as a closed subsheaf. 

\subsubsection{Affine Deligne--Lusztig varieties.}
Recall the pioneering works of Rapoport and Kottwitz in which they initiate the study of affine Deligne--Lusztig varieties \cite{MR1742251}, \cite{MR1966756}, \cite{Rap05}. 
Although at the time, the general definition of an affine Deligne--Lusztig varieties only made sense as a set, one could use the theory of $p$-divisible groups and Dieudonn\'e theory to prove that in some cases these sets could be realized as the $\bar{\bbF}_p$-points of a Rapoport--Zink space \cite{RZ96}, which are formal schemes whose reduced special fiber is a scheme locally of finite presentation over $\bar{\bbF}_p$.
In \cite{CKV15} Chen, Kisin and Viehmann found a way in which one could endow general affine Deligne--Lusztig with geometric structure, and in particular they were able to meaningfully talk about their connected components. 
With the introduction of the Witt vector affine flag variety $\calF\ell_{\calG,\bbW}$ \cite{Zhu17} and the proof of its representability \cite{BS17}, we now know that affine Deligne--Lusztig varieties always arise as the $\bar{\bbF}_p$-points of a perfect scheme perfectly of finite presentation \cite[$\mathsection$ A.2]{Zhu17} over $\bar{\bbF}_p$ \cite[Theorem 1.2.(1)]{HV20}. 
We denote this scheme by $X^{\leq \mu}_\calG(b)$ and by abuse of notation we still refer to it as the affine Deligne--Lusztig variety.
\begin{definition}
	\label{definition-aDLV}
	Let $\calG$ be a parahoric group scheme over $\bbZ_p$.
	Let $\mu\in \{\bbG_m\to G_{\bar{\bbQ}_p}\}/\sim$ be a conjugacy class of geometric cocharacters, and let $b\in G(\breve{\bbQ}_p)$. 
	The affine Deligne--Lusztig variety $X^{\leq \mu}_\calG(b)$ is the unique closed subfunctor of $\calF\ell_{\calG,\bbW}$ such that
	\[X^{\leq \mu}_\calG(b)(\bar{\bbF}_p)=\{g\cdot \calG(\breve{\bbZ}_p)\in \calF\ell_{\calG,\bbW}(\bar{\bbF}_p)\mid \calG(\breve{\bbZ}_p)\cdot g^{-1}b\phi(g) \cdot  \calG(\breve{\bbZ}_p)\in \on{Adm}(\mu)\}\]
	Here $\on{Adm}(\mu)$ denotes the $\mu$-admissible set of Kottwitz--Rapoport \cite{KR00}, \cite[Definition 3.11]{AGLR22}.
\end{definition}
By \cite[$\mathsection$ 3.1.1, Lemma 1.22]{Zhu17} \cite[Theorem 1.2.(1)]{HV20}, $X^{\leq \mu}_\calG(b)$ is representable by a perfect scheme that is locally perfectly finitely presented.
In what follows, we recall an alternative description of affine Deligne--Lusztig varieties that is more convenient for our purposes.

	Let $\calE_b$ denotes the $\varphi$-module with $\calG$-structure over $\Spec \bbW(R)[\frac{1}{p}]$ defined by the isocrystal with $\calG$-structure that $b\in G(\breve{\bbQ}_p)$ induces. 
	More precisely, $\calE_b=(\calG,\Phi_b)$ where $\calG$ is simply the trivial $\calG$-torsor and 
	\[\Phi_b:\varphi^*\calG\to \calG\]
	is the only isomorphism of $\Spec \breve{\bbQ}_p$ conjugate to $b\in \calG(\breve{\bbQ}_p)$ under the canonical identification $\varphi^*\calG\simeq \calG$.

\begin{definition}
	\label{defi:ADLVunb}
We consider a functor
\[\calS_\calG(b):\PSch^\op\to \Sets\]
with formula
\[\calS_\calG(b): \Spec R\mapsto   \{(\calT,\Phi,\lambda)\}/\sim\]
	where $\calT$ is a $\calG$-torsor over $\on{Spec}(\bbW(R))$, $\Phi:\varphi^{*}\calT\to \calT$ is an isomorphism over $\on{Spec}(\bbW(R)[\frac{1}{p}])$ and $\lambda:\calT\to \calE_b$ is a $\varphi$-equivariant isomorphism over $\on{Spec}(\bbW(R)[\frac{1}{p}])$.
\end{definition}

	\begin{proposition}
		\label{flagarietyandshtukas}
		The formula 
		\[(\calT,\Phi,\lambda)\mapsto (\calT, \lambda)\]
		defines an isomorphism $f:\calS_\calG(b)\to \calF\ell_{\calG,\bbW}$.
	\end{proposition}
	\begin{proof}
		It suffices to construct an inverse functor. 
		To do this it suffices to observe that $\Phi$ is completely determined by $\lambda$, we spell this out as follows.
		Given $(\calT,\psi)\in \calF\ell_{\calG,\bbW}(\Spec R)$
		we construct an element $(\calT',\Phi,\lambda)\in \calS_\calG(b)(\Spec R)$, to do this
		we consider the following diagram defined over $\Spec \bbW(R)[\frac{1}{p}]$
		\begin{center}
		\begin{tikzcd}
			\varphi^*\calT \arrow{r}{\varphi^*\psi} \arrow{d}  & \varphi^*\calG \ar{d}{\Phi_b} \ar{r}{\on{can}}  & \calG \ar{d}{b} \\
			\calT \arrow{r}{\psi} & \calG  \ar{r}{\on{id}} & \calG.
		\end{tikzcd}
		\end{center}
		Then $f^{-1}(\calT,\psi)=(\calT',\Phi,\lambda)$ is obtained by letting $\calT':=\calT$, letting $\Phi:=\psi^{-1}\circ \Phi_b\circ \circ (\varphi^*\psi)$, and by letting $\lambda:= \psi$.  
	\end{proof}

	Let $C$ be an algebraically closed characteristic $p$ field, suppose we are given $x=(\calT,\Phi,\lambda)\in \calS_\calG(b)(\Spec C)$.
	Note that $\calT$ is a trivial $\calG$-torsor, and after fixing a trivialization $\tau:\calT\to \calG$  
	over $\Spec \bbW(C)$
	we obtain a map \[x_\tau:\Spec C\to \calF\ell_{\calG,\bbW}.\]
	by considering 
	\[(\varphi^*\calT,\tau\circ \Phi)\in \calF\ell_{\calG,\bbW}(\Spec C).\]
	We say that $x$ is bounded by $\mu$ if $x_\tau$ factors through $\GrWmfe\subseteq \calF\ell_{\calG,\bbW}$.
	Since the admissible locus is a union of Schubert varieties it is stable under the action of $\calG(\bbW(C))$, so this condition does not depend on the trivialization chosen.

	\begin{definition}
	\label{defi:ADLV}
	We let $\calS^{\leq \mu}_\calG(b)\subseteq \calS_\calG(b)$ denote the closed subfunctor of tuples $(\calT,\Phi,\lambda)$ such that on geometric points of $\Spec R$ the induced pair $(\calT,\Phi)$ is bounded by $\mu$.
	\end{definition}

\begin{proposition}
	The isomorphism $f:\calS_\calG(b)\to \calF\ell_{\calG,\bbW}$ of \Cref{flagarietyandshtukas} restricts to an isomorphism
	\[\calS^{\leq \mu}_\calG(b)\to X^{\leq \mu}_\calG(b).\]
\end{proposition}
\begin{proof}
	Since $\calF\ell_{\calG,\bbW}$ has an ind-presentation by proper perfectly of finite presentation closed subschemes, any closed subfunctor is determined by its $\bar{\bbF}_p$-points. 	
	So it suffices to show that $f(\calS^{\leq \mu}_\calG(b)(\bar{\bbF}_p))=X^{\leq \mu}_\calG(b)(\bar{\bbF}_p)$.
	Let $(\calG,\psi_g)\in \calF\ell_{\calG,\bbW}(\bar{\bbF}_p)$ where $\psi_g:\calG\to \calG$ is the isomorphism defined over $\Spec \breve{\bbQ}_p$ associated to $g\in \calG(\breve{\bbQ}_p)$.
	Then $f^{-1}(\calG,\psi_g)=(\calG,\psi_{g^{-1}}\circ \Phi_b\circ \varphi^*\psi_{g},\psi_g)$.
	Moreover, $f^{-1}(\calG,\psi_g)\in \calS^{\leq \mu}_\calG(b)(\bar{\bbF}_p)$ if and only if $(\varphi^*\calG,\psi_{g^{-1}}\circ \Phi_b\circ \varphi^*\psi_{g})\in \GrWmfe(\bar{\bbF}_p)$.
	But the pair $(\varphi^*\calG,\psi_{g^{-1}}\circ \Phi_b\circ \varphi^*\psi_{g})$ is isomorphic to the pair $(\calG,\psi_{g^{-1}b\phi(g)})$ so they define the same point in $\calF\ell_{\calG,\bbW}$. 
	Now, for an element $a\in \calG(\breve{\bbQ}_p)$ the pair $(\calG,\psi_a)\in \GrWmfe(\bar{\bbF}_p)$ if and only if 
	\[a\in \calG(\breve{\bbZ}_p)\cdot \on{Adm}(\mu)\cdot \calG(\breve{\bbZ}_p)\]
	since by definition $\GrWmfe$ is the only closed subsheaf of $\calF\ell_{\calG,\bbW}$ satisfying that 
	\[\GrWmfe(\bar{\bbF}_p)=\calG(\breve{\bbZ}_p)\cdot \on{Adm}(\mu)\cdot \calG(\breve{\bbZ}_p)/\calG(\breve{\bbZ}_p)\subseteq \calG(\breve{\bbQ}_p)/\calG(\breve{\bbZ}_p)=\calF\ell_{\calG,\bbW}(\bar{\bbF}_p)\]
	From here it is clear that $(\calG,\psi_g)\in X^{\leq \mu}_\calG(\bar{\bbF}_p)$ if and only if $(\calG,\psi_{g^{-1}b\phi(g)})\in\GrWmfe(\bar{\bbF}_p)$. 
	The reasoning above shows that this is also equivalent to $f^{-1}(\calG,\psi_g)\in \calS^{\leq \mu}_\calG(b)(\bar{\bbF}_p)$.
\end{proof}

From now on we use $X^{\leq \mu}_\calG(b)$ to denote the functor $\calS^{\leq \mu}_\calG(b)$.
\begin{proposition}
	\label{pro:adlvasreducedfunctor}
	The following statements hold.
	\begin{enumerate}
		\item There is a natural identification $\adlvdos\xrightarrow{\sim}(\Shtm{O_{\breve{E}}})^\red$. 
		\item The adjunction map $(\adlvdos)^\dia\to \Shtm{O_{\breve{E}}}$ is injective. 
		\item $\Shtm{O_{\breve{E}}}$ is a specializing v-sheaf.
	\end{enumerate}
\end{proposition}
\begin{proof}
	We first construct a map $j:\adlvdos \to (\Shtm{O_{\breve{E}}})^\red$. 
	By adjunction, we may construct $h:(\adlvdos)^\dia \to \Shtm{O_{{\breve{E}}}}$ instead. 
	Before sheafification, a map $\Spa\Hub{T}\to (\adlvdos)^\dia$ is given by data $(\calT,\Phi,\lambda)$ over $\Spec T^+$ as in \Cref{defi:ADLVunb} and \Cref{defi:ADLV}.
	Restricting to the appropriate loci defines a map $\Spa\Hub{T}\to \Shtm{O_{{\breve{E}}}}$.  
	Indeed, $\calT$ is a $\calG$-torsor over $\Spec \bbW(T^+)$ which can be restricted to $\calY^{T^+}_{[0,\infty)}$. 
	Then $\Phi$ is a morphism in the category of $\calG$-torsors over $\Spec \bbW(T^+)[\frac{1}{p}]$ which can be restricted to a morphism of $\calG$-torsors over $\calY^{T^+}_{[0,\infty)}\setminus V(p)$ that is meromorphic along $V(p)$.
	Finally, $\lambda$ is defined over $\Spec \bbW(T^+)[\frac{1}{p}]$ which can be restricted to $\calY^{T^+}_{(0,\infty)}$.

	Since $\adlvdos$ is representable, full faithfulness of $\diamond$ on the category of perfect schemes implies that the unit of the adjunction $\adlvdos\to ((\adlvdos)^\dia)^\red$ is an isomorphism \cite[Proposition 3.16]{Gle22} \cite[Proposition 18.3.1]{SW20}. 
	To prove $j$ is injective it suffices to prove that $h$ is. 
	Indeed, $(-)^\red$ is a right adjoint functor so it preserves monomorphisms and $j$ is obtained as the following composition
	\[\adlvdos\xrightarrow{\simeq} ((\adlvdos)^\dia)^\red \xrightarrow{h^\red} (\Shtm{O_{\breve{E}}})^\red.\]
	Take maps $g_1,g_2:\Spa\Hub{R}\to (\adlvdos)^\dia$.
	Injectivity can be proved v-locally. 
	Moreover, since $\adlvdos$ is representable, we can identify the set of maps 
	\[(\adlvdos)^\dia(\Spa\Hub{R})\simeq \adlvdos(\Spec R^+)\]
	 as long as $\Spa\Hub{R}$ splits every open cover. 
	 Totally disconnected spaces form a basis for the v-topology with this property \cite[Definition 7.1, Lemma 7.18]{Sch17}.
	Without loss of generality, the $g_i$ are given by data $(\calT_i,\Phi_i,\lambda_i)$ over $\on{Spec}(\bbW(R^+))$, $\on{Spec}(\bbW(R^+)[\frac{1}{p}])$ and $\on{Spec}(\bbW(R^+)[\frac{1}{p}])$ and the $h\circ g_i$ is given by restriction to $\Yinf{R^+}$, $\Yinf{R^+}\setminus V(p)$ and $\Yint{R^+}{[r,\infty)}$ respectively. 
	It follows from $h\circ g_1=h\circ g_2$ and \Cref{thm:puncturedtoAinfked} that $g_1=g_2$. 
	This finishes showing that $h$ and $j$ are injective.
	Granted that $j$ is surjective, this also finishes the proof of the first statement.
	
	Let us prove surjectivity of $j$. 
	We show that every map $f:\Spec A\to (\Shtm{O_{{\breve{E}}}})^\red$ can be lifted uniquely to a map $\Spec A\to \adlvdos$. 
	Now, any such map factors as the composition 
	\[\Spec A\xrightarrow{\simeq} ((\Spec A)^\diamond)^\red \xrightarrow{g^\red}(\Shtm{O_{{\breve{E}}}})^\red\]
	where $g:(\Spec A)^\diamond\to \Shtm{O_{\breve{E}}}$ is the map that corresponds to $f$ by adjunction.
	It suffices to lift $g$ to a map $e:\Spec A^\diamond \to \adlvdos^\diamond$ with $h\circ e=g$, since in this case we get a diagram 
	\begin{center}
	\begin{tikzcd}
		\Spec A \arrow{r}{\simeq} \arrow{d}  & ((\Spec A)^\diamond)^\red  \arrow{d}{e^\red} \ar{rd}{g^\red} \\
		\adlvdos \arrow{r}{\simeq} & (\adlvdos^\diamond)^\red \ar{r}{h^\red} & \Shtm{O_{{\breve{E}}}})^\red.
	\end{tikzcd}
	\end{center}
	From $g$ we can construct a commutative diagram
\begin{center}
\begin{tikzcd}
	(\Spec A)^\diamond\arrow{r}{g} \arrow{d}{e}  & \Shtm{O_{{\breve{E}}}}  \\
	\adlvdos^\diamond \arrow{ur}{h} & 
\end{tikzcd}
\end{center}
by interpreting all of the objects involved as v-sheaves in the full subcategory of $\Perff$ whose objects are product of points \Cref{defi:prodpoints} and by exhibiting a natural transformation that gives rise to $e$.

Fix $t:\Spa\Hub{R}\to \on{Spec}(A)^\dia$ a map with $\Spa\Hub{R}$ a product of points. 
We construct below a unique $f'$ fitting in the following commutative diagram
			\begin{center}
				\begin{tikzcd}
					\Spa\Hub{R} \ar{r}{t}\ar{d}{f'}& \on{Spec}(A)^\dia \ar{d}{g}\\
					(\adlvdos)^\dia \ar{r}{h} & \Shtm{O_{{\breve{E}}}}.
				\end{tikzcd}
			\end{center}
	
			Fix a pseudo-uniformizer $\varpi\in R^+$, we let $\Spa\Hub{R_\infty}$ be a second product of points defined by $R_\infty^+=\prod_{i=1}^\infty R^+$ with pseudo-uniformizer $\varpi_{R_\infty}=(\varpi^i)_{i=1}^\infty$. 
			This product of points comes with a family of closed embeddings $\iota_i:\Spa\Hub{R}\to \Spa\Hub{R_\infty}$ given in coordinates by the projections onto the $i$th-factor. 
			The diagonal $\Delta_t:A\to \prod_{i=1}^\infty R^+$ induces ${\Delta_t}:\Spa\Hub{R_\infty}\to \on{Spec}(A)^\dia$ with ${\Delta_t}\circ \iota_i=t$ for every $i$. 
			Since $\Spa\Hub{R_\infty}$ is a product of points, by \Cref{pro:prodpointextendg-tors}, the map $g\circ {\Delta_t}$ can be represented by a triple $(\calT_{R_\infty},\Phi_{R_\infty},\lambda_{R_\infty})$ with $\calT_{R_\infty}$ trivial. 
			After choosing a trivialization, $\lambda_{R_\infty}$ is given by a map $\calO_\calG\to B^{{R^+_\infty}}_{[r,\infty]}$. 
			Moreover, since $g\circ {\Delta_t}\circ\iota_i=g\circ {\Delta_t}\circ\iota_j$ for all $i$ and $j$, the various $\lambda^*_i:\calO_\calG\to B^{{R^+_\infty}}_{[r,\infty]}\xrightarrow{\iota_i} B^{R^+}_{[r_i,\infty]}$  
	all lie in the same $\calG(\bbW(R^+))$-orbit. 
	Clearly, $\calG(\bbW({R^+_\infty}))=\prod_{i=1}^\infty \calG(\bbW(R^+))$. 
	By changing the trivialization, we may assume that $r_i=r_j=:r$ and $\lambda^*_i=\lambda^*_j=:\lambda^*_R$ for all $1\leq i,j<\infty$. 
	We claim that $\lambda^*_{R}$ factors through $\bbW(R^+)[\frac{1}{p}]$. 

	Take $t\in \calO_\calG$ and consider $s=\lambda^*_{R_\infty}(t)\in B^{{R^+_\infty}}_{[r,\infty]}$. 
	After enlarging $r$ if necessary, we may assume $r=n\in \bbN$. 
	In particular, $p^k\cdot s$ lies in the $p$-adic completion of $\bbW({R^+_\infty})[\frac{[\varpi_{R_\infty}]}{p^n}]$ for some $k$. 
	Let $L_i\subseteq B^{R^+}_{[\frac{1}{i},\infty]}$ denote the $p$-adic completion of $\bbW(R^+)[\frac{[\varpi^i]}{p}]$, this is a ring of definition.
	The argument above shows that for all $i$ the element $\iota_i(p^k\cdot s)$ lies in $L_i$.
	In other words, $p^k\cdot \lambda^*_R(t)\in \bigcap_{i\in \bbN} L_i$. 
	By \Cref{annoyingcomputation} below, this intersection is $\bbW(R^+)$ proving the claim. 

	Since the triple $(\calT_{R_\infty},\Phi_{R_\infty},\lambda_{R_\infty})$ is defined over $\on{Spec}(\bbW(R^+))$ and $\on{Spec}(\bbW(R^+)[\frac{1}{p}])$, it defines a map $\on{Spec}(R^+)\to \adlvdos$. 
	The composition, $\Spa\Hub{R}\to \on{Spec}(R^+)^\dia \to (\adlvdos)^\dia$, defines $f'$. 
	This finishes the proof of the first statement. 

	To prove $\Shtm{O_{{\breve{E}}}}$ is a specializing v-sheaf \Cref{defi:specializingvsheaf} we need to prove it is formally separated and v-formalizing. 
	That it is formally separated follows from \Cref{rem:diamondofreduced} (\cite[Lemma 3.30]{Gle22}) and \Cref{pro:shtukaclosedimmersion}. 
	That it is v-formalizing is \Cref{pro:shtukaboundedclosed}.
\end{proof}

\begin{lemma}
	\label{annoyingcomputation}
	Let $\Spa(R,R^+)$ be affinoid perfectoid.
	Let $\varpi\in R^+$ be a choice of pseudo-uniformizer. 
	For $n\in \bbN$ let $L_n$ denote the $p$-adic completion of $\bbW(R^+)[\frac{[\varpi^n]}{p}]$. 
	For $n\geq m$ consider the natural inclusion $L_n\subseteq L_m$.
	Then the natural map induces an isomorphism of rings
	\[\bbW(R^+)\xrightarrow{\simeq} \bigcap_{i\in \bbN} L_i.\]
\end{lemma}
\begin{proof}
	Let $A^i_n=(\bbW(R^{+})[\frac{[\varpi^i]}{p}])/(p^n,[\varpi]^n)$.
	Then $L_i=\varprojlim_{n\in \bbN}A^i_n$ and we can interpret the intersection as a limit
	\[\bigcap_{i\in \bbN} L_i=\varprojlim_{i\in\bbN}\varprojlim_{n\in \bbN} A^i_n.\] 
	We can compute $\varprojlim_{n\in \bbN}\varprojlim_{i\in \bbN} A^i_n$ instead. 
	Let $A^\infty_n$ denote the ring $\bbW(R^+)/(p^n,[\varpi]^n)$.
	We claim that the set of maps $\bbW(R^+)\to A^i_n$ induces an identification $A^\infty_n\simeq \varprojlim_{i\in \bbN} A^i_n$.
	Indeed, 
	\[A^i_n= (\bbW(R^+)[T_i]/pT_i-[\varpi^i])\otimes_{\bbW(R^+)} A^\infty_n\] 
	and the transition maps $A^{i+1}_n\to A^i_n$ are the one obtained from the rule $T_{i+1}\mapsto [\varpi]T_i$. 
	When $i\geq n$ this simplifies to $A^i_n=\bbW(R^+)[T_i]/(p^n,[\varpi]^n,pT_i)$ and for such $i$ the map $A^{i+n}_n\to A^i_n$ is the map with $T_{i+n}\mapsto 0$.  
	We see that the only polynomials in $A^i_n$ that have a preimage in $A^{i+n}_n$ are the constant ones.
	This shows $\bbW(R^+)/(p^n,[\varpi]^n)\simeq \varprojlim_{i\in \bbN} A^i_n$.
	Since $\bbW(R^+)$ is $(p,[\varpi])$-adically complete, $\bbW(R^+)\simeq \varprojlim_{n\in \bbN}\varprojlim_{i\in \bbN} A^i_n$.
\end{proof}

\begin{lemma}
	\label{lem:shtukredclosedimme}
	\label{pro:modulishtukformssmeltedkimberl}
	The adjunction map $(\adlvdos)^\dia\to \Shtm{O_{{\breve{E}}}}$ arising from the identification of \Cref{pro:adlvasreducedfunctor} is a closed immersion. 
	In particular, $(\Shtm{O_{{\breve{E}}}},\Shtm{{{\breve{E}}}})$ is a smelted kimberlite as in \Cref{definitionkimberlite}.
\end{lemma}
\begin{proof}
	Recall that $\adlvdos$ is a closed subfunctor of $\calF\ell_{\calG,\bbW}$ \cite[$\mathsection$ 3.1.1, Lemma 1.22]{Zhu17}.
	Moreover, by \cite[Corollary 9.6]{BS17} we may write 
	\[\calF\ell_{\calG,\bbW}=\cup_{n\in \bbN} \calS_n\]  
	as in increasing union of schemes $\calS_n$ each of which is the perfection of a projective scheme over $\bar{\bbF}_p$. 

	We can write \[\adlvdos^\dia=\bigcup_{n\in \bbN} (\adlvdos\cap \calS_n)^\dia.\]
	Each of the terms $(\adlvdos\cap \calS_n)^\dia$ is proper over $\Spd \bar{\bbF}_p$, since they come from a proper perfectly finitely presented scheme over $\bar{\bbF}_p$. 
	Consequently, $(\adlvdos\cap \calS_n)^\dia \to \Shtm{O_{{\breve{E}}}}$ is a closed immersion \cite[Lemma 2.1]{AGLR22}. 

	Now, $\adlvdos$ is locally perfectly of finite type \cite[Theorem 1.2]{HV20}. 
	For all $x\in |\adlvdos|$, we fix $U_x=\Spec A$ with $x\in U_x$, and such that $A$ is the perfection of a finite type algebra over $\bar{\bbF}_p$. 
	Observe that $|U_x|$ is a Noetherian topological space. 
	We claim that there is $n\in \bbN$ for which $U_x=U_x\cap \calS_n$. 
Indeed, if $U_x$ is Noetherian then every Zariski closed subset is an open subset in the constructible topology. 
For a scheme $X$ we consider its constructible topology and denote it by $X^{\on{cons}}$. 
Since $U^{\on{cons}}_x$ is compact and $U_x^{\on{cons}}\subseteq \cup_{n\in \bbN} (U_x\cap \calS_n)^{\on{cons}}$ is a nested open cover, we may find $n$ large enough such that $U_x\subseteq \calS_n$. 
	
	By \Cref{pro:adlvasreducedfunctor} and \cite[Proposition 4.14]{Gle22} we have a specialization map \[\on{sp}:|\Shtm{O_{{\breve{E}}}}|\to |\adlvdos|.\]
	Let $V_x=(\on{sp})^{-1}(U_x)$ for $x\in |\adlvdos|$ and $U_x$ as above, this forms an open cover of $\Shtm{O_{{\breve{E}}}}$. 
	By the last statement of \Cref{lemmaspecializing} above, to show that $\Shtm{O_{{\breve{E}}}}$ is a prekimberlite it suffices to show that each $V_x$ is a prekimberlite.
	This reduces us to showing that $(V_x^\red)^\dia \to V_x$ is a closed immersion. 
	Now, the adjunction map $(U_x)^\dia\to V_x$ fits in the following diagram
	\begin{center}
		\begin{tikzcd}
			U_x^\dia \ar{r}{\on{Id}}\ar{d}&	U_x^\dia \ar{r}\ar{d}& V_x \ar{d}\\
			(\calS_n\cap \adlvdos)^\dia \ar{r} & (\adlvdos)^\dia \ar{r} & \Shtm{O_{{\breve{E}}}}.
		\end{tikzcd}
	\end{center}
	By \Cref{lemmaspecializing}, the map $V_x\to \Shtm{O_{{\breve{E}}}}$ is formally adic, so the above diagram is Cartesian.
	This shows that $U_x^\dia\to V_x$ is a closed immersion. 

	Let us show the final statement. \Cref{pro:adlvasreducedfunctor} proves that $\Shtm{O_{{\breve{E}}}}$ is a specializing v-sheaf and that $(\Shtm{O_{{\breve{E}}}})^\red$ is represented by a scheme. 
	The argument above shows that $\Shtm{O_{{\breve{E}}}}$ is a prekimberlite as in \Cref{thisisaprekimberlitehurray}. 
	Since $\Shtm{O_{{\breve{E}}}}\to \on{Spd}(O_{{\breve{E}}})$ is partially proper, by \cite[Proposition 4.32, Definition 4.30]{Gle22} it is a valuative prekimberlite. 
	Finally, by \cite[Theorem 23.1.4]{SW20} $\Shtm{{{\breve{E}}}}$ is a locally spatial diamond. Consequently, $(\Shtm{O_{{\breve{E}}}},\Shtm{{{\breve{E}}}})$ is a smelted kimberlite as in \Cref{definitionkimberlite}.
\end{proof}

\subsection{Tubular neighborhoods and a local model correspondence.}
\label{tubularnighborhoodssection}
Let us fix $k$ an algebraically closed field in characteristic $p$ endowed with a fix embedding $\bar{\bbF}_p\subseteq k$.
Fix a pair $(\calD,\Phi_\calD)$ with $\calD$ a $\calG$-torsor over $\on{Spec}(\bbW(k))$ and $\Phi_\calD:\varphi^{*}\calD\to \calD$ an isomorphism over $\on{Spec}(\bbW(k)[\frac{1}{p}])$.
Fix $T\subseteq B\subseteq G_{\bar{\bbQ}_p}$ a maximal torus and a Borel respectively, fix $\mu\in X^+_*({T})$. 
Let $E$ denote the reflex field of $\mu$ with $E_0\subseteq E$ the maximal unramified extension of $\bbQ_p$ in $E$. 
We let $\breve{E}_k=E\otimes_{E_0}\bbW(k)$.

We can define “coordinate-free” versions of some moduli spaces that we studied in the previous sections (\Cref{defi:Grassmannian}, \Cref{thelocalmodeldefi}, \Cref{defi:moduliofshtukas}, \Cref{defi:moduliofshtukasbound}):
\begin{definition}
	\label{defi:coordinatefreegrass}
	We consider functors 
	\[\Gr_{\calD},\Sht_{\calD},\calM^{\leq\mu}_{\calD,O_{\breve{E}_k}},\Sht^{\leq \mu}_{\calD,O_{\breve{E}_k}}: \on{Perf}^\op_k\to \on{Sets}\]
\begin{enumerate}
	\item With $\Gr_{\calD}\Hub{R}=\{(((R^\sharp,\iota),f),\calT,\psi)\}_{/\simeq}$
		where $((R^\sharp,\iota),f)$ is an untilt of $R$ over $\Spa(\bbW(k))$, $\calT$ is a $\calG$-torsor over $\calY_{[0,\infty)}^{R^+}$ and $\psi:\calT\to \calD$ is an isomorphism defined over $\calY_{[0,\infty)}^{R^+}\setminus V(\xi_{R^\sharp})$ that is meromorphic along $\xi_{R^\sharp}$.
	\item With $\Sht_{\calD}\Hub{R}=\{(((R^\sharp,\iota),f),\calT,\Phi,\lambda)\}_{/\simeq}$
		where $((R^\sharp,\iota),f)$ is an untilt of $R$ over $\Spa(\bbW(k))$, $(\calT,\Phi)$ is a shtuka with $\calG$-structure as in \Cref{defi:singleshtukasGstructure}, and $\lambda:\calT\to \calD$ is an isogeny as in \Cref{defi:isogenies}.

	\item We define $\calM^{\leq\mu}_{\calD,O_{\breve{E}_k}}$ and $\Sht^{\leq \mu}_{\calD,O_{\breve{E}_k}}$ by requiring the same boundedness conditions as in \Cref{remark-boundednesscond}. 
		Namely, we require that the following diagrams are Cartesian.
\begin{center}
\begin{tikzcd}
	\calM^{\leq\mu}_{\calD,O_{\breve{E}_k}} \arrow{r} \arrow{d}  & \Gr_\calD \times_{\Spd {\bbW(k)}} \Spd O_{\breve{E}_k} \arrow{d} &
	\Sht^{\leq \mu}_{\calD, O_{\breve{E}_k}} \arrow{r} \arrow{d}  & \Sht_\calD\times_{\Spd \bbW(k)} \Spd O_{\breve{E}_k} \arrow{d} \\
	\on{Hk}_{\calG,O_{{\breve{E}_k}}}^{\leq \mu} \arrow{r} & \on{Hk}_{\calG}\times_{\Spd \bbZ_p}\Spd O_{\breve{E}_k} &
 \on{Hk}_{\calG,O_{{\breve{E}_k}}}^{\leq \mu} \arrow{r} & \on{Hk}_{\calG}\times_{\Spd \bbZ_p}\Spd O_{\breve{E}_k} \\
\end{tikzcd}
\end{center}
\end{enumerate}
\end{definition}

The functors $\Gr_{\calD}$ and $\Sht_{\calD}$ come with canonical sections $c_\calD:\on{Spec}(k)^\dia\to \Gr_{\calD}$ and $c_\calD:\on{Spec}(k)^\dia\to \Sht_\calD$ given by $(\varphi^{*}\calD,\Phi_\calD)$ and $(\calD,\Phi_\calD,\on{Id})$ respectively. 
Fixing an isomorphism $\tau:\calD\simeq \calG$ induces natural isomorphisms 
\[\tau:\Gr_\calD \simeq \Gr_\calG\times_{\Spd \bbZ_p} \Spd \bbW(k),\] and
\[\tau:\calM^{\leq \mu}_{\calD,O_{\breve{E}_k}}\simeq \Grm{O_{\breve{E}}}\times_{\Spd O_{\breve{E}}} \Spd O_{\breve{E}_k}.\]
Analogously, given a $\varphi$-equivariant isomorphism $\tau:\calD\simeq \calE_b$ over $\on{Spec}(\bbW(k)[\frac{1}{p}])$ it induces isomorphisms
\[\tau:\Sht_\calD\simeq \Sht_\calG(b)\times_{\Spd \breve{\bbZ}_p} \Spd \bbW(k) \] and 
\[\tau:\Sht^{\leq \mu}_{\calD,O_{\breve{E}_k}}\simeq \Shtm{O_{\breve{E}}}\times_{\Spd O_{\breve{E}}} \Spd O_{\breve{E}_k}.\]

Moreover, if we are given a section $\sigma: \on{Spec}(k)^\dia\to \Gr_\calG$ (respectively $\sigma: \on{Spec}(k)^\dia\to \Sht_\calG(b)$) we can find $(\calD,\Phi_\calD)$ and an isomorphism $\tau:\calD\simeq \calG$ over $\Spec \bbW(k)$ (respectively $\varphi$-equivariant isomorphism over $\Spec \bbW(k)[\frac{1}{p}]$ $\tau:\calD\to \calE_b$) making the diagrams below commutative
\begin{center}
\begin{tikzcd}

& \Gr_\calD \arrow{dd}{\tau} & 
& \Sht_\calD \arrow{dd}{\tau} \\

\on{Spec}(k)^\dia \arrow{ru}{c_\calD} \arrow{rd}{\sigma} &  &
\on{Spec}(k)^\dia \arrow{ru}{c_\calD} \arrow{rd}{\sigma}\\ 
						 & \Gr_\calG\times_{\Spd \bbZ_p}\Spd \bbW(k)  &

						 & \Sht_\calG(b)\times_{\Spd \breve{\bbZ}_p} \Spd \bbW(k). 
\end{tikzcd}
\end{center}

Indeed, by \Cref{pro:grassmanianreduced}, $\sigma:\on{Spec}(k)^\dia\to \Gr_\calG$ can be represented by data $(\calG,\Psi)$ with $\Psi\in \calG(\bbW(k)[\frac{1}{p}])$, we can then let $\calD=\calG$ and $\Phi_\calD$ be given as $\Phi_\calD:\varphi^*\calG\simeq \calG\xrightarrow{\Psi}\calG$. Here $\tau=\on{Id}_\calG$.

Similarly, by \Cref{pro:adlvasreducedfunctor}, a map $\sigma:\on{Spec}(k)^\dia\to \Sht_\calG(b)$ can be represented by $(\calG,\Phi,\rho)$. In this case we let $\calD=\calG$ and $\Phi_D=\Phi$ and $\tau=\rho$.

Since $k$ is algebraically closed every closed point of $(\Grm{O_{\breve{E}}}\times_{\Spd O_{\breve{E}}} \Spd O_{\breve{E}_k})^\red=\calA_{\calG,\mu,k}$ arises from a section 
\[\sigma:\Spec k^\diamond\to \Grm{O_{\breve{E}}}\times_{\Spd O_{\breve{E}}} \Spd O_{\breve{E}_k},\]
and similarly every closed point of $(\Sht^{\leq \mu}_\calG(b)\times_{\Spd O_{\breve{E}}} \Spd O_{\breve{E}_k})^\red=X^{\leq \mu}_{\calG,k}(b)$ arises from a section
\[\sigma:\Spec k^\diamond\to \Sht^{\leq \mu}_\calG(b)\times_{\Spd O_{\breve{E}}} \Spd O_{\breve{E}_k}.\] 

Constructing $\tau$ and $(\calD,\Phi_\calD)$ as above and passing to formal neighborhoods of closed points \Cref{formalnbbhoosoughtotbedefined} we get identifications
\begin{center}
\begin{tikzcd}
	\Tf{\Sht^{\leq \mu}_\calG(b)\times_{\Spd O_{\breve{E}}} \Spd O_{\breve{E}_k}}{\sigma}\simeq^\tau \Tf{\Sht^{\leq \mu}_{\calD,O_{\breve{E}_k}}}{c_\calD}  &  \Tf{\Grm{O_{\breve{E}}}\times_{\Spd O_{\breve{E}}} \Spd O_{\breve{E}_k}}{\sigma}\simeq^\tau     \Tf{\calM^{\leq \mu}_{\calD,O_{\breve{E}_k}}}{c_\calD}
\end{tikzcd}
\end{center}

This reduces our study of formal neighborhoods at closed points to the study of the spaces
\[\Tf{\calM^{\leq \mu}_{\calD,O_{\breve{E}_k}}}{c_\calD} \text{    and    } \Tf{\Sht^{\leq \mu}_{\calD,O_{\breve{E}_k}}}{c_\calD}\]
as we vary the pair $(\calD,\Phi_\calD)$.
We may think of these spaces as ``deformation spaces'' of the pair $(\calD,\Phi_\calD)$. 
We have the following conjecture.

\begin{conjecture}
	\label{mainconjecture}
With the notation as above, there exists an isomorphism of v-sheaves
\[\Tf{\calM^{\leq \mu}_{\calD,O_{\breve{E}_k}}}{c_\calD} \simeq \Tf{\Sht^{\leq \mu}_{\calD,O_{\breve{E}_k}}}{c_\calD}.\]
\end{conjecture}

\begin{remark}
During the revision process of this article there has been a lot of progress in proving \Cref{mainconjecture} when $\mu$ is assumed to be minuscule. 
Notably, \cite{PR22} for the local abelian type case, \cite{bartling2022mathcalgmudisplayslocalshtuka}, \cite{ito2023deformation} for hyperspecial case and \cite{takaya2024moduli} for the unramified group case. 
\end{remark}

In what follows we formulate and prove \Cref{thm:comparetub} which allows us to bypass \Cref{mainconjecture}. 
We first need to set some notation.
\begin{definition}
	\label{defi:moduliofmatrices}
	We define small v-sheaves as follows
	\[\calGd,\Lgr,\Lgrm,\Lsht,\Lshtm:\on{Perf}^\op_{k}\to \on{Sets}.\]
\begin{enumerate}
	\item With $\calGd\Hub{R}=\{(((R^\sharp,\iota),f),g)\}$ where $((R^\sharp,\iota),f)$ is an isomorphism class of untilts of $R$ over $\Spa \bbW(k)$, $g:\calD\to \calD$ is an automorphism over $\on{Spec}(\bbW(R^+))$ subject to the following condition. 
We require that there exists a pseudo-uniformizer $\varpi_g\in R^+$, depending of $g$, such that the restriction of $g$ to $\on{Spec}(\bbW(R^+)/[\varpi_g])$ is the identity. 
We define $\calGphid$ exchanging the role of $\calD$ for $\varphi^{*}\calD$.
\item With $\Lgr\Hub{R}=\{((R^\sharp,\iota,)f),\calT,\psi,\sigma\}_{/\simeq}$ 
	where $((R^\sharp,\iota),f)$ is an untilt of $R$ over $\Spa \bbW(k)$, $\calT$ is a $\calG$-torsor over $\on{Spec}(\bbW(R^+))$, $\psi:\calT\to \calD$ is an isomorphism over $\on{Spec}(\bbW(R^+)[\frac{1}{\xi}])$ and $\sigma:\calT\to \varphi^{*}\calD$ is an isomorphism over $\on{Spec}(\bbW(R^+))$ subject to the following condition. 
	We require that there is exists a pseudo-uniformizer $\varpi\in R^+$ depending on the data for which $\Phi_\calD\circ \sigma=\psi$ when the data is restricted to $\on{Spec}(\bbW(R^+)/[\varpi])$. 
\item With $\Lsht\Hub{R}=\{((R^\sharp,\iota),f),\calT,\Phi,\lambda,\sigma\}_{/\simeq}$ 
	$((R^\sharp,\iota), f)$ is an untilt of $R$ over $\Spa \bbW(k)$, $\calT$ is a $\calG$-torsor over $\on{Spec}(\bbW(R^+))$, $\Phi:\varphi^{*}\calT\to \calT$ is an isomorphism over $\on{Spec}(\bbW(R^+)[\frac{1}{\xi}])$, $\lambda:\calT\to \calD$ is an isogeny over $\Yint{R^+}{[r,\infty]}$ and $\sigma:\calT\to \calD$ is an isomorphism over $\on{Spec}(\bbW(R^+))$ subject to the following condition. 
	We require that there exists a pseudo-uniformizer $\varpi\in R^+$ depending on the data for which $\sigma=\lambda$ when restricted to $\on{Spec}(B^{R^+}_{[r,\infty]}/[\varpi])$. 
	\item We define $\Lgrm$ or $\Lshtm$ by requiring the same boundedness conditions as in \Cref{remark-boundednesscond}.
		Namely, we require that the following diagrams are Cartesian.
\begin{center}
\begin{tikzcd}
	\Lgrm \arrow{r} \arrow{d}  & \Lgr\times_{\bbW(k)}  \Spd O_{\breve{E}_k} \arrow{d} &
	\Lshtm\arrow{r} \arrow{d}  & \Lsht\times_{\Spd \bbW(k)} \Spd O_{\breve{E}_k} \arrow{d} \\
	\on{Hk}_{\calG,O_{{\breve{E}_k}}}^{\leq \mu} \arrow{r} & \on{Hk}_{\calG}\times_{\Spd \bbZ_p}\Spd O_{\breve{E}_k} &
 \on{Hk}_{\calG,O_{{\breve{E}_k}}}^{\leq \mu} \arrow{r} & \on{Hk}_{\calG}\times_{\Spd \bbZ_p}\Spd O_{\breve{E}_k} \\
\end{tikzcd}
\end{center}
\end{enumerate}
\end{definition}

\begin{theorem}
\label{thm:comparetub}
Given $(\calD,\Phi_\calD)$ and $\mu\in X_*({T}_{\bar{\bbQ}_p})$ as above, and with notation as in \Cref{defi:moduliofmatrices} we have a natural identification 
\[\Lshtm\simeq \Lgrm\]
and consequently we have a correspondence 
\begin{center}
\begin{tikzcd}
&  \Lshtm\simeq \Lgrm\arrow{rd} \arrow{dl}  &  \\
\Tf{\Sht^{\leq \mu}_{\calD,O_{\breve{E}_k}}}{c_\calD} &   & \Tf{\calM^{\leq \mu}_{\calD,O_{\breve{E}_k}}}{c_\calD}
\end{tikzcd}
\end{center}
Moreover, both arrows are $\calGd$-torsors. 
\end{theorem}

\begin{remark}
We emphasize that although the maps $\Lshtm\to \Tf{\Sht^{\leq \mu}_{\calD,O_{\breve{E}_k}}}{c_\calD}$ and $\Lgrm\to \Tf{\calM^{\leq \mu}_{\calD,O_{\breve{E}_k}}}{c_\calD}$ are $\calGd$-torsors, the natural isomorphism $\Lshtm\simeq \Lgrm$ is not equivariant for the $\calGd$-action. 
\end{remark}

Standard arguments using \Cref{pro:descendalongcurve} will prove that the objects in \Cref{defi:moduliofmatrices} are v-sheaves. 
In what follows, we will systematically suppress the untilt of $R$ over $\Spa \bbW(k)$ from the notation. 
There are natural maps $\Lgr\to \Gr_\calD$ and $\Lsht\to \Sht_\calD$ over $\Spd \bbW(k)$. 
After suppressing the untilt the map can be described as follows.
The first one takes a tuple $(\calT,\psi, \sigma)$ and assign $(\calT,\psi)$. 
The second one takes $(\calT,\Phi,\lambda,\sigma)$) and assigns $(\calT,\Phi,\lambda)$. 

The first map is $\calGphid$-equivariant when we endow $\Lgr$  with the action 
\[g\star(\calT,\psi,\sigma)\mapsto (\calT,\psi,g\circ \sigma)\] 
and when we endow $\Gr_\calD$ with the trivial action.

Similarly, the second map is $\calGd$-equivariant when we endow $\Lsht$ with the action 
\[g\star(\calT,\Phi,\lambda,\sigma)\mapsto (\calT,\Phi,\lambda,g\circ \sigma)\]
and where we endow $\Sht_\calD$ with the trivial action.
\Cref{lem:torsorgrassmanian} and \Cref{lem:torsorshtuka} below explain the structure of the maps $\Lgr\to \Gr_\calD$ and $\Lsht\to \Sht_\calD$ below. 

\begin{lemma}
\label{lem:torsorgrassmanian}
The map $\Lgr\to \Gr_\calD$ factors surjectively onto $\Tf{\Gr_\calD}{c_\calD}$. Moreover, the map 
\[\Lgr\to \Tf{\Gr_\calD}{c_\calD}\] is a $\calGphid$-torsor.
\end{lemma}
\begin{proof}
	Observe that if $m:\Spa(A,A^+)\to \Lgr$ is a map with the source being an affinoid perfectoid then $m$ formalizes, see \Cref{defi:formalizing}.
	Indeed, given $(\calT,\psi,\sigma)\in \Lgr(A,A^+)$ and a map $f:\Spa\Hub{B}\to \Spd(A^+,A^+)$, we form $(f^*\calT,f^*\psi,f^*\sigma)$ with $f^*\calT$ defined over $\on{Spec}(\bbW(B^+))$, $f^*\psi:f^*\calT\to \calD$ defined over $\on{Spec}(\bbW(B^+)[\frac{1}{f(\xi)}])$ and $f^*\sigma:f^*\calT\to f^*\varphi^{*}\calD$ defined over $\on{Spec}(\bbW(B^+))$. 
To prove this triple satisfies the constraints, take $\varpi_A\in A^+$ such that after restricting the data to $\on{Spec}(\bbW(A^+)/[\varpi_A])$ the identity $\Phi_\calD\circ \sigma=\psi$ holds. 
Then $f(\varpi_A)$ is topologically nilpotent.
In particular, for any pseudo-uniformizer $\varpi_B\in B^+$ that divides $f(\varpi_A)$ when one restricts the appropriate objects to $\on{Spec}(\bbW(B^+)/[\varpi_B])$ the equation $\Phi_\calD\circ f^*\sigma=f^*\psi$ holds as we wanted to show.

To prove that $\Lgr\to \Gr_\calD$ factors through $\Tf{\Gr_\calD}{c_\calD}$ it suffices to show that given $\Spd(A^+,A^+)\to \Lgr$ described by a triple $(\calT,\psi,\sigma)$ as above the induced map $\Spd(A^+,A^+)^\red\to \Gr_\calD^\red$ factors through the canonical section $c_\calD:\on{Spec}(k)^\diamond\to \Gr_\calD^\red$. 
	Let ${A_\red^+}=A^+/A^{\circ \circ}$. 
	After restricting $(\calT,\psi,\sigma)$ to $\on{Spec}\bbW({A_\red^+})$ we get $\Phi_\calD\circ \sigma =\psi$, and $\on{Spec}({A_\red^+})^\dia\to\Gr_\calD$ is given by $(\calT,\psi)$. 
	Now, one can use $\sigma$ to construct an isomorphism $(\calT,\psi)\simeq(\varphi^{*}\calD,\Phi_\calD)$, so the map factors through $c_\calD:\on{Spec}(k)^\dia\to \Gr_\calD$. 

To prove $\Lgr\to \Tf{\Gr_\calD}{c_\calD}$ is surjective it suffices to lift maps valued on product of points. 
Let $S=\Spa(R,R^+)$ be a product of points \Cref{defi:prodpoints}. 
After suppressing the untilt from the notation, a $S$-valued point of $\Tf{\Gr_\calD}{c_\calD}$ is given by data $(\calT,\psi)$ with $\calT$ defined over $\on{Spec}(\bbW(R^+))$ and $\psi:\calT\to \calD$ defined over $\on{Spec}(\bbW(R^+)[\frac{1}{\xi_{R^\sharp}}])$ such that $(\calT,\psi)$ becomes isomorphic to $(\varphi^{*}\calD, \Phi_\calD)$ when one restricts the data to $\bbW(R_\red^+)$ and $\bbW({R_\red^+})[\frac{1}{p}]$. Indeed, this follows from \Cref{pro:prodpointextendg-tors}.   
The isomorphism ${\sigma_\red}:(\calT,\psi)\to (\varphi^{*}\calD, \Phi_\calD)$ over $\bbW({R_\red^+})$ is unique and given by ${\sigma_\red}=\Phi_\calD^{-1}\circ \psi$ since it has to satisfy the commutative diagram
\begin{center}
\begin{tikzcd}
	
	\calT \arrow{rd}{\psi} \arrow{dd}{{\sigma_\red}}   \\
 & \calD \\
 \varphi^{*}\calD \arrow{ru}{\Phi_\calD}.
\end{tikzcd}
\end{center}
We define $\tilde{\sigma}:=\Phi_\calD^{-1}\circ \psi: \calT\to \varphi^{*}D$ over $\Yint{R^+}{[r,\infty]}$ for $r$ sufficiently big (avoiding $V(\xi)$), clearly $\tilde{\sigma}$ restricts to ${\sigma_\red}$. 
Using \Cref{lem:smallintegraldecomptorsors} we find $\sigma:\calT\to \varphi^{*}\calD$ such that $\sigma=\tilde{\sigma}$ when restricted to $\on{Spec}(B^{R^+}_{[r,\infty]}/[\varpi'])$ for some pseudo-uniformizer $\varpi'\in R^+$. 
In particular $\Phi_\calD\circ \sigma=\psi$ over $\on{Spec}(\bbW(R^+)/[\varpi'])$. 
The data $(\calT,\psi,\sigma)$ defines $\Spa\Hub{R}\to \Lgr$ lifting the original map $\Spa\Hub{R}\to \Tf{\Gr_\calD}{c_\calD}$.

To prove $\Lgr\times_{\Gr_\calD} \Lgr\simeq \calGphid\times_{\Spd \bbW(k)} \Lgr$, take two sets of data $(\calT_i,\psi_i,\sigma_i)$ over $\Spa\Hub{A}$ with $(\calT_1|_{\calY_{[0,\infty)}^{A^+}},\psi_1|_{\calY_{[0,\infty)}^{A^+}\setminus V(\xi)})\simeq (\calT_2|_{\calY_{[0,\infty)}^{A^+}},\psi_2|_{\calY_{[0,\infty)}^{A^+}\setminus V(\xi)})$. 
The isomorphism must be given by $\psi_1^{-1}\circ \psi_2$ and by the full faithfulness of \Cref{thm:puncturedtoAinfked} it extends to $\on{Spec}(\bbW(A^+))$. 
Let 
\[g=\sigma_1\circ \psi_1^{-1}\circ \psi_2\circ \sigma_2^{-1}:\varphi^{*}\calD\to \varphi^{*}\calD.\] 
By hypothesis, $\sigma_i\circ \psi_i^{-1}=\Phi_\calD^{-1}$ on $\on{Spec}(\bbW(A^+)/[\varpi_i])$ for suitable choices of $\varpi_i\in A^+$. 
Consequently, \[(g,\calT_2,\psi_2,\sigma_2)\in [\calGphid\times_{\Spd \bbW(k)} \Lgr]\Hub{A},\] since the construction is functorial we get a map 
$\Lgr\times_{\Gr_\calD} \Lgr\to \calGphid\times_{\Spd \bbW(k)} \Lgr$.
On the other hand, to $(g,\calT,\psi,\sigma)$ we associate the pair of tuples $(\calT,\psi,g\circ\sigma )$ and $(\calT,\psi,\sigma)$. 
These constructions are inverses of each other.
\end{proof}

\begin{lemma}
\label{lem:torsorshtuka}
The map $\Lsht\to \Sht_\calD$ factors surjectively onto $\Tf{\Sht_\calD}{c_\calD}$. Moreover, $\Lsht\to \Tf{\Sht_\calD}{c_\calD}$ is a $\calGd$-torsor.
\end{lemma}
\begin{proof}
That $\Lsht\to \Sht_\calD$ factors surjectively onto $\Tf{\Sht_\calD}{c_\calD}$ follows closely the argument of \Cref{lem:torsorgrassmanian}, and we omit the details.
To prove that $\Lsht\times_{\Sht_\calD} \Lsht\simeq \calGd\times_{\Spd \bbW(k)} \Lsht$, take two sets of data $(\calT_i,\Phi_i,\lambda_i,\sigma_i)$ over $\Spa\Hub{A}$ with $(\calT_1|_{\Yinf{A^+}},\Phi_1|_{\Yinf{A^+}\setminus V(\xi)},\lambda_1)\simeq (\calT_2|_{\Yinf{A^+}},\Phi_2|_{\Yinf{A^+}\setminus V(\xi)},\lambda_2)$. 
The isomorphism must be the unique lift of $\lambda_1^{-1}\circ \lambda_2:\calT_2\to \calT_1$ to $\Yinf{A^+}$. 
Glueing along the $\lambda_i$ and by the fully-faithfulness part of \Cref{thm:puncturedtoAinfked} $\lambda_1^{-1}\circ \lambda_2$ extends to $\on{Spec}(\bbW(A^+))$. 
Moreover, letting 
\[g=\sigma_1\circ \lambda_1^{-1}\circ \lambda_2\circ \sigma_2^{-1}:\calD\to \calD\]
we have $\sigma_1\circ \lambda_1^{-1}=\on{Id}=\lambda_2\circ \sigma_2^{-1}$ over $\on{Spec}(B^{A^+}_{[r,\infty]}/[\varpi_A])$ for suitable $\varpi_A\in A^+$. 
We associate to the original data the tuple $(g,\calT_2,\Phi_2,\lambda_2,\sigma_2)\in [\calGd\times_{\Spd \bbW(k)}\Lsht]\Hub{A}$ giving the left to right map. 
One can construct the inverse using the action map as in the proof of \Cref{lem:torsorgrassmanian}. 
\end{proof}

We can now prove \Cref{thm:comparetub}.

\begin{proof}[Proof. (of \Cref{thm:comparetub})]
Let $\tau=(\varphi)^{-1}$. 
Observe that $\theta:\calGd\to\calGphid$ given by $g\mapsto \varphi^*g$ is an isomorphism with inverse $h\mapsto \tau^*h$. 
Using $\theta$ we can endow $\Lgr$ with a $\calGd$ action, and the projection $\pi:\Lgr\to \Tf{\Gr_\calD}{c_\calD}$ of \Cref{lem:torsorgrassmanian} becomes a $\calGd$-torsor.

We construct an isomorphism $\Theta:\Lgr\to \Lsht$, given on $(A,A^+)$-valued points by $$(\calT,\psi,\sigma)\mapsto (\tau^*\calT,\Phi,\lambda,\tau^{*}\sigma).$$ Here $\Phi:\calT\to \tau^*\calT$ is defined by $\Phi= (\tau^*\sigma)^{-1}\circ \psi$, and $\lambda:\tau^*\calT \to \calD$ is constructed as follows. 
Consider the following (non-commutative!) diagram,
\begin{center}
\begin{tikzcd}

\calT \arrow{r}{\sigma} \arrow{d}{\Phi} & \varphi^*\calD \arrow{d}{\Phi_\calD} \\
\tau^*\calT \arrow{r}{\tau^*\sigma} & \calD
\end{tikzcd}
\end{center}
defined over $\Yint{A^+}{[r,\infty]}$ for big enough $r$ avoiding $V(\xi)$. 

By hypothesis, there is $\varpi\in A^+$ with $\psi=\Phi_\calD\circ \sigma$ over $\on{Spec}(\bbW(R^+)/[\varpi])$. 
Consequently, $\tau^*\sigma\circ\Phi=\Phi_\calD\circ \sigma$ over $\on{Spec}(B^{A^+}_{[r,\infty]}/[\varpi])$. 
By \Cref{lem:liftabilityofiso}, we can construct $\lambda$ as the unique isogeny over $\Yint{A^+}{[r,\infty]}$ lifting $\tau^*\sigma$ with $\lambda=\tau^*\sigma$ over $\on{Spec}(B^{A^+}_{[r,\infty]}/[\varpi])$. 
The uniqueness of $\lambda$ makes this construction functorial so that $\Theta:\Lgr\to \Lsht$ is well-defined. 

The inverse $\Omega=\Theta^{-1}$ is given on $(A,A^+)$-valued points by
\[(\calT,\Phi,\lambda,\sigma)\mapsto (\varphi^*\calT,\sigma\circ \Phi,\varphi^*\sigma).\]
Direct computations show $\Omega\circ \Theta=\on{Id}$, and that $\Theta\circ \Omega(\calT,\Phi,\lambda,\sigma)=(\calT,\Phi,\lambda',\sigma)$ for some $\lambda'$ with $\lambda'=\sigma=\lambda$ over $B^{A^+}_{[r,\infty]}/[\varpi]$. 
The uniqueness part of \Cref{lem:liftabilityofiso} proves $\lambda=\lambda'$ and $\Theta\circ \Omega=\on{Id}$.

One shows directly that both $\Theta$ and $\Omega$ preserve the boundedness conditions so that 
\[\Theta:\Lgrm\to \Lshtm\]
is also an isomorphism. More precisely, we have a commutative diagram.
\begin{center}
\begin{tikzcd}
	\Lgr\arrow[rr, bend left=5, "\Theta"] \arrow{dr}&  & \Lsht \arrow{dl} \arrow[ll, bend left=5, "\Omega"] \\
				  & \on{Hk}_{\calG}\times_{\Spd \bbZ_p}\Spd \bbW(k) &
\end{tikzcd}
\end{center}
and the isomorphism $\Theta:\Lgrm\simeq \Lshtm$ is obtained from the diagram above by pulling back along the closed immersion 
\[ \on{Hk}^{\leq \mu}_{\calG,O_{\breve{E}_k}} \to \on{Hk}_{\calG}\times_{\Spd \bbZ_p}\Spd O_{\breve{E}_k}. \qedhere \]
\end{proof}

In what follows we will apply \Cref{thm:comparetub} to deduce \Cref{thm:specializtheorem}. 
For this we have to understand how local behavior (at formal neighborhoods of closed points) glues to describe global behavior over the rest of our v-sheaf. 
For this we recall some terminology introduced in \cite{Gle22} that try to capture finiteness hypothesis on a v-sheaf.
\begin{definition}
	\label{defi:cJ}
	We say that a locally spatial diamond $X$ is \textit{constructibly Jacobson} if the subset of rank $1$ points are dense for the constructible topology of $|X|$. We refer to them as \textit{cJ-diamonds}.
\end{definition}

\begin{definition} 
\label{defi:kimberlite} 
Let $\calK=(\calF, \calD)$ be a smelted kimberlite.
\begin{enumerate}
	\item We say that $\calK$ is \textit{rich} if $\calD$ is a cJ-diamond, $|\calF^\red|$ is locally Noetherian and $\on{sp}_{\calD}:|\calD|\to |\calF^\red|$ is surjective. 
	\item If $\calK$ is rich we say it is \textit{topologically normal} if for every closed point $x\in |\calF^\red|$ the tubular neighborhood (see \Cref{definitionkimberlite}) ${\calK}^{\circledcirc}_{/x}$ is connected.
	\item We say that a kimberlite $\calF$ is rich (respecitvely topologically normal) if the pair $(\calF,\calF^{\on{an}})$ is rich (respectively topologically normal).
\end{enumerate}
\end{definition}

\begin{theorem}[\cite{AGLR22},\cite{gleason2024tubula}]
	\label{nromalityofshcuberts}
	If $C/E$ is an algebraically closed non-Archimedean field extension, then $\Grm{O_C}$ is rich and topologically normal kimberlite.
\end{theorem}
\begin{proof}
	That $(\Grm{O_C})^{\on{an}}$ is a cJ-diamond follows from \cite[Corollary 4.7]{AGLR22} and \cite[Proposition 4.51]{Gle22}.
	Since $(\Grm{O_C})^\red$ is the perfection of a projective scheme (see \Cref{thm:AGLR}), it is locally Noetherian. 
	To show that $\Grm{O_C}$ is rich it suffices to see that the specialization map is surjective, but this follows from \cite[Proposition 2.34, Proposition 4.14]{AGLR22}.
Topological normality (or unibranchness) is \cite[Theorem 1.3]{gleason2024tubula}.	
\end{proof}

We now show that moduli spaces of $p$-adic shtukas are rich and topologically normal as in \Cref{defi:kimberlite}.

\begin{theorem}
\label{thm:shtukasisOrapian}
The pair $\calK=(\Shtm{O_{\breve{E}}},\Shtm{\breve{E}})$ is a rich and topologically normal smelted kimberlite. 
\end{theorem}
\begin{proof}
\Cref{pro:modulishtukformssmeltedkimberl} proves that this pair is a smelted kimberlite. 
By \Cref{nromalityofshcuberts}, $\Grm{C}$ and consequently $\Grm{\breve{E}}$ are cJ-diamonds, \cite[Proposition 4.46.(1)]{Gle22}.
In \cite[Proposition 23.3.3]{SW20} it is proven that the period morphism $\Shtm{\breve{E}}\to \Grm{\breve{E}}$ is \'etale, so by \cite[Proposition 4.46]{Gle22} $\Shtm{\breve{E}}$ is also a cJ-diamond. 
By \cite[Theorem 1.1]{HV20}, we know that $\adlvdos$ is locally Noetherian. 
To show that $\calK$ is rich it suffices to show that the specialization map is surjective.

By \cite[Lemma 5.23]{Gle22}, we are reduced to showing that for any nonarchimedean field extension $C/\bbW(k)[\frac{1}{p}]$ with $C$ algebraically closed, the specialization map of the base change $\Shtm{O_C}$ is surjective on closed points. 
To show that $\calK$ is rich it is then enough to prove that for any such $C$ the $p$-adic tubular neighborhoods of $\Shtm{O_C}$ are non-empty.
Moreover, to show that $\calK$ is topologically normal it suffices to show that the tubular neighborhoods are connected. 
Non-emptiness and connectedness follow from \Cref{thm:comparetub} and \Cref{nromalityofshcuberts}. 
Indeed, if $k=k_C$ denotes the residue field of $O_C$ we may apply \Cref{thm:comparetub} to compare
\[({\Shtm{O_{\breve{E}_{k}}}})^\circledcirc_{/x}\times_{\Spd \breve{E}_k} \Spd C \text{   with   } ({\Grm{O_{\breve{E}_{k}}}})^\circledcirc_{/y} \times_{\Spd \breve{E}_k} \Spd C\] for some $y$. 

More precisely, we have the following diagram 
\begin{center}
\begin{tikzcd}
\Lshtm\times_{\Spd O_{\breve{E}_k}}\Spd C 	\ar{d}  & \simeq & \Lgrm  \times_{\Spd O_{\breve{E}_k}}\Spd C\ar{d}    \\
	(\Shtm{O_C})^\circledcirc_{/x}\simeq ({\Sht^{\leq \mu}_{\calD,O_C}})^\circledcirc_{/c_\calD}  & & (\calM^{\leq \mu}_{\calD,O_C})^\circledcirc_{/c_\calD}\simeq (\calM^{\leq \mu}_{\calG,O_{C}})^\circledcirc_{/y}
\end{tikzcd}
\end{center}

To show that $(\Shtm{O_{C}})^\circledcirc_{/x}$ is non-empty it suffices to know that $(\calM^{\leq \mu}_{\calG,O_{C}})^\circledcirc_{/y}$ is non-empty, but this directly follows from \Cref{nromalityofshcuberts}.
To show that $(\Shtm{O_{C}})^\circledcirc_{/x}$ is connected it suffices by \Cref{lem:torsorshtuka} to show that $\Lshtm \times_{\Spd O_{\breve{E}_k}} \Spd C$ is connected. 
That $\calM^{\leq \mu}_{\calG,O_C}$ is topologically normal, means that $(\calM^{\leq \mu}_{\calG,O_{C}})^\circledcirc_{/y}$ is connected for all closed points $y$.
Consequently, $(\calM^{\leq \mu}_{\calD,O_{C}})^\circledcirc_{/c_\calD}$ is connected and since  
\[\Lgrm  \times_{\Spd O_{\breve{E}_k}}\Spd \breve{E}_k\to (\calM^{\leq \mu}_{\calD,O_{\breve{E}_{k}}})^\circledcirc_{/c_\calD}\]
is a $\calGd$-torsors and the group $\calGd$ is connected \cite[Proposition 3.4.7]{PR24} we conclude that the space $\Lshtm  \times_{\Spd O_{\breve{E}_k}}\Spd C$ is also connected. 
\end{proof}

We finish this section with the proof of one of our main theorem. 
For the convenience of the reader we recall the notation. 
Let $\calG$ be a parahoric group scheme over $\bbZ_p$, with generic fiber $G$ over $\bbQ_p$. 
Let $b\in \calG(\breve{\bbQ}_p)$ be an element and let $\mu:\bbG_m\to \calG_{\bar{\bbQ}_p}$ be a conjugacy class of geometric cocharacters with reflex field $E$.
Let $\breve{E}$ be the compositum of $E$ and $\breve{\bbQ}_p$ in $\bbC_p$.
\begin{theorem}
\label{thm:specializtheorem}
For every nonarchimedean field extension $F$ of $\breve{E}$ contained in $\bbC_p$ the following hold.  
\begin{enumerate}
\item There is a continuous specialization map
\[\on{sp}:|\Sht_{G,b,\mu,\calG(\bbZ_p)}\times_{\Div^1_E}\times \Spd F|\to |X^{\leq \mu}_\calG(b)|.\]
This map is a specializing and spectral map of locally spectral topological spaces. It is a quotient map and $G_b(\bbQ_p)$-equivariant.
\item The specialization map induces a bijection on connected components 
\[\on{sp}:\pi_0(\Sht_{G,b,\mu,\calG(\bbZ_p)}\times_{\Div^1_E} \Spd F)\to \pi_0(X^{\leq \mu}_\calG(b)).\]
\end{enumerate}
\end{theorem}
\begin{proof}
	By \Cref{thm:shtukasisOrapian} and \Cref{pro:adlvasreducedfunctor} $(\Shtm{O_{{F}}},\Shtm{{F}})$ is a rich smelted kimberlite with reduction $X_\calG^{\leq \mu}(b)$. 
	This implies by \cite[Theorem 4.40]{Gle22} that the specialization map 
\[\on{sp}:|\Sht_{G,b,\mu,\calG(\bbZ_p)}\times_{\Div^1_E} \Spd F|\to |X^{\leq \mu}_\calG(b)|,\]
	 is a specializing spectral map of locally spectral spaces. 
Moreover, by \cite[Lemma 4.53]{Gle22} it is a quotient map. 
Now, $G_b(\bbQ_p)$-equivariance follows from functoriality of the specialization map. 
Since $\Shtm{O_{F}}$ is rich and topologically normal we can apply \cite[Proposition 4.55]{Gle22} to show that the specialization map induces a bijection on connected components. 
\end{proof}

\section{Geometric connected components of moduli spaces of shtukas at infinite level.}
\label{thiisghesecondsection}
\subsection{Moduli spaces of $p$-adic shtukas at infinite level.}
\label{Moduliofshtukassection}
\subsubsection{A conjecture of Rapoport and Viehmann.}
\label{sectionoftheconjecture}
	As we mentioned in the introduction, given $(G,b,\mu)$ a local Shimura datum \cite[Definition 5.1]{RV14} with reflex field $E$, Rapoport and Viehmann conjectured the existence of a tower $\bbM_K$ of smooth rigid-analytic spaces over $\breve{E}$ indexed by the compact open subgroups $K\subseteq G(\bbQ_p)$ whose cohomology should realize instances of the local Langlands correspondence \cite[$\mathsection$ 7]{RV14}. 

	Inspired by the works of de Jong \cite{dJ95}, of Strauch \cite{MR2441699} and of Chen \cite{Chen}, Rapoport and Viehmann make a general conjecture describing the geometric connected components of $\bbM_K$.
Let $\bbM^{\on{geo}}_K$ denote $\bbM_K\times_{\on{Spa} \breve{E}}\on{Spa}\bbC_p$.
Given $(G,b,\mu)$ a local Shimura datum, we let $(G^\ab,b^\ab,\mu^\ab)$ denote the local Shimura datum induced by the determinant map $\det:G\to G^\ab=G/G^\der$.
In other words, $b^\ab=\on{det}(b)$ and $\mu^\ab=\on{det}\circ \mu$. 
We let $\bbM^{\ab,\on{geo}}_{L}$ denote the tower of rigid spaces over $\bbC_p$ attached to $(G^\ab,b^\ab,\mu^\ab)$ and indexed by compact open subgroups of $L\subseteq G^\ab(\bbQ_p)$ \cite{Che13} \cite[$\mathsection$ 5.2]{RV14}.

\begin{conjecture}[{\cite[Conjecture 4.30]{RV14}}]
	\label{conj-rapoport-viehmann}
	Assume that $G^\der$ is simply connected. Then there is a morphism of towers indexed by compact open subgroups $K\subseteq G(\bbQ_p)$ 
	\[\bbM^{\on{geo}}_K\to \bbM^{\on{ab},\on{geo}}_{\on{det}(K)}\]
	compatible with the group action on the towers.
	Furthermore, whenever $(b,\mu)$ is HN-irreducible then this morphism induces bijections
	\[\pi_0(\bbM^{\on{geo}}_K)\simeq \pi_0(\bbM^{\ab,\on{geo}}_{\on{det}(K)})\simeq G^\ab(\bbQ_p)/\det(K). \]
\end{conjecture}

\begin{remark}
	Since $G^\ab$ is a torus, the local Shimura variety $\bbM^\ab_L$ is a tower of $0$-dimensional rigid spaces over $\Spa \breve{E}$. Consequently, $\bbM^{\ab,\on{geo}}_{\on{det}(K)}\simeq \coprod_{G^\ab(\bbQ_p)/\det(K)} \Spa \bbC_p$. 
\end{remark}

The existence of the tower was proven to exist by Scholze and Weinstein \cite{SW13} \cite[$\mathsection$ 24.1]{SW20} by letting $\bbM_K$ denote the unique smooth rigid-analytic space representing $\Sht_{G,b,\mu,K}$ as in $\mathsection$ \ref{section:finitelevelshtukas} (see \cite[Definition 24.1.3]{SW20}).  
From this perspective, the functoriality of towers with respect to the morphism $\det:G\to G^\ab$ together with the compatibility with group actions becomes evident.
So the only thing left to prove is that when $(G,b,\mu)$ is HN-irreducible \Cref{defi:BGmu} the following identity holds
\begin{equation}
	\label{RZbijection}
\pi_0(\bbM^{\on{geo}}_K)\simeq \pi_0(\bbM^{\ab,\on{geo}}_{\on{det}(K)}).
\end{equation}

The purpose of this section is to show that \Cref{RZbijection} holds whenever $G$ is an unramified group over $\bbQ_p$ and $(G,b,\mu)$ is HN-irreducible. This shows that \Cref{conj-rapoport-viehmann} holds in the unramified case.  

The sketch of the proof is as follows. 
Since by definition $(\bbM^{\on{geo}}_K)^\lozenge=\Sht^{\on{geo}}_{G,b,\mu,K}$, it suffices by \cite[Lemma 15.6]{Sch17} to compute the connected components of the latter. 
In turn, it is more convenient to work with the infinite level moduli space of shtukas $\Sht^{\on{geo}}_{G,b,\mu,\infty}$. 
Our main result in this section \Cref{thm2:mainTheorem} computes $\pi_0(\Sht^{\on{geo}}_{G,b,\mu,\infty})$ as a right $G(\bbQ_p)\times G_b(\bbQ_p)\times W_E$-set whenever $G$ is a unramified group and $(G,b,\mu)$ is HN-irreducible. 
We use \Cref{thm:maintheorem} as a stepping stone which is already enough to show \Cref{conj-rapoport-viehmann} (see \Cref{corollaryRVconj}).
The proof of \Cref{thm:maintheorem} is heavily inspired by the arguments in \cite{Chen} except that we avoid the use of de Jong's fundamental group and we also crucially use \Cref{thm:specializtheorem} instead of de Jong's theorem on normality.

For future reference, the first few subsections review the geometry of moduli spaces of $p$-adic shtukas at infinite level. 

	\subsubsection{Divisors and $G$-bundles on the Fargues--Fontaine curve.}
	We start by recalling how to discuss the Weil group action.
	For $S\in \Perf$ we let ${X}_{\on{FF},S}$ denote the relative Fargues--Fontaine curve over $S$ \cite[Definition II.1.15]{FS21}. 
	Given $S\in\Perf$ and $S^\sharp$ an untilt of $S$ we obtain a closed Cartier divisor $S^\sharp \to {X}_{\on{FF},S}$ \cite[Definition 5.3.2]{SW20}, \cite[Proposition II.1.18]{FS21}. 
	Given two vector bundles $\calE_1$ and $\calE_2$ defined over $X_{\on{FF},S}$ and an untilt $S^\sharp$ of $S$, we refer to isomorphisms \[\alpha:\calE_1\dashrightarrow \calE_2\] 
	defined over ${X}_{\on{FF},S}\setminus S^\sharp$ and meromorphic along $S^\sharp$ \cite[Definition 5.3.5]{SW20} as \textit{modifications}.
	We let
	\[\Div_{\bbQ_p}^1:\Perff^\op\to \Sets\]
	denote the moduli space of degree $1$ divisors on the relative Fargues--Fontaine curve \cite[Definition II.1.19]{FS21}. 
	Recall that we have an identity $\Spd \breve{\bbQ}_p=\Spd \bbQ_p\times \Spd \bar{\bbF}_p$.
	Here the projection map to $\Spd \bar{\bbF}_p$ is given by the formula
	\[((S^\sharp,i),f)\mapsto [f^\flat\circ i^{-1}:S\to \Spd \bar{\bbF}_p].\]
	Here, when $S^\sharp=\Spa(R^\sharp,R^{\sharp,+})$, $f^\flat$ denotes the unique map $f^\flat:(S^\sharp)^\flat\to \Spd \bar{\bbF}_p$ inducing the following commutative diagram
	\begin{center}
	\begin{tikzcd}
		S^\sharp	\arrow{r}{\infty} \arrow{rd}{f}  & \calY_{(0,\infty)}^{(R^{\sharp,+})^\flat} \arrow{d}{\bbW(f^\flat)} \\
						      & \Spa \bbW(\bar{\bbF}_p)[\frac{1}{p}].
	\end{tikzcd}
	\end{center}
	Recall that the ring map $\phi:\breve{\bbQ}_p\to \breve{\bbQ}_p$ induces an isomorphism $\Spd(\phi):\Spd \breve{\bbQ}_p\to \Spd \breve{\bbQ}_p$ which we denote by $\varphi$.
	One can verify that 
	\[\varphi=\on{Id}_{\Spd \bbQ_p}\times \on{Frob}_{ \Spd \bar{\bbF}_p}\] as automorphisms of $\Spd \bbQ_p\times \Spd \bar{\bbF}_p$. 
	We let $\tau=\on{Frob}^{-1}_{\Spd \breve{\bbQ}_p}\circ \varphi:\Spd \breve{\bbQ}_p\to \Spd \breve{\bbQ}_p$, then 
	\[\tau=\Frob^{-1}_{\Spd \bbQ_p}\times \on{Id}_{\Spd \bar{\bbF}_p}\]
	is an automorphism of $\Spd \breve{\bbQ}_p$ over $\Spd \bar{\bbF}_p$. 
	This gives a right action of $\tau^\bbZ$ on $\Spd \breve{\bbQ}_p$ which we can make more explicit. 
	It is
	\[\Spd \breve{\bbQ}_p \times \tau^\bbZ \to \Spd \breve{\bbQ}_p\]
	with formula
	\[((S^\sharp,i),f)\star \tau^s \mapsto ((S^\sharp,\Frob_S^{s}\circ i), \varphi^{s} \circ f  ).\]
	We observe that the projection map $\pi_{\Spd \bar{\bbF}_p}:\Spd \breve{\bbQ}_p\to \Spd \bar{\bbF}_p$ takes the shape 
	\[\pi_{\Spd \bar{\bbF}_p}(((S^\sharp,i),f)\star \tau^s)=\pi_{\Spd \bar{\bbF}_p}((S^\sharp,\Frob_S^{s}\circ i),\varphi^{s}\circ f)=\]
	\[(\varphi^s\circ f)^\flat \circ i^{-1}\circ \Frob_S^{-s}=(\Frob_{\Spd \bar{\bbF}_p}^s\circ f^\flat) \circ i^{-1} \circ \Frob_S^{-s}= f^\flat \circ i^{-1}.\]
	Let $\on{deg}:W_{\bbQ_p}\to \bbZ$ be the degree map with $\on{deg}(w)=s$ if $w$ acts by $\phi^s$ on $\breve{\bbQ}_p$.
	We can endow $\Spd \bbC_p$ with a ``modified'' right action that respects the structure map to $\Spd \bar{\bbF}_p$, it is
	\[\Spd \bbC_p \times W_{\bbQ_p} \to \Spd \bbC_p\]
	with formula
	\begin{equation}
		\label{action-spdcp}
	((S^\sharp,i),f)\star w \mapsto ((S^\sharp,\Frob_S^{\deg(w)}\circ i),\Spa(w)\circ f).
	\end{equation}
We note that the modified action agrees with the standard action from \Cref{sec:Notation} when we restrict to the inertia subgroup $I_E\subseteq W_E$. 
	Let $r:\Spa \bbC_p\to \Spa \breve{\bbQ}_p$ be the natural map induced from the inclusion $\breve{\bbQ}_p\subseteq \bbC_p$.
	We can compute the projection map $\pi_{\Spd \bar{\bbF}_p}:\Spd \bbC_p\to \Spd \bar{\bbF}_p$ to be given by the formula
	\[((S^\sharp,i),f)\mapsto [(r\circ f)^\flat\circ i^{-1}:S\to \Spd \bar{\bbF}_p].\]
	With this formula we can compute 
	\[\pi_{\Spd \bar{\bbF}_p}(((S^\sharp,i),f)\star w)=\pi_{\Spd \bar{\bbF}_p}((S^\sharp,\Frob_S^{\on{deg}(w)}\circ i),\Spa(w)\circ f)=(r\circ \Spa(w)\circ f)^\flat \circ i^{-1}\circ \Frob_S^{-{\on{deg}}(w)}=\]
	\[ (\varphi_{\breve{\bbQ}_p}^{\on{deg}(w)}\circ r\circ f)^\flat \circ i^{-1}\circ \Frob_S^{-{\on{deg}}(w)}=
	\Frob_{\Spd \bar{\bbF}_p}^{\on{deg}(w)}\circ (r\circ f)^\flat \circ i^{-1} \circ \Frob_S^{-\on{deg}(w)}= (r\circ f)^\flat \circ i^{-1}.\]
	This computation shows concretely that the action defined above respects the structure map $\Spd \bbC_p\to \Spd \bar{\bbF}_p$.
	
	We have identifications
	\[\Div_{\bbQ_p}^1=[\Spd \breve{\bbQ}_p/\tau^\bbZ]=[\Spd \bbC_p/\underline{W_{\bbQ_p}}].\]
	This construction is the same as \cite[$\mathsection$ IV.7]{FS21}.

	More generally, let $E\subseteq \bbC_p$ be a finite field extension of $\bbQ_p$ with absolute Galois group $\Gamma_E\subseteq \Gamma_{\bbQ_p}$ and Weil group $W_E=\Gamma_E\cap W_{\bbQ_p}$. 
	Let $\bbQ_{p^s}\subseteq E$ denote the maximal unramified extension.
	We let $\breve{E}\subseteq \bbC_p$ denote the compositum $E\cdot \breve{\bbQ}_p$, the natural map $E\otimes_{\bbQ_{p^s}}\breve{\bbQ}_p\to \bbC_p$ induces an isomorphism onto $\breve{E}$. 
	We let $\phi_{\breve{E}}:\breve{E}\to \breve{E}$ be the automorphism conjugate to $\on{Id}\otimes \phi^s$.
	This induces an automorphism $\varphi_{\breve{E}}:\Spd \breve{E}\to \Spd \breve{E}$, and we have a commutative diagram 
\begin{center}
\begin{tikzcd}
	\Spd \breve{E} \arrow{d} \arrow{r}{\varphi_{\breve{E}}}  & \Spd \breve{E} \arrow{d} \\
	\Spd \bar{\bbF}_p \arrow{r}{\Frob^s} & \Spd \bar{\bbF}_p
\end{tikzcd}
\end{center}
We let $\tau_{\breve{E}}=\Frob^{-s}_{\Spd \breve{E}}\circ \varphi_{\breve{E}}$, then $\tau_{\breve{E}}$ is an automorphism of $\Spd \breve{E}$ over $\Spd \bar{\bbF}_p$.
If $\Div^1_E:\Perff^\op\to \Sets$ denotes the space of degree $1$ divisors on the relative Fargues--Fontaine curve attached to $E$ then 
\[\Div_E^1=[\Spd \breve{E}/\tau_{\breve{E}}^\bbZ]=[\Spd \bbC_p/\underline{W_{E}}].\]
There is a natural map $\Div^1_E\to \Div^1_{\bbQ_p}$, if we have $D_E\in \Div_E^1(S)$ we will denote by $D_{\bbQ_p}$ its image in $\Div^1_{\bbQ_p}(S)$. 
Note that a degree $1$ divisor $D\in \Div_{\bbQ_p}^1(S)$ gives rise to a perfectoid space $D^\flat$ that is isomorphic to $S$, but $D$ is not an untilt of $S$ since the data $D\subseteq {X}_{\on{FF},S}$ does not specify an isomorphism between $D^\flat$ and $S$. 
Similarly, fix $D_{\bbQ_p}\in \Div_{\bbQ_p}^1(S)$, this is a perfectoid space over $\Spa \bbQ_p$. 
The set of $D_E\in \Div^1_E(S)$ mapping to $D_{\bbQ_p}$ under the map $\Div^1_E(S)\to \Div_{\bbQ_p}^1(S)$ correspond to the different ways of endowing $D_{\bbQ_p}$ with the structure of an adic space over $\Spa E$.

\subsubsection{The group $G_b(\bbQ_p)$}
\label{subsection-Gb}
Recall that $G$ denotes a connected reductive group over $\bbQ_p$.
Recall that given $b\in G(\breve{\bbQ}_p)$ one can define another reductive group $G_b$ (or $J$ in the notation of \cite[$\mathsection 3.3$]{Kot97}) over $\Spec \bbQ_p$ that satisfies the following property. 
If we denote by $V_b$ the isocrystal with $G$-structure attached to $b$, then $G_b$ satisfies    
\[G_b(\bbQ_p)=\on{Aut}(V_b).\]

Fargues--Scholze consider a more geometric version of $G_b$.
When $V$ is an isocrystal with $G$-structure and $S\in \Perff$ one can construct a $G$-bundle over $X_{\on{FF},S}$ that we denote $\calE_{V,S}$ as in \cite[$\mathsection$ III.2.1]{FS21}.
	When $V=V_b$ and if $S$ is clear from the context we simply write $\calE_b$ for $\calE_{V_b,S}$.
	Following \cite[$\mathsection$ III.5.1]{FS21}, we let $\widetilde{G}_b$ denote the functor 
	\[\widetilde{G}_b:\Perff^\op\to \on{Groups}\] given by the formula
	\[\widetilde{G}_b(S)=\on{Aut}_{X_{\on{FF,S}}}(\calE_b).\]
	There is a v-surjective map of groups $\widetilde{G}_b\to \underline{G_b(\bbQ_p)}$, and if $\widetilde{G}_b^U$ denotes the kernel then $\widetilde{G}_b= \widetilde{G}_b^U\rtimes \underline{G_b(\bbQ_p)}$ \cite[Proposition III.5.1]{FS21}. 
	The inclusion $\underline{G_b(\bbQ_p)}\subseteq \widetilde{G}_b$ is the subgroup of automorphisms of $\calE_b$ that come from an automorphism of $V_b$ \cite[$\mathsection$ III.5.1, Proposition III.4.7]{FS21}.
	Whenever $b$ is basic the identity $\widetilde{G}_b= \underline{G_b(\bbQ_p)}$ holds and when $b=\mathds{1}$ we have an identification $G_{\mathds{1}}=G$.
	In particular, $\widetilde{G}_\mathds{1}=\underline{G(\bbQ_p)}$ \cite[$\mathsection$ III.4, Proposition III.4.2]{FS21}.

\subsubsection{Definition of $\Sht^{\on{geo}}_{G,b,\mu,\infty}$}
\label{some-section-defining-shtukas}
	We fix $(G,b,\mu)$ a $p$-adic shtuka datum over $\bbQ_p$ \cite[Definition 5.1]{RV14}, \cite[Definition 23.1.1]{SW20}. 
	That is, $G$ is a connected reductive group over $\bbQ_p$, $b$ is an element of the Kottwitz set $B(G)$ \cite{Kot97}, and $\mu$ is a conjugacy class of geometric cocharacters $\mu\in\{\bbG_m\to G_{\bar{\bbQ}_p}\}/\sim$.
	From now on $E/\bbQ_p$ will denote the reflex field of $\mu$ i.e. the smallest Galois extension of $\bbQ_p$ for which $\Gamma_E$ preserves the conjugacy class $\mu$.

	Let $S=\Spa C$ with $C$ an algebraically closed field. 
	Let $D_E\in \Div_E^1(S)$ and let $\alpha:\calE_1\dashrightarrow \calE_2$ be a modification of $G$-bundles defined over $X_{\on{FF},S}\setminus D_{\bbQ_p}$.
	In what follows we will explain what it means for $\alpha$ to be bounded by $\mu$.
	The completion of $X_{\on{FF},S}$ along $D_{\bbQ_p}$ gives rise to a discrete valuation ring $B^+_{D_{\bbQ_p}}$ with fraction field $B_{D_{\bbQ_p}}$ and residue field $C_D$ \cite[Definition 12.4.3]{SW20} \cite[$\mathsection$ VI.1]{FS21}. 
	The ring $B^+_{D_{\bbQ_p}}$ is non-canonically isomorphic to $C_D\pot{ \xi}$, and the choice of an untilt of $C^\sharp$ of $C$ projecting to $D_{\bbQ_p}$ induces an isomorphism $B^+_{D_{\bbQ_p}}\simeq B^+_{\on{dR}}(C^\sharp)$ where the latter is the de-Rham period ring of Fontaine (see \cite[Definition 12.4.3]{SW20}).

	By Beauville-Laszlo glueing \cite[Lemma 5.2.9]{SW20}, the modification $\alpha$ gives rise to a well-defined element of 
	\[\alpha_D\in G(B^+_{D_{\bbQ_p}})\backslash G(B_{D_{\bbQ_p}}) /G(B^+_{D_{\bbQ_p}})= G(C_D\pot{ \xi})\backslash G(C_D\rpot{ \xi}) /G(C_D\pot{ \xi}).\]
	After fixing a choice $T\subseteq B\subseteq G_{C_D}$ of $T$ a maximal torus and $B$ a Borel, it follows from the Cartan decomposition that 
	\[ X_*^+(T)\simeq G(C_D\pot{ \xi})\backslash G(C_D\rpot{ \xi}) /G(C_D\pot{ \xi})\]
	where $X_*^+(T)$ is the set of dominant cocharacters. 
	Since we started with data $D_E\in \Div_E^1(S)$, $C_D$ comes equipped with $E$-structure. 
	Since $C_D$ is algebraically closed we can find an $E$-linear embedding $\iota:\bar{\bbQ}_p\to C_D$. 
	The conjugacy class $\iota(\mu)\in \{\bbG_m\to G_{C_D}\}/\sim$ gives rise to a unique element $\mu_0\in X_*^+(T)$ belonging to this conjugacy class. 
	Moreover, since $E$ is the reflex field of $\mu$ the formation of $\mu_0$ does not depend on $\iota$. 
	We say that $\alpha$ is \textit{bounded} by $\mu$ if $\alpha_D\leq \mu_0$ in the Bruhat order \cite[Definition 19.2.]{SW20}.
	More generally, if $S\in \Perf$, $D_E\in \Div^1_E(S)$ and $\alpha:\calE_1\dashrightarrow \calE_2$ is modification, then we say that $\alpha$ is \textit{bounded} by $\mu$ if for every geometric point of $\bar{s}\to S$ the induced modification $\alpha_{\bar{s}}:\calE_{1,\bar{s}}\dashrightarrow \calE_{2,\bar{s}}$ is bounded by $\mu$. 
	Our main object of study in this section is $\Sht^{\on{geo}}_{G,b,\mu,\infty}$ which is defined as follows (see \cite[$\mathsection$ 23.2]{SW20}).

	\begin{definition}
		\label{infinitelevelmoduli}
	We consider the \textit{moduli space of $p$-adic shtukas at infinite level} of type $(G,b,\mu)$,
	as a functor
	\[\Sht_{G,b,\mu,\infty}:\Perf^\op_{\bar{\bbF}_p} \to \Sets,\]
	that is given by the formula
	\[\Sht_{G,b,\mu,\infty}(S)=\{(D_E,\alpha)\}\]
	where $D_E\in \Div^1_E(S)$ and $\alpha$ is a modification  
	of $G$-bundles 
	\[\alpha:\calE_\mathds{1}\dashrightarrow \calE_b\]
	over $X_{FF,S}\setminus D_{\bbQ_p}$ meromorphic along $D_{\bbQ_p}$ that is bounded by $\mu$.
	We let 
\[\Sht^{\on{geo}}_{G,b,\mu,\infty}:=\Sht_{G,b,\mu,\infty}\times_{\Div^1_E}\Spd \bbC_p.\]
	\end{definition}

	\subsubsection{Definining the three actions.}
	\label{definingsection}
Let us make the right action of $\underline{G(\bbQ_p)}\times \underline{G_b(\bbQ_p)}\times \underline{W_E}$ on $\Sht^{\on{geo}}_{G,b,\mu,\infty}$ explicit.
An $S$-point of $\Sht^{\on{geo}}_{G,b,\mu,\infty}$ is given by data $(((S^\sharp,i),f),\alpha)$ with $((S^\sharp,i),f)$ an untilt over $\Spa \bbC_p$ and $\alpha:\calE_\mathds{1}\dashrightarrow \calE_b$ a modification over $X_{FF,S}\setminus S^\sharp$ meromorphic along $S^\sharp$ and bounded by $\mu$.
	There is an evident left action map \[\widetilde{G}_b(S)\times \Sht^{\on{geo}}_{G,b,\mu,\infty}(S)\to \Sht^{\on{geo}}_{G,b,\mu,\infty}(S)\]
	with formula 
	\[g_b\star ((S^\sharp,i),f),\alpha)\mapsto ((S^\sharp,i),f),g_b\circ \alpha)\]
	which we can restrict along the inclusion $\underline{G_b(\bbQ_p)}(S)\subseteq \widetilde{G}_b(S)$ and turn it into a right action with formula
	\[((S^\sharp,i),f),\alpha)\star g_b\mapsto ((S^\sharp,i),f),g^{-1}_b\circ \alpha).\]

	There is an evident right action map 
	\[\Sht^{\on{geo}}_{G,b,\mu,\infty}(S)\times \underline{G(\bbQ_p)}(S)\to \Sht^{\on{geo}}_{G,b,\mu,\infty}(S)\]
	with formula 
	\[((S^\sharp,i),f),\alpha)\star g_1\mapsto ((S^\sharp,i),f),\alpha\circ g_1).\]
Finally, the right action 
\[ \Sht^{\on{geo}}_{G,b,\mu,\infty}(S) \times  \underline{W_E}(S) \to \Sht^{\on{geo}}_{G,b,\mu,\infty}(S)\]
is the one obtained from \Cref{action-spdcp}.

From the explicit formula, we see that the actions of $\underline{G(\bbQ_p)}$, $\widetilde{G}_b(S)$ and $\underline{W_E}$ commute so we can gather them to define an action of $\underline{G(\bbQ_p)}\times \underline{G_b(\bbQ_p)} \times \underline{W_E}$ on $\Sht^{\on{geo}}_{G,b,\mu,\infty}$. 
Since $G(\bbQ_p)\times G_b(\bbQ_p) \times {W_E}$ is a locally profinite group we have an identification of topological spaces
\[|\Sht^{\on{geo}}_{G,b,\mu,\infty}\times \underline{G(\bbQ_p)}\times \underline{G_b(\bbQ_p)} \times \underline{W_E} |\simeq |\Sht^{\on{geo}}_{G,b,\mu,\infty}| \times G(\bbQ_p)\times G_b(\bbQ_p)\times W_E    \] 
and the action map 
\[|\Sht^{\on{geo}}_{G,b,\mu,\infty}| \times G(\bbQ_p)\times G_b(\bbQ_p)\times W_E  \to |\Sht^{\on{geo}}_{G,b,\mu,\infty}|    \] 
is continuous since the association of a topological space to a small v-sheaf is functorial \cite[Proposition 12.7]{Sch17}.
\begin{remark}
	\label{functoriality-wtih-group}
	One can see from the definition of $\Sht^{\on{geo}}_{G,b,\mu,\infty}$, and the definition of the right group action by $\underline{G(\bbQ_p)}\times \underline{G_b(\bbQ_p)} \times \underline{W_E}$ that this construction is functorial in the triple $(G,b,\mu)$.  
	More precisely, let $f:G\to H$ be a map of reductive groups, let $b_H=f(b)$, let $\mu_H=f\circ \mu$ and let $F\subseteq E$ denote the reflex field of $\mu_H$.
	We have canonical group maps 
\[\underline{G(\bbQ_p)}\times \underline{G_b(\bbQ_p)} \times \underline{W_E}\to \underline{H(\bbQ_p)}\times \underline{H_{b_H}(\bbQ_p)} \times \underline{W_F}.\]
In this way, pre-composition endows $\Sht^{\on{geo}}_{H,b_H,\mu_H,\infty}$ with a right action by $\underline{G(\bbQ_p)}\times \underline{G_b(\bbQ_p)} \times \underline{W_E}$. 
Moreover, the canonical map
\[\Sht^{\on{geo}}_{G,b,\mu,\infty}\to \Sht^{\on{geo}}_{H,b_H,\mu_H,\infty}\]
obtained by sending a $G$-torsor to its associated $H$-torsor is equivariant for the $\underline{G(\bbQ_p)}\times \underline{G_b(\bbQ_p)} \times \underline{W_E}$-action that we just defined.
\end{remark}

\subsubsection{Moduli spaces of $p$-adic shtukas at finite level.} 
\label{section:finitelevelshtukas}
For future reference we record here the relation between $\Shtm{O_{\breve{E}}}$ (defined and studied in \Cref{sectionspecialization}) and $\Sht^{\on{geo}}_{G,b,\mu,\infty}$ \Cref{infinitelevelmoduli}. 
We also set some notation.

For any compact open subgroup $K\subseteq G(\bbQ_p)$ we let
\[\Sht_{G,b,\mu,K}:=\Sht_{G,b,\mu,\infty}/\underline{K} \text{  and  } \Sht^{\on{geo}}_{G,b,\mu,K}:=\Sht^{\on{geo}}_{G,b,\mu,\infty}/\underline{K} .\]
	This notation is justified by \cite[Definition 23.1.1, Proposition 23.3.1]{SW20}. 
For field extensions $\breve{E}\subseteq F\subseteq \bbC_p$ we let 
\[\Sht^F_{G,b,\mu,K}=\Sht_{G,b,\mu,K}\times_{\Div^1_E} \Spd F.\] 
We will use the following limit formulas (see \cite[$\mathsection$ 23.3]{SW20})
\[\Sht_{G,b,\mu,\infty}^{\on{geo}}=\varprojlim_{\breve{E}\subseteq F\subseteq \bbC_p} \Sht_{G,b,\mu,\infty}^{F} \text{  and  } \Sht_{G,b,\mu,\infty}^{F}=\varprojlim_{K\subseteq G(\bbQ_p)} \Sht_{G,b,\mu,\infty}^{F}.\]
Recall that throughout the text $\calG$ denotes a parahoric model of $G$. 
When $K=\calG(\bbZ_p)$ we will crucially use the following identities which can be deduced from \cite[Proposition 23.3.1]{SW20} 
\[\Sht_{G,b,\mu,\calG(\bbZ_p)}\times_{\Div^1_E}\Spd \breve{E}=\Shtm{O_{\breve{E}}}\times_{\Spd O_{\breve{E}}} \Spd \breve{E}\] 
and 
\[\Sht^F_{G,b,\mu,\calG(\bbZ_p)}=\Shtm{O_{\breve{E}}}\times_{\Spd O_{\breve{E}}} \Spd F.\]

\subsubsection{The period morphism.}
Recall the Grothendieck--Messing period map \cite[\S 23.3]{SW20} 
\[\pi_{\on{GM}}:\Sht_{G,b,\mu,\infty}\to \Gr^{\leq \mu}_{G,\Div^1_E}\] 
defined over $\Div_E^1$. 
Here we use Beauville--Laszlo glueing to identify $\Gr^{\leq \mu}_{G,\Div^1_E}$ with the moduli space classifying tuples 
\[S\mapsto \{(D_E,\calE,\alpha)\}_{/\simeq}\]
where $D_E\in \Div^1_E(S)$, $\calE$ is a $G$-bundle on the relative Fargues--Fontaine curve $X_{\on{FF},S}$ and 
\[\alpha:\calE\dashrightarrow \calE_b\]
is a modification.
For $(D_E,\alpha)\in \Sht_{G,b,\mu,\infty}(S)$ as in \Cref{infinitelevelmoduli} the formula for $\pi_{\on{GM}}(D_E,\alpha)$ is 
\[\pi_{\on{GM}}(D_E,\alpha)=(D_E,\calE_\mathds{1},\alpha).\]
Its image is an open subset $\Gr^{\leq \mu,b-\on{adm}}_{G,\Div^1_E}\subseteq \Gr^{\leq \mu}_{G,\Div^1_E}$ \cite[Theorem 22.6.2]{SW20} 
and the map
\begin{center}
	\begin{align}
		\pi_{\on{GM}}:\Sht_{G,b,\mu,\infty}\to \Gr^{\leq \mu,b-\on{adm}}_{G,\Div^1_E}
	\end{align}
\end{center}
is a $\underline{G(\bbQ_p)}$-torsor \cite[$\mathsection$ 23.2, Theorem 22.5.2]{SW20} for the pro-\'etale topology.
For a given tuple $(D_E,\calE,\alpha)\in \Gr^{\leq \mu,b-\on{adm}}_{G,\Div^1_E}(S)$ one can describe the fiber $\pi_{\on{GM}}^{-1}(S)$ as the moduli space of isomorphisms $\tau:\calE_{\mathds{1}}\xrightarrow{\simeq} \calE$.

\begin{remark}
	\label{torsors-vs-galoisreps}
	Recall that any pro-\'etale $\underline{G(\bbQ_p)}$-torsor over $\Spd \bbC_p$ is trivial. 
	Indeed, by \cite[Lemma 10.13]{Sch17} the torsor is representable by a perfectoid space and since $\Spd \bbC_p$ is strictly totally disconnected, by \cite[Corollary 7.22]{Sch17} $\Spd \bbC_p$ splits every pro-\'etale cover.
In particular, if $\breve{E}\subseteq F \subseteq \bbC_p$ is a $p$-complete field extension every $\underline{G(\bbQ_p)}$-torsor for the v-topology over $\Spd F$ comes via descent datum for the quasi-pro-\'etale cover $\Spd\bbC_p\to \Spd F$.
Standard arguments relating Galois theory to the theory of descent will show that the groupoid of $\underline{G(\bbQ_p)}$-torsors over $\Spd F$ is equivalent to the groupoid of continuous group homomorphisms $\rho:\Gamma_F\to G(\bbQ_p)$. 
Furthermore, given a continuous group map $\rho:\Gamma_F\to G(\bbQ_p)$ with associated $\underline{G(\bbQ_p)}$-torsor $\calF_\rho$ we have a $G(\bbQ_p)$-equivariant identification $|\calF_\rho|\simeq G(\bbQ_p)/\rho(\Gamma_F)$. In particular, for any $x\in |\calF_\rho|$ the stabilizer of $x$ in $G(\bbQ_p)$, $G_x\subseteq G(\bbQ_p)$, is conjugate to the image of $\rho$. 
\end{remark}

Recall the category of crystalline Galois representations \cite[$\mathsection$ 5.1]{MR657238}. 

\begin{proposition}
	\label{crystallinerepsongrassmanian}
	Fix a field extension $[F:\breve{E}]<\infty$ and a crystalline representation with $G$-structure 
	$\rho:\Gamma_F\to G(\bbQ_p)$. 
	Let $D_{\on{cris}}(\rho)=(D_\rho,\Phi_\rho,\on{Fil}_{F,\rho}^\bullet)$ denote the filtered isocrystal with $G$-structure associated to $\rho$.
	Assume that the isocrystal with $G$-structure $(D_\rho,\Phi_\rho)$ is isomorphic to $(V_b,\Phi_b)$, and that the filtration type of $\on{Fil}_{F,\rho}$ is $\mu$.
	Then there exists a map $f_\rho:\Spd F\to \Gr^{\leq \mu,b-\on{adm}}_{G,\Div^1_E}$ such that 
	\[\pi_{\on{GM}}:\Sht_{G,b,\mu,\infty}\times_{\Gr^{\leq \mu,b-\on{adm}}_{G,\Div^1_E}} \Spd F\to \Spd F\] 
	is isomorphic to the $\underline{G(\bbQ_p)}$-torsor corresponding to $\rho$ under the equivalence of \Cref{torsors-vs-galoisreps}.
\end{proposition}
\begin{proof}
	This boils down to \cite[Proposition 10.5.3, 10.5.6, Lemme 10.5.4]{FF18}, let us elaborate.	
	The reference attaches to any filtered isocrystal $(D,\Phi,\on{Fil}_F^\bullet)$ a modification of vector bundles over $X_{\on{FF},\bbC_p^\flat}$ along $\infty$, 
	\[\calE(D,\Phi,\on{Fil}_F^\bullet)\dashrightarrow\calE(D,\Phi).\]
	The vector bundle $\calE(D,\Phi,\on{Fil}_F^\bullet)$ comes equipped with $\Gamma_F$-action. Moreover, as \cite[Proposition 10.5.6]{FF18} explains, whenever $(D,\Phi,\on{Fil}_F^\bullet)$ is a weakly admissible filtered isocrystal the vector bundle $\calE(D,\Phi,\on{Fil}_F^\bullet)$ is semi-stable of slope $0$ and the rule
	\[(D,\Phi,\on{Fil}_F^\bullet)\mapsto \on{H}^0(X_{\on{FF},\bbC^\flat_p},\calE(D,\Phi,\on{Fil}_F^\bullet))\]
together with its $\Gamma_F$-action agrees with Fontaine's functor $V_{\on{cris}}$.
Finally, as \cite[Lemme 10.5.4]{FF18} explains, the type of the modification agrees with the type of the filtration $\on{Fil}_F^\bullet$.
Passage to objects with $G$-structure is formal since the relevant construction are $\otimes$-exact.

Starting with $\rho:\Gamma_F\to G(\bbQ_p)$, we consider $(D_\rho,\Phi_\rho,\on{Fil}^\bullet_{F,\rho}):=D_{\on{cris}}(\rho)$ as a filtered isocrystal with $G$-structure. 
After fixing an isomorphism $(D_\rho,\Phi_\rho)\simeq (V_b,\Phi_b)$ of isocrystals with $G$-structure, we may transfer the weakly admissible filtration $\on{Fil}^\bullet_{F,\rho}$ on $D_\rho\otimes_{\breve{\bbQ}_p} F$ to a weakly admissible filtration $\on{Fil}^\bullet_{F,\rho,b}$ on $V_b\otimes_{\breve{\bbQ}_p} F$. 
By the above, we get a modification 
\[\calE(V_b,\Phi_b,\on{Fil}^\bullet_{F,\rho,b})\dashrightarrow \calE_b,\]
of $G$-bundles over $X_{\on{FF},\bbC_p}$ where $\calE(V_b,\Phi_b,\on{Fil}^\bullet_{F,\rho,b})$ is isomorphic to the trivial $G$-torsor over $X_{\on{FF},\bbC_p^\flat}$.
This is precisely a $\Spd \bbC_p$-point of $\Gr^{\leq \mu,b-\on{adm}}_{G,\Div^1_E}$.
The $\Gamma_F$-action on $\calE(V_b,\Phi_b,\on{Fil}^\bullet_{F,\rho,b})$ descends this to a $\Spd F$-point.  
The $\underline{G(\bbQ_p)}$-torsor over $\Spd F$ is the moduli space of trivializations of the form 
\[\tau:\calE_{\mathds{1}}\xrightarrow{\simeq} \calE(V_b,\Phi_b,\on{Fil}^\bullet_{F,\rho,b}).\]
Since $\on{H}^0(X_{\on{FF},\bbC^\flat_p},\calE(D,\Phi,\on{Fil}_{F,\rho,b}^\bullet))\simeq V_{\on{cris}}\circ D_{\on{cris}}(\rho)\simeq \rho$ the $\underline{G(\bbQ_p)}$-torsor over $\Spd F$ is isomorphic to the one associated to $\rho$. 
\end{proof}

\subsection{Geometric connected components in the $G^\der=G^\scn$ case.}
\label{thesimplestcase}

Recall that a connected reductive group over $\bbQ_p$ is said to be an \emph{unramified} group if it is quasi-split and it splits over an unramified extension of $\bbQ_p$. 
Equivalently, $G$ is unramified if and only if $G$ admits a reductive integral model over $\bbZ_p$ \cite[Proposition 4.2, 4.3]{MR4060875}.

For the rest of this section, we let our fixed group $G$ be an unramified group over $\Spec \bbQ_p$.
We fix data $T\subseteq B\subseteq G$ with $B\subseteq G$ a Borel subgroup and $T\subseteq B$ a maximally split maximal torus defined over $\bbQ_p$. 
We let $A\subseteq T$ denote the maximal split subtorus.

\subsubsection{The set $B(G,\mu)$ and HN-irreducibility}
\label{HN-irreducible}
Recall that $\Rep_G$ denotes the Tannakian category of algebraic representations of $G$ on finite $\bbQ_p$ vector spaces.
Recall that the category of isocrystals $\Isoc_{\bar{\bbF}_p}$ is a semisimple Tannakian category and that one can naturally endow it with a $\bbQ$-grading given by the slope decomposition \cite[$\mathsection$ 2.1, 3.2]{Kot97}, \cite[\S 3.3, Th\'eor\`eme 3.3.2]{saavedra_rivano_categories_tannakiennes}.
Consequently, any isocrystal with $G$-structure
\[V_b:\Rep_G\to \Isoc_{\bar{\bbF}_p}\]
gives rise to a slope morphism
\[\nu_b:\bbD\to G_{\bar{\bbQ}_p}\]
where $\bbD$ denotes the pro-torus with character group $\bbQ$. 
We let $\nu^\dom_b$ denote the unique dominant cocharacter factoring through $T$ and conjugate to $\nu_b$. 
Then $\nu^\dom_b$ is defined over $\bbQ_p$ and factors through $A\subseteq T$ (see \cite[$\mathsection$ 4]{MR2484281} \cite[(1.1.3.1)]{MR761308}.
This construction defines a map 
\[\nu^\dom_{(-)}:B(G)\to X_*^+(A)_\bbQ\]
\[b\mapsto \nu_b^\dom\]
which is usually referred to as the Newton map.

Recall Borovoi's algebraic fundamental group $\pi_1(G)$ \cite{borovoi1989algebraic} \cite[$\mathsection$ 1.13]{RR96} which can be defined as the quotient of $X_*(T)$ by the co-root lattice. 
This group comes equipped with $\Gamma_{\bbQ_p}$ action and Kottwitz constructs a map 
\[\kappa_G:B(G)\to \pi_1(G)_{\Gamma_{\bbQ_p}}\]
that is usually referred to as the Kottwitz map \cite{Kot85} \cite[Theorem 1.15]{RR96}. 

An important result of Kottwitz \cite[$\mathsection$ 4.13]{Kot97} states that the map of sets 
\[(\nu^\dom_{(-)},\kappa_G):B(G)\to X_*^+(A)_\bbQ \times \pi_1(G)_{\Gamma_{\bbQ_p}}\]
is injective. 
Now, if we are given a conjugacy class of cocharacters $\mu\in\{\bbG_m\to G\}/\sim$, its reflex field $E$ is necessarily an unramified extension of $\bbQ_p$ since we are assuming that $G$ splits over an unramified extension. 
There is a unique dominant cocharacter $\mu_0\in X^+_*(T)$ conjugate to $\mu$.
Moreover, since $B$ and $T$ are defined over $\bbQ_p$, $\Gamma_{\bbQ_p}$ acts on $X^+_*(T)$ and we can construct an element 
\[\bar{\mu}\in X_*^+(A)_\bbQ= X_*^+(T)_\bbQ^{\Gamma_{\bbQ_p}}\]
by averaging over the Galois orbit $\on{Gal}(E/\bbQ_p)\cdot \mu_0$.
More precisely,

\[\bar{\mu}=\frac{1}{[E:\bbQ_p]} \sum_{\gamma\in {\on{Gal}}(E/\bbQ_p)} \mu_0^\gamma.\]
We can now recall Kottwitz' definition of the set $B(G,\mu)\subseteq B(G)$ for unramified groups $G$ \cite[$\mathsection$ 6]{Kot97}.

\begin{definition}
	\label{defi:BGmu}
	Fix notation $A\subseteq T\subseteq B\subseteq G$ as above.
	\begin{enumerate}
		\item The set $B(G,\mu)$ consists of those conjugacy classes $b\in B(G)$ for which $\kappa_G(b)=\mu$ in $\pi_1(G)_{\Gamma_{\bbQ_p}}$ and for which $\bar{\mu}-\nu_b^{\dom}\in X_*^+(A)_\bbQ$ is a non-negative $\bbQ$-linear combination of positive co-roots.
		\item We say that $b\in B(G,\mu)$ is \emph{Hodge--Newton irreducible} with respect to $\mu$ if all the coefficients $\bar{\mu}-\nu_b^{\dom}\in X_*^+(A)_\bbQ$ as $\bbQ$-linear combination of simple coroots are strictly positive.
		\item We say that the $p$-adic shtuka datum $(G,b,\mu)$ is HN-irreducible if $b\in B(G,\mu)$ and $b$ is HN-irreducible with respect to $\mu$. If the group $G$ is clear from the context we also say that the pair $(b,\mu)$ is HN-irreducible when the triple $(G,b,\mu)$ is HN-irreducible.
	\end{enumerate}
\end{definition}

\subsubsection{The $G(\bbQ_p)$ action in the $G^\der=G^\scn$ case}
Recall that $G$ is an unramified reductive group over $\bbQ_p$. 
From now on, we require that $\calG$ is a hyperspecial parahoric model of $G$ over $\bbZ_p$ i.e. $\calG$ is a smooth affine group over $\bbZ_p$ whose generic fiber is $G$ and whose special fiber is also connected reductive.
We fix a $p$-adic shtuka datum $(G,b,\mu)$ and assume that it is HN-irreducible as in \Cref{defi:BGmu}.
We let $G^\der$, $G^\scn$ and $G^\ab$ denote respectively the derived subgroup of $G$, the simply connected cover of $G^\der$ and the maximal abelian quotient of $G$.
All of these algebraic groups also have hyperspecial parahoric models that we denote by $\calG^\der$, $\calG^\scn$ and $\calG^\ab$.
We let 
\[\det:\calG\to \calG^\ab \text{ and } \det:G\to G^\ab\]
denote the natural quotient maps. 
Following \cite{Che13} we call them \textit{determinant maps}.

Let $(G^\ab,b^\ab,\mu^\ab)$ denote the shtuka datum with $b^\ab=\det(b)$, $\mu^\ab=\det\circ \mu$.
By functoriality, \Cref{functoriality-wtih-group}, we obtain morphism of spaces
\[\det:\Sht^{\on{geo}}_{G,b,\mu,\infty}\to \Sht^{\on{geo}}_{G^\ab,b^\ab,\mu^\ab, \infty} \text{ and } \det: \Shtm{{O_{\breve{E}}}}\to \Sht^{\mu^\ab}_{\calG^\ab,O_{\breve{E}}}(b^\ab).\]

Until further notice, we assume that the derived subgroup of $G$ is simply connected i.e. that $G^\der=G^\scn$.
We will need the following input from the study of connected components of affine Deligne--Lusztig varieties \cite{CKV15}, \cite{HZ20}, \cite{Nie18} for unramified groups.   
\begin{proposition}
	\label{cor:adlvconnected}
	Suppose that $G^\der=G^\scn$ and that $(b,\mu)$ is HN-irreducible. Then the natural map 
	\[\on{det}:\pi_0(X^{\leq \mu}_\calG(b))\to \pi_0(X^{\mu^\ab}_{\calG^\ab}(b^\ab))\] 
	is bijective.
\end{proposition}
\begin{proof}
	One can verify that a pair $(b,\mu)$ is HN-irreducible if and only if for every $\bbQ_p$-simple factor $G_i$ of $G^\ad$ with projection map $\pi_i:G\to G_i$ the pair $(b_i,\mu_i):=(\pi_i(b), \pi_i \circ \mu)$ is HN-irreducible. 
	Indeed, the coefficient of $\mu^{\on{dom}}-\nu_b^{\on{dom}}$ associated to a positive root can be computed on the $\bbQ_p$-simple factors of the adjoint quotient (see proof of \cite[Corollary 4.1.16]{CKV15}). 
	Moreover, it also follows directly from the definitions that the formation of affine Deligne--Lusztig varieties open products of groups i.e. if $\calG=\calG_1\times \calG_2$, $b=(b_1,b_2)$ and $\mu=(\mu_1,\mu_2)$ then, 
	\[{X^{\leq \mu}_{\calG}(b)}\simeq {X^{\leq \mu_1}_{\calG_1}(b_1)}\times {X^{\leq \mu_2}_{\calG_2}(b_2)}.\]
	From \cite[Corollary 2.4.3]{CKV15} and with the notation $\omega_{G_i}$ used as in loc. cit. we get a commutative diagram with Cartesian squares
	\begin{center}
		\begin{tikzcd}
			\pi_0(\adlvdos) \ar{r}\ar{d}{w_G}& \pi_0({X^{\leq \mu^{\on{ad}}}_{\calG^{\on{ad}}}(b^{\on{ad}})}) \ar{r}{\simeq}\ar{d}{\omega_{G^{\on{ad}}}} &  \pi_0(X^{\leq \mu_1}_{\calG_1}(b_1))\times \dots \times \pi_0(X^{\leq \mu_n}_{\calG_1}(b_n))\ar{d}{w_{G_i}}\\
			c_{b,\mu} \pi_1(G)^{\Gamma_{{\bbQ_p}}} \ar{r} & c_{b^{\on{ad}},\mu^{\on{ad}}} \pi_1(G^{\on{ad}})^{\Gamma_{{\bbQ_p}}}\ar{r}{\simeq}  &  c_{b_{1},\mu_{1}}\pi_1(G_{1})^{\Gamma_{{\bbQ_p}}}\times \dots \times c_{b_{n},\mu_{n}}\pi_1(G_{n})^{\Gamma_{{\bbQ_p}}}.
		\end{tikzcd}
	\end{center}
	The vertical right hand map is a bijection by \cite[Theorem 1.1]{CKV15}, \cite[Theorem 8.1]{HZ20}, \cite[Theorem 1.1]{Nie18} which implies the vertical left hand map is also a bijection. 
	The result follows from showing that in the commutative diagram below the bottom horizontal arrow and the vertical right hand arrow are both bijective. 
	\begin{center}
		\begin{tikzcd}
			\pi_0(\adlvdos) \ar{r} \ar{d}{w_G} & \pi_0(X^{\mu^\ab}_{\calG^\ab}(b^\ab)) \ar{d}{w_{G^\ab}}\\
			c_{b,\mu}\pi_1(G)^{\Gamma_{\bbQ_p}} \ar{r} & c_{b_{\ab},\mu_{\ab}}\pi_1(G^{\ab})^{\Gamma_{\bbQ_p}} 
		\end{tikzcd}
	\end{center}
	Since $G^\der$ is simply connected we have a $\Gamma_{\bbQ_p}$-equivariant identification $\pi_1(G)\to \pi_1(G^{\on{ab}})$ so the bottom map is easily seen to be a bijection. 
	Moreover, the adjoint quotient of $G^\ab$ is $\{e\}$ and \cite[Corollary 2.4.3]{CKV15} says that 
\[w_{G^{\on{ab}}}:\pi_0(X^{\mu^\ab}_{\calG^\ab}(b^\ab))\to c_{b_\ab,\mu_\ab}\pi_1(G^\ab)^{\Gamma_{\bbQ_p}}\] is an isomorphism in this case.
\end{proof}

\begin{remark}
	Since $X^{\mu^\ab}_{\calG^\ab}(b^\ab)\subseteq \calF\ell_{\calG^\ab,\bbW}$ is $0$-dimensional and a disjoint union of copies of $\Spec \bar{\bbF}_p$, we can reformulate \Cref{cor:adlvconnected} as saying that the map $\adlvdos\to X^{\mu^\ab}_{\calG^\ab}(b^\ab)$ has geometrically connected fibers.
\end{remark}

\begin{lemma}
	\label{remarktori}
The space $\Sht^{\on{geo}}_{G^\ab,b^\ab,\mu^\ab,\infty}$ is a trivial $\underline{G^\ab(\bbQ_p)}$-torsor over $\Spd \bbC_p$.
\end{lemma}
\begin{proof}
	We clam that when $G$ is abelian (a torus) the structure map produces an identification \cite[$\mathsection$ 25.1]{SW20} 
\[\Gr^{\leq \mu,b-\on{adm}}_{G,\Div^1_E}\simeq \Div^1_E.\]
Indeed, the map $f:\Gr^{\leq \mu}_{G,\Div^1_E}\to \Div^1_E$ is proper and one can verify from the definitions and the Cartan decomposition that it is bijective on geometric points.  
By \cite[Lemma 12.5]{Sch17}, $f$ is an isomorphism (see also \cite[Proposition 21.3.1]{SW20}). 

Moreover, $\Gr^{\leq \mu,b-\on{adm}}_{G,\Div^1_E}\subseteq \Gr^{\leq \mu}_{G,\Div^1_E}$ is a non-empty open subset. 
Since $|\Gr^{\leq \mu}_{G,\Div^1_E}| \simeq |\Div^1_E|$ is precisely one point, we must also have $\Gr^{\leq \mu,b-\on{adm}}_{G,\Div^1_E}=\Gr^{\leq \mu}_{G,\Div^1_E}$. 

Now, $\pi_{\on{GM}}$ makes $\Sht^{\on{geo}}_{G,b,\mu,\infty}$ into a $\underline{G(\bbQ_p)}$-torsor over $\Gr^{\leq \mu,b-\on{adm}}_{G,\Div^1_E}\times_{\Div_E^1}\Spd \bbC_p=\Spd \bbC_p$.
It suffices to show that every $\underline{G(\bbQ_p)}$-torsor is trivial over $\Spd \bbC_p$. 
But any pro-\'etale cover of $\Spd \bbC_p$ splits \cite[Corollary 7.22]{Sch17}.
\end{proof}

We set some notation. 
Let $\calK=\calG(\bbZ_p)$, let $\calK^\der=\calG^\der(\bbZ_p)$ and let $\calK^\ab=\calG(\bbZ_p)$. 
An application of Lang's theorem shows that $\calK^\ab=\det(\calK)$.
For any $m\in \Sht^{\on{geo}}_{G^\ab,b^\ab,\mu^\ab,\infty}(\Spd \bbC_p)$ we let $X_m$ denote the geometric fiber defined by the following Cartesian diagram
\begin{center}
\begin{tikzcd}
	X_m \arrow{r} \arrow{d}  & \Spd \bbC_p \arrow{d}{m} \\
 \Sht^{\on{geo}}_{G,b,\mu,\infty} \arrow{r} & \Sht^{\on{geo}}_{G^\ab,b^\ab,\mu^\ab,\infty}. 
\end{tikzcd}
\end{center}

\begin{lemma}
	\label{transitiveofgder}
	Let $(G,b,\mu)$ and $X_m$ be as above. Then the following statements hold.
	\begin{enumerate}
		\item $G(\bbQ_p)$ acts transitively on $\pi_0(\Sht^{\on{geo}}_{G,b,\mu,\infty})$.
	\label{transitiveofGder}
		\item $\calK^\der$ acts transitively on $\pi_0(X_m)$ for all $m\in \Sht^{\on{geo}}_{G^\ab,b^\ab,\mu^\ab,\infty}(\Spd \bbC_p)$.
	\label{transitiveofKder}
	\end{enumerate}
\end{lemma}
\begin{proof}
	Let us show the first claim.
	Observe that since $\Sht^{\on{geo}}_{G,b,\mu,\infty}\to \Sht^{\on{geo}}_{G,b,\mu,\calK}$ is a $\underline{\calK}$-torsor,
	the group $\calK$ acts transitively on the fibers of $f:\pi_0(\Sht^{\on{geo}}_{G,b,\mu,\infty})\to \pi_0(\Sht^{\on{geo}}_{G,b,\mu,\calK})$.
	Let $x_\infty\in \pi_0(\Sht^{\on{geo}}_{G,b,\mu,\infty})$ and let $x_\calK\in \pi_0(\Sht^{\on{geo}}_{G,b,\mu,\calK})$. 
	It suffices to find $g\in G(\bbQ_p)$ such that $f(x_\infty \star g)=x_\calK$. 
	By functoriality of the specialization map, \Cref{functoriality-of-special}, we have the following commutative diagram 
	\begin{equation}
	\begin{tikzcd}
		\label{commutativediagramconnectedop}
		\pi_0(\Sht^{\on{geo}}_{G,b,\mu,\infty}) \arrow{r}{\on{det}} \arrow{d}{f}  & \pi_0(\Sht^{\on{geo}}_{G^\ab,b^\ab,\mu^\ab,\infty}) \arrow{d}{f^\ab} \\
			\pi_0(\Sht^{\on{geo}}_{G,b,\mu,\calK}) \arrow{r}{\on{det}} \arrow{d}{\pi_0(\on{sp})}  & \pi_0(\Sht^{\on{geo}}_{G^\ab,b^\ab,\mu^\ab,\calK^\ab}) \arrow{d}{\pi_0(\on{sp})} \\
		\pi_0(\adlvdos) \arrow{r} & \pi_0(X^{\mu^\ab}_{\calG^\ab}(b^\ab).
	\end{tikzcd}
	\end{equation}
	By \Cref{thm:specializtheorem} and \Cref{cor:adlvconnected} all the maps in the lower square are bijective, so by using \Cref{functoriality-wtih-group} it suffices to find $g\in G(\bbQ_p)$ with 
	\[f^\ab(\det(x_\infty))\star \det(g)=\det\circ f(x_\infty\star g)=\det(x_\calK).\] 
	Recall that $G^\ab(\bbQ_p)$ acts transitively on $\pi_0(\Sht^{\on{geo}}_{G^\ab,b^\ab,\mu^\ab,\infty})$ by \Cref{remarktori}. 
	To find $g$ as above it suffices to show that the map $G(\bbQ_p)\to G^\ab(\bbQ_p)$ is surjective.
	But this follows from Kneser's theorem \cite[Satz 1]{Kne65}. 

	Let us show the second statement, take $x_\infty,y_\infty\in \pi_0(X_m)$ and denote by $x_\calK$ and $y_\calK$ their images in $\pi_0(\Sht^{\on{geo}}_{G,b,\mu,\calK})$.
	By hypothesis $\det(x_\infty)=\det(y_\infty)=m$, so it follows from the bijectivity of the lower square in diagram \ref{commutativediagramconnectedop} that $x_\calK=y_\calK$.
	In particular, there is $g\in \calK$ with $x_\infty\star g=y_\infty$.
	But then $\det(x_\infty)=\det(x_\infty)\star \det(g)$.
	By \Cref{remarktori}, the action of $G^\ab(\bbQ_p)$ is simple.
	This implies that $g\in G^\der(\bbQ_p)\cap \calK=\calK^\der$ as we wanted to show.
\end{proof}

\begin{theorem}
\label{thm:maintheorem}
Suppose that $G$ is an unramified group over $\bbQ_p$, that $G^{\on{der}}=G^\scn$ and that $(G,b,\mu)$ is HN-irreducible. 
Then the determinant map 
\[\det:\Sht_{G,b,\mu,\infty} \to \Sht_{G^\ab,b^\ab,\mu^\ab,\infty}\] has connected geometric fibers.
\end{theorem}
	\begin{proof}
		By \Cref{remarktori}, $\Sht_{G^\ab,b^\ab,\mu^\ab,\infty}\times_{\Div^1_{E(\mu^\ab)}} \Spd \bbC_p$ is isomorphic to $\underline{G^\ab(\bbQ_p)}\times \Spd{\bbC_p}$.
		Consequently, we may prove instead that the determinant map induces a bijection 
		\[\pi_0(\on{det}): \pi_0(\Sht^{\on{geo}}_{G,b,\mu,\infty}) \to \pi_0(\Sht^{\on{geo}}_{G^\ab,b^\ab,\mu^\ab,\infty}).\]
		Let $x\in \pi_0(\Sht^{\on{geo}}_{G,b,\mu,\infty})$. 
		Given $F$ a finite extension of $\breve{E}$ we let $x_F$ denote the image of $x$ on $\pi_0(\Sht^F_{G,b,\mu,\infty})$. 

		Recall Chen's result on generic crystalline representations for HN-irreducible datum $(G,b,\mu)$ \cite[Th\'eor\`eme 5.0.6]{Chen}. 
		A consequence of this theorem is the existence of a finite field extension $[K:\breve{E}]<\infty$ and a crystalline representation with $G$-structure $\xi_{b,\mu}:\Gamma_K\to G(\bbQ_p)$ whose image contains an open subset of $G^{\der}(\bbQ_p)$.
		Fix such a representation and let $f:=f_{\xi_{(b,\mu)}}$ denote one of the maps associated to $\xi_{b,\mu}$ under \Cref{crystallinerepsongrassmanian} of the form 
		\[f:\Spd{K}\to \Gr^{\leq \mu,b-\on{adm}}_{G,\Div^1_E}.\] 
		Let \[S_f:=\Sht_{G,b,\mu,\infty}\times_{\Gr^{\leq \mu,b-\on{adm}}_{G,\Div^1_E}} \Spd K\] 
		denote the fiber $\pi^{-1}_{\on{GM}}(f)$.
Let $s\in \pi_0(S_f)$ be an element mapping to $x_K$. 
In summary, we have taken a commutative diagram of sets,  
		\begin{center}
		\begin{tikzcd}
			* \ar{r}{x} \ar{d}{s}\ar[swap]{rd}{x_K} & \pi_0(\Sht^{\on{geo}}_{G,b,\mu,\infty})\ar{d} \\
			\pi_0(S_f) \ar{r}{f}  & \pi_0(\Sht^{K}_{G,b,\mu,\infty}).
		\end{tikzcd}
		\end{center}

		We let $G^\der_x$ (respectively $G^{\on{der}}_{x_K}$ and $G^{\on{der}}_s$) denote the stabilizer of $x$ (respectively of $x_K$ and of $s$) in $G^{\on{der}}(\bbQ_p)$ of its action on $\pi_0(\Sht^{\on{geo}}_{G,b,\mu,\infty})$ (respectively $\pi_0(\Sht^K_{G,b,\mu,\infty})$ and $\pi_0(S_f)$).
		By $G^\der(\bbQ_p)$-equivariance, we have inclusions $G^{\on{der}}_x,G^{\on{der}}_s\subseteq G^{\on{der}}_{x_K}$.
Moreover, by \Cref{torsors-vs-galoisreps}, $G^\der_s$ is conjugate to $\xi_{b,\mu}(\Gamma_K)\cap G^\der(\bbQ_p)$ and by Chen's theorem \cite[Th\'eor\`eme 5.0.6]{Chen} $G^{\on{der}}_s$ is an open subgroup of $G^{\on{der}}(\bbQ_p)$. 
This also implies that $G^{\on{der}}_{x_K}$ is an open subgroup of $G^{\on{der}}(\bbQ_p)$. 

By \Cref{transitiveofgder}.(\ref{transitiveofKder}), $G^{\on{der}}_x\cdot \calK^{\on{der}}=G^{\on{der}}(\bbQ_p)$ which implies that $G^{\on{der}}_{x_K}\cdot \calK^{\on{der}}=G^{\on{der}}(\bbQ_p)$ as well. 
In particular, the projection map $\calK^{\on{der}}\to G^{\on{der}}(\bbQ_p)/G^{\on{der}}_{x_K}$ is surjective. 
Since $G^{\on{der}}(\bbQ_p)/G^{\on{der}}_{x_K}$ has the discrete topology and $\calK^{\on{der}}$ is compact, we get that $G^{\on{der}}_{x_K}$ is closed and of finite index within $G^{\on{der}}(\bbQ_p)$. 
In particular, it is closed and of finite covolume.
We may apply Margulis's theorem \cite[Chapter II, Theorem 5.1]{Marg} to conclude that $G^{\on{der}}_{x_K}=G^{\on{der}}({\bbQ_p})$. 
Indeed, we just verified $G^{\on{der}}_{x_K}$ is closed and of finite covolume, and since $G^{\on{der}}$ is quasi-split (even unramified) all of the simple factors of $G^{\on{der}}$ are isotropic. 
Moreover, since $G^\der=G^\scn$ all of the simple factors are simply connected.

Since the argument doesn't depend on the choice of $x$ the action of $G^{\on{der}}(\bbQ_p)$ on $\pi_0(\Sht^K_{G,b,\mu,\infty})$ is trivial.
Now, $\Spd{\bbC_p}=\varprojlim_{\breve{E}\subseteq K\subseteq \bbC_p} \Spd{K}$ and we may use \cite[Lemma 11.22]{Sch17} to compute the action map
\[|\Sht^{\on{geo}}_{G,b,\mu,\infty}|\times {G^{\on{der}}({\bbQ_p})} \to |\Sht^{\on{geo}}_{G,b,\mu,\infty}|\]
		as the limit of the action maps
\[\varprojlim_{\breve{E}\subseteq K\subseteq \bbC_p}|\Sht^{K}_{G,b,\mu,\infty}|\times {G^{\on{der}}({\bbQ_p})} \to |\Sht^{K}_{G,b,\mu,\infty}|.\]
Since in the transition maps $|\Sht^{K_1}_{G,b,\mu,\infty}|\to |\Sht^{K_2}_{G,b,\mu,\infty}|$ every connected component on the source surjects onto a connected component on the target, we get that 
\[\pi_0(\Sht^{\on{geo}}_{G,b,\mu,\infty})=\varprojlim_{\breve{E}\subseteq K\subseteq \bbC_p} \pi_0(\Sht^{K}_{G,b,\mu,\infty}).\] 
This proves that $G^{\on{der}}(\bbQ_p)$ acts trivially on $\pi_0(\Sht^{\on{geo}}_{G,b,\mu,\infty})$. 
Consequently, the $G(\bbQ_p)$-action on $\pi_0(\Sht^{\on{geo}}_{G,b,\mu,\infty})$ factors through $G^{\on{ab}}(\bbQ_p)$. 
Since the action is transitive by \Cref{transitiveofgder}, $\pi_0(\on{det})$ must be bijective. 
	\end{proof}

\begin{corollary}
	\label{corollaryRVconj}
	\Cref{conj-rapoport-viehmann} holds whenever $G$ is unramified.
\end{corollary}
\begin{proof}
	It suffices to show the formula $\pi_0(\bbM^{\on{geo}}_K)\simeq \pi_0(\bbM^{\ab,\on{geo}}_{\on{det}(K)})\simeq G^\ab(\bbQ_p)/\det(K)$ since representability and compatibility with group actions was already settled in \cite{SW20}.
	When $G=G^\der$ it follows from \Cref{thm:maintheorem} that $\pi_0(\Sht^{\on{geo}}_{G,b,\mu,\infty})\simeq \pi_0(\Sht^{\on{geo}}_{G^\ab,b^\ab,\mu^\ab,\infty})$ as $G^\ab(\bbQ_p)$-torsors. 
	We also have that 
	\[\pi_0(\Sht^{\on{geo}}_{G,b,\mu,K})= \pi_0(\Sht^{\on{geo}}_{G,b,\mu,\infty}/\underline{K})\simeq \pi_0(\Sht^{\on{geo}}_{G,b,\mu,\infty})/{K}\simeq \pi_0(\Sht^{\on{geo}}_{G,b,\mu,\infty})/\on{det}{(K)}\]
	and 
	\[\pi_0(\Sht^{\on{geo}}_{G^\ab,b^\ab,\mu^\ab,\on{det}(K)})= \pi_0(\Sht^{\on{geo}}_{G^\ab,b^\ab,\mu^\ab,\infty}/\underline{\on{det}(K)})\simeq \pi_0(\Sht^{\on{geo}}_{G^\ab,b^\ab,\mu^\ab,\infty})/\on{det}({K}).\]
	This gives 
\[\pi_0(\bbM^{\on{ab},\on{geo}}_{\on{det}(K)})\simeq \pi_0(\Sht^{\on{geo}}_{G^\ab,b^\ab,\mu^\ab,\on{det}(K)})\simeq \pi_0(\Sht^{\on{geo}}_{G,b,\mu,K})\simeq \pi_0(\bbM^{\on{geo}}_K)\]
as we wanted to show.
\end{proof}

\subsubsection{The $G_b(\bbQ_p)\times W_E$-action in the $G^\der=G^\scn$ case.}
\label{subsubsectionactions}
On this subsection we keep the notation of the previous one. Namely, $G$ is an unramified reductive group over $\bbQ_p$, $(G,b,\mu)$ is HN-irreducible, and $G^\der=G^\scn$. 
We know from \Cref{thm:maintheorem} that the action of $G^\ab(\bbQ_p)$ on $\pi_0(\Sht^{\on{geo}}_{G,b,\mu,\infty})$ is simple and transitive. 
In particular, for any element $j\in G_b(\bbQ_p)$ (respectively $\gamma\in W_E$) and an element $x\in \pi_0(\Sht^{\on{geo}}_{G,b,\mu,\infty})$ there is a unique element $g_j\in G^\ab(\bbQ_p)$ (respectively $g_\gamma\in G^\ab(\bbQ_p)$) such that 
\[x\star_{G_b} j=x\star g_j \text{ (respectively } x\star_{W_E} \gamma=x\star g_\gamma\text{).}\]
The rules $j\mapsto g_j$ and $\gamma\mapsto g_\gamma$ define group homomorphisms 
\[\rho_{G_b}:G_b(\bbQ_p)\to G^\ab(\bbQ_p) \text{ and } \rho_{W_E}:W_E\to G^\ab(\bbQ_p)\]
that do not depend of the choice of $x$.
The purpose of this subsection is to make $\rho_{G_b}$ and $\rho_{W_E}$ as explicit as possible.
The basic principle that allows us to compute $\rho_{G_b}$ and $\rho_{W_E}$ is that the rule 
\[(G,b,\mu)\mapsto \pi_0(\Sht^{\on{geo}}_{G,b,\mu,\infty})\]
is functorial \Cref{functoriality-wtih-group}.

We first study $\rho_{G_b}$.
Functoriality gives a map 
$G_b(\bbQ_p)\to G^\ab_{b^\ab}(\bbQ_p)$. 
Moreover, if we regard $G^\ab(\bbQ_p)$ and $G^\ab_{b^\ab}(\bbQ_p)$ as the subgroups of elements of $g\in G^\ab(\breve{\bbQ}_p)$ with the property $g^{-1}\mathds{1}\phi(g)=\mathds{1}$ and $g^{-1}b\phi(g)=b$, then commutativity of $G^\ab$ readily implies $G^\ab_{b^\ab}(\bbQ_p)=G^\ab(\bbQ_p)$.
Overall, this gives rise to a map 
\[\on{det}_b:G_b(\bbQ_p)\to G^\ab(\bbQ_p).\]

\begin{proposition}
	\label{theactionofJ}
	With the notation as above, for all $j\in G_b(\bbQ_p)$ we have that $\rho_{G_b}(j)=({\on{det}_b}(j))^{-1}$.	
\end{proposition}
\begin{proof}
	Functoriality and \Cref{thm:maintheorem} reduces us to show that for all $g\in G^\ab(\bbQ_p)$ and $x\in \pi_0(\Sht^{\on{geo}}_{G^\ab,b^\ab,\mu^\ab,\infty})$ the following identity holds  
	\[x\star_{G^\ab} g=x\star_{G^\ab_{b^\ab}}{g^{-1}}.\]
		In turn, by our definition of the right action of $G^\ab_{b^\ab}$ ($\mathsection$ \ref{definingsection}) this is the same as showing
		\[g\star_{G_{b^\ab}^\ab} x = x\star_{G^\ab} g\]
		for the natural left action.
		Since $G^\ab$ is abelian any modification $\alpha:\calE_{\mathds{1}}\dashrightarrow \calE_b$ induces an isomorphism 
		\[G^\ab(\bbQ_p)=\on{Aut}(\calE_{\mathds{1}})\to \on{Aut}(\calE_{b^{\ab}})=G^\ab_{b^\ab}(\bbQ_p) \text{ with } g\mapsto \alpha \circ g\circ\alpha^{-1}\]
		and this isomorphism does not depend on the choice of $\alpha$.
		For such an $\alpha$, if follows immediately that 
		\[x\star_{G^\ab} g=[\alpha\circ g\circ \alpha^{-1}] \star_{\on{Aut}(\calE_{b^{\ab}})} x\]
		By considering, $\calE_\mathds{1}$ and $\calE_b$ as trivial $\varphi$-modules with $G^\ab$-structure over $\calY^{\bbC_p^\flat}_{(0,\infty)}$ one can show that the identification between $\on{Aut}(\calE_{\mathds{1}})$ and $\on{Aut}(\calE_b)$ produced by any $\alpha$ as above agrees with the standard one (as subgroups of $G^\ab(\breve{\bbQ}_p)\subseteq G^\ab(B^{\bbC_p^\flat}_{(0,\infty)})$).  
\end{proof}

We now compute $\rho_{W_E}$.
We learned the following line of reasoning from \cite[Lemma 1.22]{RZ96}, which in turn is an elaboration of an argument in \cite[$\mathsection$ 12]{Kot92}. 
	Let $E\subseteq \bar{\bbQ}_p$ denote a finite field extension of $\bbQ_p$ which is a subfield of $\bar{\bbQ}_p$. 
	Let $\{\on{Tori}_{\bbQ_p}\}$ denote the category of tori defined over $\bbQ_p$. 
	Recall the functor $X_*(-):\{\on{Tori}_{\bbQ_p}\}\to \on{Sets}$ given by the set of group homomorphisms $\bbG_m\to T_{{\overline{\bbQ}_p}}$. 
	Consider the subfunctor $X_*^E(-)\subseteq X_*(-)$ given by the subset of maps that are already defined over $E$. 
	This functor is co-representable by $\on{Res}_{E/\bbQ_p}\bbG_m$ and comes equipped with a universal cocharacter $\mu_{\on{univ}}\in X_*^E(\on{Res}_{E/\bbQ_p}\bbG_m)$. 
	In other words, given a torus $T\in  \{\on{Tori}_{\bbQ_p}\}$ and $\mu_0\in X^E_*(T)$ there is a unique map $\on{Nm}^E_{\mu_0}:\on{Res}_{E/\bbQ_p}\bbG_m\to T$ of algebraic groups over $\bbQ_p$ such that ${\on{Nm}}_{\mu_0}\circ \mu_{\on{univ}}=\mu_0$ in $X_*(T)$. 
	Representability of $X_*^E(-)$ can be obtained by recalling that the functor $X_*(-)$ defines an equivalence between the category $\{\on{Tori}_{\bbQ_p}\}$ and the category of free finite rank abelian groups endowed with $\Gamma_{\bbQ_p}$-action.
	The functor $X^E_*(-)$ corresponds to taking $\Gamma_{E}$-fixed points of $X_*(-)$ which is co-represented by 
	\[\on{Ind}^{\Gamma_E}_{\Gamma_{\bbQ_p}} \bbZ\simeq X_*(\on{Res}_{E/\bbQ_p}\bbG_m).\]
	The universal cocharacter is obtained from the unit of the adjunction 
	\[\bbZ\to \on{Res}_{\Gamma_E}^{\Gamma_{\bbQ_p}}\on{Ind}_{\Gamma_E}^{\Gamma_{\bbQ_p}}\bbZ\simeq \bbZ[\Gamma_{\bbQ_p}/\Gamma_E].\]
It can be expressed on $\bar{\bbQ}_p$-points as follows
\[\mu_{\on{univ}}:\bar{\bbQ}_p^\times \to (\bar{\bbQ}_p\otimes_{\bbQ_p} E)^\times\simeq  (\prod_{\iota:E\to \bar{\bbQ}_p} \bar{\bbQ}_p)^\times.\]
with $\mu_{\on{univ}}(e)=(e,1,\dots,1)$ where the first entry denotes the identity embedding $E\subseteq \bar{\bbQ}_p$ coming from the fact that we took $E$ to be a subfield of $\bar{\bbQ}_p$ and the other entries correspond to the various embeddings $\iota$ after ordering them. 
	We call the maps of tori obtained in this way \textit{norm} maps in analogy with the classical norm map form Galois theory
	\[\on{Nm}^E_{\bbQ_p}:E^\times \to \bbQ_p^\times\]
	with formula 
	\[e\mapsto \prod_{\iota:E\to \bar{\bbQ}_p} \iota(e),\]
	which can be regarded as the $\bbQ_p$-points of
	\[\on{Nm}^E_{\on{Id}_{\bbG_m}}:\on{Res}_{E/\bbQ_p}\bbG_m\to \bbG_m.\]

	Suppose now that $\mu\in\{\bbG_m\to G_{\bar{\bbQ}_p}\}/\sim$ is a conjugacy class of cocharacters with field of definition $E$. 
	Since the group $G^\ab$ is abelian, the conjugacy class $\mu^\ab\in\{\bbG_m\to G^\ab_{\bar{\bbQ}_p}\}/\sim$ defines a unique element $\mu^\ab\in X_*(G^\ab)$. Moreover, $\mu^\ab\in X^E_*(G^\ab)\subseteq X_*(G^\ab)$, and gives rise to a norm map 
	\[\on{Nm}^E_\mu:\on{Res}_{E/\bbQ_p}\bbG_m\to G^\ab.\]
	Recall the reciprocity character from local class field theory 
	\[\on{Art}_E:W_E\to E^\times.\]
	More precisely, denote by $E^{\on{ab}}$ the maximal abelian extension of $E$.
	Let $\on{rec}_E:E^\times \to \on{Gal}(E^{\on{ab}}/E)$ the local reciprocity map (norm residue symbol) of local class field theory.
	This map induces an isomorphism onto the image of the natural map $W_E\to \on{Gal}(E^\ab/E)$, and $\on{Art}_E$ is the unique map making the following diagram commute. 
	\begin{center}
	\begin{tikzcd}
	 & W_E\ar{rd} \arrow[ld, "\on{Art}_E", swap] & \\
		E^\times \arrow{rr}{\on{rec}_E} & & \on{Gal}(E^\ab/E)
	\end{tikzcd}
	\end{center}

\begin{proposition}
	\label{theactionofW}
	With the notation as above for all $\gamma\in W_E$ we have that $\rho_{W_E}(\gamma)=[\on{Nm}^E_\mu\circ \on{Art}_{E}(\gamma)]$.	
\end{proposition}
\begin{proof}
	\Cref{functoriality-wtih-group} and \Cref{thm:maintheorem} reduces us to show the statement in the case $G=G^\ab$.
	Since the map $\mu:\bbG_m\to G^\ab$ is defined over $E$ it factors as $\mu=\on{Nm}^E_\mu\circ \mu_{\on{univ}}$.
	By functoriality, we can reduce to the case $G=\on{Res}_{E/\bbQ_p}\bbG_m$ and $\mu=\mu_{\on{univ}}$.
	In this case, we should show that $\rho_{W_E}=\on{Art}_{E}$.
	The same argument as in \cite[Proposition 3.1.4]{PR24} proves this case. 
	The argument loc. cit. identifies 
	\[\Sht^{\on{geo}}_{(\on{Res}_{E/\bbQ_p}\bbG_m,\mu_{\on{univ}},b_{\on{univ}},\infty)}\] with the limit of a Rapoport--Zink tower of EL-type using \cite[Corollary 24.3.5]{SW20} and uses \cite{Che13} to describe the action explicitly. 

	Alternatively, we can argue as follows. 
	Recall that since $G=\on{Res}_{E/\bbQ_p}\bbG_m$ is a torus, $B(G,\mu_{\on{univ}})$ is a singleton since it is determined by the image of $\mu_{\on{univ}}$ in $\pi_1(G)_\Gamma=X_*(G)_\Gamma$ \cite[$\mathsection$ 7.6, (7.6.2)]{Kot97}.
	Representatives $b\in G(\breve{\bbQ}_p)=\breve{E}^\times$ of the element of $B(G,\mu_{\on{univ}})$ correspond to uniformizers in $\breve{E}$.
	Moreover, the groupoid of $G$-bundles over $X_{\on{FF},S}$ is equivalent to the groupoid of line bundles over $X_{\on{FF},E,S}$ the relative Fargues--Fontaine curve with respect to $E$ (instead of $\bbQ_p$) \cite[Definition II.1.15]{FS21}.
	Let us pick a uniformizer $\pi\in {E}$, this gives rise to a line bundle $\calO_\pi(1)$.
	We have a commutative diagram
	\begin{center}
	\begin{tikzcd}
		\Sht_{(\on{Res}_{E/\bbQ_p}\bbG_m,\mu_{\on{univ}},b_\pi,\infty)} \arrow{r} \arrow{d}{\simeq}  & \Gr^{\mu_{\on{univ}}}_{G,\Div^1_E} \arrow{d}{\simeq} \\
		\calB\calC(\calO_{E,\pi}(1))\setminus \{0\} \arrow{r}{m} & \Div^1_E 
	\end{tikzcd}
	\end{center}
	and an identification $\calB\calC(\calO_E(1))\setminus \{0\}/\underline{E^\times}\simeq \Div^1_E$ as in \cite[Corollary 2.4]{FS21}.
	It is proven in \cite[$\mathsection$ II.2.1]{FS21} that the $\underline{E^\times}$-torsor above is the one associated to rational Tate module of $G_{\on{LT}}$ the Lubin--Tate formal group of $E$. 
	This is the unique $1$-dimensional formal $O_E$-group over $O_{\breve{E}}$ for which the two $O_E$ actions on the Lie algebra coincide. 
	The relation of local class field theory with Lubin--Tate theory shows that $\rho_{W_E}=\on{Art}_{E}$ \cite{LT65}. 
\end{proof}

\subsection{The general case.}
\label{generalcase}
\subsubsection{z-extensions.}
In this section we assume that $G$ is an unramified reductive group over $\bbQ_p$ and that $(G,b,\mu)$ is HN-irreducible, but we no longer assume that $G^\der=G^\scn$.
We study this case using \Cref{thm:maintheorem}, \Cref{theactionofJ}, \Cref{theactionofW} and z-extensions techniques. 
Recall the following definition used extensively by Kottwitz \cite[$\mathsection$ 1]{RatConjug}.
\begin{definition}
	Let $f:G'\to G$ a map of connected reductive groups, let $Z=\on{Ker} f$. We say that $f$ is a z-extension if the following hold.
	\begin{enumerate}
		\item $f$ is surjective and $Z$ is central in $G'$.
		\item $Z$ is isomorphic to a product of tori of the form $\on{Res}_{F_i/\Q_p}\bbG_m$ for some finite extensions $F_i\subseteq \bar{\bbQ}_p$.
		\item $G'$ has simply connected derived subgroup. 
	\end{enumerate}
\end{definition}
By \cite[Lemma 1.1]{RatConjug} whenever $G$ is an unramified group over ${\bbQ_p}$ that splits over $\bbQ_{p^s}$, there exists a z-extension $G'\to G$ with $Z$ isomorphic to a product of tori of the form $\on{Res}_{\bbQ_{p^s}/{\bbQ_p}}\bbG_m$. 
In particular, $G'$ chosen in this way is unramified as well.
Recall from \cite[$\mathsection$ 6.5]{Kot97} that for any reductive group $G$ and cocharacter $\mu$ the natural morphism $B(G)\to B(G^{\on{ad}})$ induces a bijection $B(G,\mu)\simeq B(G^{\on{ad}},\mu^{\on{ad}})$.
From here one can deduce the following statement which we will use later on.
\begin{lemma}
	\label{lem:zextensionlemma}
	Fix $A\subseteq T\subseteq B\subseteq G$ as in \Cref{HN-irreducible}. 
	Assume that $\bbQ_{p^s}$ is a splitting field for $G$. 
	Let $\mu_0\in X_*^+(T)$, $b\in B(G,\mu)$, and fix $f:G'\to G$ a $z$-extension with $Z=\on{Ker}(f)$ isomorphic to $\prod_{i=1}^n\on{Res}_{\bbQ_{p^s}/{\bbQ_p}}\bbG_m$. Let $T'=f^{-1}(T)$ denote the maximal torus of $G'$ projecting onto $T$. Then the following hold.
	\begin{enumerate}
		\item For any choice of $\mu'\in X_*(T')^+$ lifting $\mu$ there is a unique lift $b'\in B(G')$ lifting $b$ with $b'\in B(G',\mu')$. 
		\item For $b'$ and $\mu'$ as in the previous claim $(b,\mu)$ is HN-irreducible if and only if $(b',\mu')$ is HN-irreducible.  
		\item If $E$ is the field of definition of $\mu$ with ${\bbQ_p}\subseteq E\subseteq \bbQ_{p^s}$ then there is a lift $\mu'\in X_*(T')^+$ with field of definition $E$.
	\end{enumerate}
\end{lemma}
\begin{proof}
	The first claim follows directly from the identifications $B(G,\mu)=B(G^{\on{ad}},\mu^{\on{ad}})=B(G',\mu')$.
	The second claim follows from the first claim, from the fact that $Z:=\on{Ker}(f)$ is central and from the fact that HN-irreducibility can be checked on the adjoint quotient (see proof of \cite[Corollary 4.1.16]{CKV15}). 
For the third claim consider the exact sequence of $\Gamma_{\bbQ_p}$-modules:
$$e\to X_*(Z)\to X_*(T')\to X_*(T)\to e$$
One can use Shapiro's lemma to prove $X_*(T')^{\Gamma_E}\to X_*(T)^{\Gamma_E}$ is surjective \cite[$\mathsection$ 2.3, \S 2.5 Proposition 10]{Ser02}. 
Indeed, 
\begin{align*}
	\on{H}^1(\Gamma_{E}, X_*(\on{Res}_{\bbQ_{p^s}/\bbQ_p}\bbG_m))& \simeq  \on{H}^1(\Gamma_{E},\bbZ[\Gamma_{\bbQ_p}/\Gamma_{\bbQ_{p^s}}])) \\ 
								     & \simeq \on{H}^1(\Gamma_{E},\bbZ[\Gamma_E/\Gamma_{\bbQ_{p^s}}])^{[E:\bbQ_p]} \\
								     & \simeq \on{H}^1(\Gamma_{\bbQ_{p^s}},\bbZ)^{[E:\bbQ_p]} \\
								     & \simeq \on{Hom}_{\on{cont}}(\Gamma_{\bbQ_{p^s}},\bbZ)^{[E:\bbQ_p]} \\
								     & \simeq 0
\end{align*}
\end{proof}

\subsubsection{Extended norm and determinant maps}
\label{extendedversionsofmpassection}
It is well-known that the image of the natural map $G^\scn(\bbQ_p)\to G(\bbQ_p)$, which we denote $\on{Im}G^\scn(\bbQ_p)\subseteq G(\bbQ_p)$, agrees with the commutator subgroup of $G(\bbQ_p)$ (see \Cref{drevidesubgroupsbusiness}).
\begin{remark}
	\label{drevidesubgroupsbusiness}
	For the convenience of the reader, we recall how this works. 
For a reductive group $G$, we let $G(\bbQ_p)^0\subseteq G(\bbQ_p)$ denote the subgroup generated by $\bbQ_p$-points contained in the unipotent radical of the $\bbQ_p$-rational parabolic subgroups of $G$.  
If $f:G^\scn\to G$ is the natural map then $f(G^{\scn}(\bbQ_p)^0)=G(\bbQ_p)^0$. 
It suffices to show that $G^{\scn}(\bbQ_p)^0=G^\scn(\bbQ_p)$ and that $G(\bbQ_p)^0$ is the commutator subgroup of $G(\bbQ_p)$.
The first statement follows from the affirmative resolution of the Kneser--Tits problem for non-archimedean fields (see \cite[$\mathsection$ 2]{PR85}). 
By \cite[Main Theorem]{Tit64} the commutator group of $G(\bbQ_p)$ contains $G(\bbQ_p)^0$ and by \cite[$\mathsection$ 1.4]{Tit64} this group is the commutator subgroup. 
\end{remark}

We let 
\[{G(\bbQ_p)_{\circ}}:=G(\bbQ_p)/\on{Im}G^\scn(\bbQ_p)\]
and denote by 
\[{\det}^\circ:G(\bbQ_p)\to {G(\bbQ_p)_{\circ}}\]
the natural quotient map. 
Whenever $G^\scn=G^\der$ we have that ${G(\bbQ_p)_{\circ}}=G^\ab(\bbQ_p)$ and the map $\det^\circ:G(\bbQ_p)\to G^\ab(\bbQ_p)$ is simply the evaluation on $\bbQ_p$-points of the algebraic map $\det:G\to G^\ab$ discussed above.
Nevertheless, when $G^\scn\neq G^\der$ the map $\det^\circ$ does not come directly from a map of algebraic groups.
We now explain a different perspective on $\det^\circ$.
Let $f:G'\to G$ be a z-extension such that $Z=\on{Ker}(f)$ is a product of tori of the form $\on{Res}_{F/\bbQ_p}\bbG_m$ with $F$ unramified over $\bbQ_p$.
By construction, $(G')^\der=G^\scn$.
This allows us to construct the following commutative diagram of topological groups.

\begin{center}
\begin{tikzcd}
	e\arrow{r} &G^\scn(\bbQ_p)\arrow{r}\arrow{d}  &G'(\bbQ_p)  \arrow{r}{\det} \arrow{d}  & (G')^\ab(\bbQ_p) \arrow{d} \arrow{r} & e \\
	e \arrow{r} & \on{Im}(G^\scn(\bbQ_p)) \arrow{r} & G(\bbQ_p) \arrow{r}{{\det}^\circ} & {G(\bbQ_p)_{\circ}} \arrow{r} & e 
\end{tikzcd}
\end{center}
Using Hilbert's 90 theorem and Shapiro's lemma one can show that $\on{H}_{\et}^1(\Spec \bbQ_p,Z)$ is trivial.
As a consequence, the maps $G'(\bbQ_p)\to G(\bbQ_p)$ and $(G')^\ab(\bbQ_p)\to {G(\bbQ_p)_{\circ}}$ are surjective.

To construct extended versions $\det^\circ_b$ of $\det_b$ we use z-extensions, so let us fix for the remainder of this subsection $f:G'\to G$ and $Z$ as above.
We define 
\[{\det}_b^\circ:G_b(\bbQ_p)\to {G(\bbQ_p)_{\circ}}\]
as follows.
We have a short exact sequence of the form 
\[e\to Z(\breve{\bbQ}_p)\to G'(\breve{\bbQ}_p)\to G(\breve{\bbQ}_p)\to e.\]
By abuse of notation, we let $b\in G(\breve{\bbQ}_p)$ denote a representative of $b\in B(G)$, we let $b'\in G'(\breve{\bbQ}_p)$ denote a lift of $b$. 
By definition, $G'_{b'}(R)=\{g\in G'(R\otimes_{\bbQ_p}\breve{\bbQ}_p)\mid g^{-1}b'\phi(g)=b'\}$ and analogously for $G_b$.  
Note that different lifts $b'$ of $b$ induce the same group $G'_{b'}$ since they differ by an element of the center of $G'$.
This induces a sequence of algebraic groups
\[e\to Z\to G'_{b'}\to G_b\to e.\]
In turn, this induces a sequence 
\[e\to Z(\bbQ_p)\to G'_{b'}(\bbQ_p)\to G_b(\bbQ_p)\to e\]
which is again exact since $\on{H}_{\et}^1(\Spec \bbQ_p,Z)$ is trivial.
We can consider the following commutative diagram
	\begin{center}
		\begin{tikzcd}
			Z(\bbQ_p)\ar{r}\ar{d}	& G'(\bbQ_p) \ar{r}\ar{d}{\det}& G(\bbQ_p)\ar{d}{{\det}^\circ}  \\
			G'_{b'}(\bbQ_p)	\ar{r}{\det_{b'}}\ar{dd}	& (G')^{\on{ab}}(\bbQ_p) \ar{r} & {G(\bbQ_p)_{\circ}}  \\
										& &  \\
			G_b(\bbQ_p)  \ar[bend right, swap, dotted]{rruu}{ \exists ! \on{det}_b^\circ}	& &   
		\end{tikzcd}
	\end{center}
	Since the image of $Z(\bbQ_p)$ in ${G(\bbQ_p)_{\circ}}$ is trivial, the surjectivity of the map $G'_{b'}(\bbQ_p)\to G_b(\bbQ_p)$ induces a unique map
	\[{\det}_b^\circ:G_b(\bbQ_p)\to {G(\bbQ_p)_{\circ}}.\]
	Finally, by \cite[Lemma 1.1.(3)]{RatConjug}, for every pair of z-extensions $G'_1\to G$ and $G'_2\to G$ we can find a third z-extension $G'_3$ and a commutative diagram 
	\begin{center}
	\begin{tikzcd}
		G'_3\arrow{r}{p_1} \arrow{d}{p_2}  & G'_1 \arrow{d} \\
	 G'_2\arrow{r} & G.
	\end{tikzcd}
	\end{center}
	It follows that the definition of ${\det}_b^\circ$ does not depend on the z-extension $G'\to G$ chosen above. 

	The extended versions of norm maps are easier to define.
	Fix the notation as in \Cref{defi:BGmu}. 
	Suppose we are given $\mu\in \{\bbG_m\to G\}/\sim$ with field of definition $E$ and let $\mu_0\in X^+_*(T)$ the unique dominant representative.
	Since both $B$ and $T$ are defined over $\bbQ_p$, it follows that $\mu_0$ is itself defined over $E$. 
	From our considerations in \Cref{subsubsectionactions} we obtain a commutative diagram of maps of algebraic groups
	\begin{center}
	\begin{tikzcd}
		\on{Res}_{E/\bbQ_p}\bbG_m \ar{rrd}{\Nm_\mu^E} \arrow{r}{\Nm_{\mu_0}^E}   & T  \ar{r} &  G \ar{d} \\
		& & G^\ab,
	\end{tikzcd}
	\end{center}
from which we can obtain a map of topological groups
\[\Nm^{E,\circ}_{\mu}:E^\times \to T(\bbQ_p)\to G(\bbQ_p)\to {G(\bbQ_p)_{\circ}}.\]
Since the group ${G(\bbQ_p)_{\circ}}$ is abelian this does not depend on the choice of rational Borel and torus $T\subseteq B\subseteq G$ fixed.

\subsubsection{The second main theorem.}
On this section $G$ is an unramified reductive group over $\bbQ_p$ and $(G,b,\mu)$ is HN-irreducible, but we no longer assume that $G^\der=G^\scn$.
We let ${G(\bbQ_p)_{\circ}}$, $\on{det}_b^\circ$, $\Nm^{E,\circ}_\mu$ and $\on{Art}_E$ be as above ($\mathsection$ \ref{extendedversionsofmpassection}, \ref{subsubsectionactions}). 
The following is our second main theorem.

\begin{theorem}
	\label{thm2:mainTheorem}
	Let $(G,b,\mu)$ be a $p$-adic shtuka datum such that $G$ an unramified reductive group over $\bbQ_p$ and such that the pair $(b,\mu)$ is HN-irreducible. 
	Let $E$ denote the reflex field of $\mu$.
	Then the following hold.
	\begin{enumerate}
		\item The right $G(\bbQ_p)$ action on $\pi_0(\Sht^{\on{geo}}_{G,b,\mu,\infty})$ is trivial on $\on{Im}(G^\scn(\bbQ_p))$ and the corresponding $G(\bbQ_p)_{\circ}$-action is simply-transitive. 
		\item If $s\in  \pi_0(\Sht^{\on{geo}}_{G,b,\mu,\infty})$ and $j\in G_b(\bbQ_p)$ then 
			\[s\star_{G_b} j=s\star_{{G(\bbQ_p)_{\circ}}} {\on{det}}_b^\circ(j^{-1}) \]		
		\item If $s\in  \pi_0(\Sht^{\on{geo}}_{G,b,\mu,\infty})$ and $\gamma \in W_{{E}}$ then 
			\[s\star_{{W_{E}}} \gamma =s\star_{{{G(\bbQ_p)_{\circ}}}} [{\on{Nm}}^{E,\circ}_{\mu}\circ {\on{Art}}_{{E}}(\gamma)].\] 
	\end{enumerate}
\end{theorem}
\begin{proof}
Let $f:G'\to G$ be a z-extension as in \Cref{lem:zextensionlemma}.	
Let $Z$ denote the kernel of $f$.
By \Cref{lem:zextensionlemma}, we may find a HN-irreducible triple $(G',b',\mu')$ with $\mu=f\circ \mu'$, $f(b')=b$ and such that $\mu'$ has the same reflex field as $\mu$. 
We claim that 
\[\Sht^{\on{geo}}_{G,b,\mu,\infty}=\Sht^{\on{geo}}_{G',b',\mu',\infty} \overset{\underline{G'(\bbQ_p)}}{\times} \underline{G(\bbQ_p)}.\]
Indeed, this is the content of \cite[Proposition 3.1.1, 3.1.2]{PR24}. 
Fix an element $x'\in \pi_0(\Sht^{\on{geo}}_{G',b',\mu',\infty})$ with image $x\in \pi_0(\Sht^{\on{geo}}_{G,b,\mu,\infty})$. 
It follows from \Cref{thm:maintheorem} that  
\[\pi_0(\Sht^{\on{geo}}_{G,b,\mu,\infty})\simeq \pi_0(\Sht^{\on{geo}}_{G',b',\mu',\infty} \overset{\underline{G'(\bbQ_p)}}{\times} \underline{G(\bbQ_p)})\simeq (G')^\ab(\bbQ_p)\overset{{G'(\bbQ_p)}}{\times} G(\bbQ_p)\]
as right $G(\bbQ_p)$-sets.
We can further write this as 
\[(G')^\ab(\bbQ_p)\overset{{G'(\bbQ_p)}}{\times} G(\bbQ_p)\simeq e\overset{{(G')^\der(\bbQ_p)}}{\times}G(\bbQ_p)\simeq e\overset{{G^{\on{sc}}(\bbQ_p)}}{\times}G(\bbQ_p)\simeq G(\bbQ_p)_\circ\]
This finishes the proof of the first statement.

Given $j\in G_b(\bbQ_p)$ we may find a lift to an element $j'\in G'_b(\bbQ_p)$ $\mathsection$ \ref{extendedversionsofmpassection}.
Then by \Cref{theactionofJ} and functoriality \Cref{functoriality-wtih-group} 
\[x\star_{G_b}j=f(x'\star_{G'_b} j')=f(x'\star_{(G')^\ab} \on{det}_{b'}(j'^{-1}))).\] 
Now, the element $\on{det}_{b'}(j'^{-1})\in (G')^\ab(\bbQ_p)\overset{{G'(\bbQ_p)}}{\times} G(\bbQ_p)$ is precisely $\on{det}_b^\circ(j^{-1})$ under the natural identification $(G')^\ab(\bbQ_p)\overset{{G'(\bbQ_p)}}{\times} G(\bbQ_p)\simeq {G(\bbQ_p)_{\circ}}$.
This finishes the proof of the second statement.
 
Similarly, given $\gamma\in W_E$ it follows from \Cref{theactionofW} and functoriality
\[x\star_{W_E}\gamma=f(x'\star_{W_E} \gamma)=f(x'\star_{(G')^\ab} \on{Nm}^E_{\mu'}\circ \on{Art}_E(\gamma)).\] 
So it suffices to show that the image of $\on{Nm}^E_{\mu'}\circ \on{Art}_E(\gamma)\in (G')^\ab(\bbQ_p)$ maps to $\on{Nm}^{E,\circ}_\mu\circ \on{Art}_E(\gamma)$ in ${G(\bbQ_p)_{\circ}}$.
Choose $T\subseteq B\subseteq G$ and let $T'=f^{-1}(T)$.
Consider the following diagram
\begin{center}
\begin{tikzcd}
	\Sht^{\on{geo}}_{G',b',\mu',\infty} \arrow{r} \arrow{d}  & \Sht^{\on{geo}}_{((G')^\ab,b'^\ab,\mu'^\ab,\infty)} & \ar{l} \Sht^{\on{geo}}_{(T',b'_{T'},\mu'_{T'},\infty)} \arrow{d} \\
	\Sht^{\on{geo}}_{G',b',\mu',\infty} & &  \Sht^{\on{geo}}_{(T,b_{T},\mu_{T},\infty)}.
\end{tikzcd}
\end{center}
Here we have taken the dominant representatives $\mu_T$ and $\mu'_T$ of $\mu$ and $\mu'$ in $X_*^E(T)$ and $X^E_*(T')$ respectively.
After passing to connected components and fixing elements $x'\in \pi_0(\Sht^{\on{geo}}_{G',b',\mu',\infty})$ and $t'\in \pi_0(\Sht^{\on{geo}}_{(T',b'_{T'},\mu'_{T'},\infty)})$ with a common image in $\pi_0(\Sht^{\on{geo}}_{((G')^\ab,b'^\ab,\mu'^\ab,\infty)})$ the diagram becomes
\begin{center}
\begin{tikzcd}
	(G')^\ab(\bbQ_p)	\arrow{r}{\simeq} \arrow{d} & (G')^\ab(\bbQ_p)  & T'(\bbQ_p) \arrow{l} \arrow{d} \\
	{G(\bbQ_p)_{\circ}}	 & & T(\bbQ_p). 
\end{tikzcd}
\end{center}
We see that $\on{Nm}^E_{\mu'}\circ \on{Art}_E(\gamma)$ is the image of $\on{Nm}^E_{\mu'_{T'}}\circ \on{Art}_E(\gamma)$ under $T'(\bbQ_p)\to (G')^\ab(\bbQ_p)$.  
But the map $T'(\bbQ_p)\to {G(\bbQ_p)_{\circ}}$ factors through $T(\bbQ_p)$. 
The image of $\on{Nm}^E_{\mu'_{T'}}\circ \on{Art}_E(\gamma)$ in $T(\bbQ_p)$ is $\on{Nm}^E_{\mu_{T}}\circ \on{Art}_E(\gamma)$ and the image of this element in ${G(\bbQ_p)_{\circ}}$ is by definition $\on{Nm}^{E,\circ}_\mu\circ \on{Art}_E(\gamma)$ as we wanted to show.
\end{proof}

\begin{remark}
	A formulation of \Cref{thm2:mainTheorem} appears in the work of Chen as \cite[Th\'eor\`{e}me 7.0.2]{Chen}.
	Chen works directly with Rapoport--Zink spaces as rigid-analytic spaces and at the time the theory of diamonds had not been developed. In particular, it was not completely established how to consider spaces at infinite level. 
	Nevertheless, it still made sense to take the colimit of the top degree cohomology groups with compact support, which by Poincar\'e duality captures the behavior of connected components.
	In rough terms, the result of Chen is related to \Cref{thm2:mainTheorem} by appealing to \cite[Corollary 24.3.5]{SW20} and by passing from connected components to cohomology.
\end{remark}

\begin{remark}
	\label{whereHn-comesin}
	\Cref{thm2:mainTheorem} is optimal for unramified groups in the following sense. 
	One can prove that the action of $G(\bbQ_p)$ on $\pi_0(\Sht^{\on{geo}}_{G,b,\mu,\infty})$ only factors through ${G(\bbQ_p)_{\circ}}$ when $(b,\mu)$ is HN-irreducible. 
	Moreover, we expect that combining the methods of \cite{gaisin2022nonsemistablelociheckestacks} and \cite{MR4331441} with the methods in the present article one can express the general formula for the $G(\bbQ_p)\times G_b(\bbQ_p)\times W_E$ action on $\pi_0(\Sht^{\on{geo}}_{G,b,\mu,\infty})$ in terms of parabolic induction of $\pi_0(\Sht^{\on{geo}}_{M,b_M,\mu_M,\infty})$ for HN-irreducible data $(M,b_M,\mu_M)$ associated to Levi subgroups $M\subseteq G$ appearing in the Hodge-Newton decomposition of $(b,\mu)$.
\end{remark}

\begin{remark}
	\label{Remark-unramifiedness-rmove}
	During the revision process of this article, in a joint work with Lim and Xu \cite{gleason2023connectedcomponentsaffinedelignelusztig}, we found a method to generalize \Cref{thm:maintheorem} to groups $G$ that are no longer assumed to be unramified. On that collaboration, we build on the methods developed on this article to show the general case. As a corollary, we showed almost all cases of \Cref{conj-rapoport-viehmann} excluding only cases where the group has anisotropic semisimple factors.
\end{remark}

	\bibliography{biblio.bib}
	\bibliographystyle{alpha}
	
\end{document}